\documentclass[preprint]{elsarticle}

\usepackage{amssymb}
\usepackage{indentfirst}
\usepackage{graphicx}
\usepackage{epstopdf}
\usepackage{epsfig}
\usepackage{overpic}
\usepackage{array}
\usepackage{color}
\usepackage{multirow}
\usepackage{hyperref}
\usepackage{float}
\usepackage{subfigure}
\usepackage{amssymb}
\usepackage{amsmath,amssymb, amsthm}
\usepackage[mathscr]{euscript}
\usepackage{latexsym,bm}
\usepackage{booktabs}
\usepackage{bbm}
\usepackage{diagbox}
 \usepackage[all]{xy}
\usepackage{diagbox}
\usepackage{algorithm}  
\usepackage{algorithmicx} 
\usepackage{algpseudocode} 
\usepackage{mathtools}

\newtheorem{theorem}{Theorem}
\newtheorem{lemma}{Lemma}
\newtheorem{proposition}{Proposition}
\newtheorem{remark}{Remark}
\newtheorem{example}{\textup{\textbf{Example}}}

\newcommand\cf{\color{red}}

\newcommand\D{\textup{d}}

\def\me{\mathrm{e}}

\def\bx{\mathbf{x}}

\makeatletter
\@addtoreset{equation}{section}
\makeatother

\begin{document}

\newsavebox{\tablebox}

\begin{frontmatter}



\title{ Efficient structure-preserving scheme for chemotaxis PDE system with singular sensitivity in crime and epidemic modeling}

\author[label]{Rui Wang}
\ead{rwang0913@mail.bnu.edu.cn}

\author[label,label2]{Yunfeng Xiong\corref{cor}}
\ead{yfxiong@bnu.edu.cn}

\author[label,label2]{Zhengru Zhang}
\ead{zrzhang@bnu.edu.cn}

\affiliation[label]{organization={School of Mathematical Sciences, Beijing Normal University},
          postcode={100875}, 
            state={Beijing},
            country={P. R. China}
 }
 \affiliation[label2]{organization={Laboratory of Mathematics and Complex Systems, Ministry of Education, Beijing Normal University, Beijing},
            postcode={100875},
            state={Beijing},
            country={China}}
 
\cortext[cor]{To whom correspondence should be addressed}

\begin{abstract}
The chemotaxis PDE system with singular sensitivity was originally proposed by Short et al. [Math. Mod.
Meth. Appl. Sci., 18:1249–1267, 2008] as the continuum limit of a biased random walk model to account for the formation of crime hotspots and environmental feedback successfully. Recently, this idea has also been applied to epidemiology to model the impact of human social behaviors on disease transmission. In order to characterize the phase transition, pattern formation and statistical properties in the long-term dynamics, a stable and accurate numerical scheme is urgently demanded, which still remains challenging due to the positivity constraint on the singular sensitivity and the absence of an energy functional. In particular, the loss of positivity may produce nonphysical states and even cause spurious blow-up.
To address these numerical challenges, this paper constructs an efficient positivity-preserving, implicit-explicit scheme with second-order accuracy. A rigorous error estimation is provided with the Lagrange multiplier correction to deal with the singular sensitivity. The whole framework is extended to a multi-agent epidemic model with degenerate diffusion, in which both positivity and mass conservation are achieved. 
Numerical experiments are performed to validate the theoretical results and demonstrate the necessity 
of the correction strategy. Our simulations reveal rich dynamical behaviors, including the phase transition between 
aggregation-dominated and dissipative regimes, as well as the nucleation, spread, and dissipation of crime hotspots. 
For the epidemic model, the results further show that spatial clustering of population density may accelerate 
virus transmission and significantly amplify the infectious wave.

\end{abstract}






\end{frontmatter}

\section{Introduction}

We investigate the chemotaxis PDE system with singular sensitivity:
\begin{equation}\label{cross_diffusion_PDE}
	\left\{
	\begin{aligned}
		&\frac{\partial \phi(\bx, t)}{\partial t}   =   \mathcal{R} \phi(\bx, t) - (p(\bx, t)+p^0(\bx))\phi(\bx,t) + \gamma(\bx), \quad &\mathbf{x}\in \Omega, ~t>0, \\
		&\frac{\partial p(\bx,   t)}{\partial t} = \eta D \Delta p(\bx,   t) + (p(\bx, t)+p^0(\bx))\phi(\bx,t) - p(\bx,t) ,\quad & \mathbf{x}\in \Omega, ~t>0,\\
		& \phi(\mathbf{x},0) = \phi_0(\mathbf{x})\geq 0,\quad p(\mathbf{x},0) = p_0(\mathbf{x})\geq 0,\quad &\mathbf{x}\in\Omega,
	\end{aligned}
	\right.
\end{equation}
where $p(\bx, t)$ denotes the spatially varying component of the attractiveness field.  The vector-valued fields $\phi(\bx, t)$ are densities of multiple agents driven by the chemotaxis-type operator $\mathcal R$, 
\begin{equation}\label{chemotaxis_operator}
	\mathcal{R}  \phi(\bx,   t ) : = \dfrac{{D}}{4}\nabla\cdot \left ( \nabla \phi(\bx,   t ) - \dfrac{2  \phi(\bx,   t )}{   p(\bx,  t) + p^0(\bx)} \nabla   (p(\bx,  t) +  p^0(\bx))\right),
\end{equation} 
so that agents move up gradients of $\log (p + p^0)$ but simply diffuse in the absence of a risk gradient \cite{ShortBrantinghamBertozziTita2010}.
Here $p^0(\mathbf{x})$ is a fixed,  smooth and small environmental field to prevent the blow-up of $\ln p(\mathbf{x},t)$ when $p(\mathbf{x},t)$ approaches zero, which  satisfies
$0 < m_p \le p^0(\mathbf{x}) \le M_p
$ for some constants $m_p,M_p$. The source term $\gamma(\bx)$ alters the equilibrium states. The coefficients $D$ and $\eta$ quantify the diffusion rates of agents and the attraction field, respectively. It naturally requires to impose the positivity constraints $\phi(\bx, t) \ge 0, p(\bx, t) \geq 0$ \cite{LiXie2024}.

The above chemotaxis model, as the continuum limit of a biased random walk model, was originally proposed in pioneering works \cite{ShortDorsognaPasourTitaBrantinghamBertozziChayes2008,ShortBrantinghamBertozziTita2010} to explain the formation of crime hotspots and environmental feedback, where its chemotaxis-type operator $\mathcal R$ describes the complex competition between crowd aggregation and diffusion. It deeply reveals the mechanism of generation and evolution of crime hotspots by simultaneously characterizing the bidirectional coupling of individual (or behavioral) diffusion and environmental feedback at a multi-scale level \cite{BellomoOutadaSolerTaoWinkler2022,CaiWangZhang2022}. This idea has recently been applied to epidemiology to investigate how human social behaviors, like crowding and public awareness, influence disease transmission. Studies have shown that the chemotaxis mechanism may produce long-term stable spatial clustering \cite{CSIAM-LS-1-2} and significantly influence spatiotemporal epidemic waves \cite{XiongWangZhang2024}. Moreover, it is a promising tool for other practical applications, such as opinion dynamics, due to its strong explanatory power and deep connection to statistical mechanics.

Despite the profound successes, both the theoretical and numerical aspects of the chemotaxis PDE system \eqref{cross_diffusion_PDE} face significant challenges. The existence of global-in-time classical solution has been established only for a specified range of $\eta$ \cite{Winkler2010,TaoWinkler2021}, due to the absence of an energy functional. Only recently has the dependence on $\eta$ been lifted for global generalized solutions in the 2-D setting \cite{LiXie2024}. The singular chemotactic sensitivity also poses significant challenges for numerical computation, as preserving the positivity of the solution is essential to ensure its well-posedness. This requirement is particularly critical in epidemic models, as negative solution values may induce numerical instabilities and severely hinder long-term simulations \cite{blanes2022positivity}. Short et al. first proposed a first-order semi-implicit scheme to study the emergence, evolution, and steady-state characteristics of crime hotspots \cite{ShortBrantinghamBertozziTita2010}. For partially
heterogeneous model and realistic urban geometries, a finite element framework with a first-order backward differentiation formula for time discretization has also been developed to achieve an efficient and robust numerical resolution \cite{HaoMilyQuainiZhong2026}.  Regarding the accuracy of long-time simulations, however, a higher-order positivity-preserving scheme is more desirable.
Recently, a delicate structure-preserving scheme has been proposed for susceptible-infected-susceptible models based on the time splitting strategy to achieve second-order accuracy, and a linear stability analysis has been conducted to study the impact of chemical sensitivity on the equilibrium \cite{ding2025structure}. The cost to pay is the introduction of a stabilizer, making the scheme fully implicit and thus demanding highly efficient iterative solvers.

It is worth noting that the chemotaxis PDE system \eqref{cross_diffusion_PDE} can be viewed as a variant of the classical Keller-Segel (KS) model \cite{KELLER1970399, Patlak1953,CAO2016382, arumugam2021keller}, so that the same kind of numerical difficulties is shared.
Regarding the stable evolution of the classical KS model, Shen and Xu constructed a class of schemes that preserves mass conservation, uniqueness, positivity, and energy dissipation \cite{Shenjie2020}. The first variational structural scheme that simultaneously achieves second-order accuracy, positivity preservation, and original energy dissipation has been proposed, and its uniqueness, convergence, and robustness have been verified \cite{Ding2025SecondOrder}. For spatial numerical resolution, a fourth-order finite difference discretization with positivity preserving and energy dissipation was designed to achieve stable long-time simulation \cite{hu2023positivity}. The conservative upwind finite volume method and central-upwind positivity-preserving scheme were also applied to the KS model \cite{filbet2006finite,chertock2008second}. A simplified linear finite volume method that satisfies positivity and mass conservation has been designed, and a discrete free energy inequality and $L^p$ error estimation are given \cite{zhou2017finite}. By applying a log‑transformation to preserve positivity and incorporating a recovery strategy to ensure mass conservation, a linear and decoupled finite element method that attained optimal error bounds was given in \cite{wang2024optimal}. However, the numerical analysis of classical KS models, which relies heavily on the energy functional, might not be readily generalized to chemotaxis models lacking an energy structure. Recently, another structure-preserving framework, termed the Lagrange multiplier method \cite{van2019positivity}, was introduced to efficiently handle the constraints of positivity, mass conservation, boundary conservation, and energy dissipation \cite{cheng2022new,cheng2022bound,cheng2025new,TongFenghua2024Positivity}. The basic idea is to force the numerical solution to satisfy physical constraints by projection (or correction), without reliance on the energy structure. Moreover, it can be coupled with explicit integrators to boost efficiency. 

This paper designs a structure-preserving implicit-explicit (SPIMEX) scheme for Eq.~\eqref{cross_diffusion_PDE} to achieve second-order convergence, structure-preserving property, and efficiency simultaneously, paving the way for accurate description of pattern formulation and statistical property for singular chemotaxis movements. SPIMEX adopts a predictor-corrector framework. First, an intermediate solution is computed using the finite‑difference method (FDM), with the linear diffusion term implicitly discretized by the Crank–Nicolson scheme to eliminate stiffness, and the nonlinear terms explicitly advanced in time to improve computational efficiency.
Second, the intermediate solution is projected onto the constraint manifold, which has been formulated as a convex $L^2$-$H^1$ minimization problem. For numerical analysis, the key to address the difficulty arising from the singular chemotactic sensitivity is to utilize a post-processing Lagrange multiplier method to provide a posteriori error estimation \cite{BeckerRannacher2001}. For the crime model \cite{ShortDorsognaPasourTitaBrantinghamBertozziChayes2008,ShortBrantinghamBertozziTita2010}, a rigorous error bound is established, guaranteeing second‑order convergence in both space and time under the $L^2$-$H^1$ norm. This framework can be readily extended to a multi‑agent epidemic model with degenerate diffusion in the hospitalized agents \cite{XiongWangZhang2024}. Numerical experiments on the crime model provide a reliable basis for understanding the spatial dynamics of crime hotspots under various influencing factors. Simulations of the epidemic model reveal that local population density plays a key role in transmission dynamics: High‑density clusters significantly enhance the peak intensity of the second wave, quantitatively confirming that crowding aggravates epidemic fluctuations by accelerating virus transmission.

The rest of this paper is organized as follows. Notations are introduced in \ref{sec:notation}. The setting of SPIMEX with the $L^2$-$H^1$ projection strategy for Eq.~\eqref{cross_diffusion_PDE} and a rigorous numerical analysis are provided in \ref{sec:cross}. In \ref{sec:epi}, the scheme is generalized to a multi‑agent epidemic model to maintain positivity‑preserving and mass conservation properties, and refine the error estimation. In \ref{sec:numerical}, the theoretical convergence order of the proposed schemes is verified through typical numerical experiments, and phase transition, pattern formation and statistical characteristics of both crime and epidemic models are demonstrated. Finally, conclusions are drawn in \ref{sec:conclusion}.

\section{Notations}\label{sec:notation}

Suppose the computational domain $\Omega = [a_{\min}, a_{\max}] \times [b_{\min}, b_{\max}]$ is uniformly partitioned into a grid
\[
\Omega_h := \left\{ (x_i, y_j) \mid x_i = a_{\min} + i h_x,\ y_j = b_{\min} + j h_y,\ 0 \le i \le N_x,\ 0 \le j \le N_y \right\},
\]
with spacings
$h_x = \frac{a_{\max} - a_{\min}}{N_x}$, $h_y = \frac{b_{\max} - b_{\min}}{N_y}$, where $N_x$ and $N_y$ are numbers of grid points.

We consider discrete grid functions $u_h: \Omega_h \to \mathbb{R}$, which can be represented by a set of values $u_{i,j} = u_h(x_i, y_j)$ for $(x_i, y_j) \in \Omega_h$. Under periodic boundary conditions, the discrete function space for $u_h$ is defined as
\[
X := \left\{ u_h: \Omega_h \to \mathbb{R} \mid u_{0,j} = u_{N_x,j},\ u_{i,0} = u_{i,N_y},\ \forall 0 \le i \le N_x,\ 0 \le j \le N_y \right\}.
\]
For simplicity, we assume that the domain is square, i.e., $a_{\max} - a_{\min} = b_{\max} - b_{\min}$, and that the grid is uniform in both directions, so that $h_x = h_y = h$ and $N_x = N_y = N$.

Let $u(\mathbf{x}, t)$ be the exact solution (scalar- or vector-valued) on $\Omega \times [0, T]$, and $u_h(t_k)$ denotes the exact solution on the grid mesh at time $t_k$. The predicted and corrected numerical solutions are denoted by $\tilde{u}_h^k$ and $u_h^k \in X$, respectively. For vector-valued fields, we use bold notations. 

For $u_h\in  X $, we denote the averaging operator $\mathcal{A}_x(\mathcal{A}_y):  X \rightarrow  X $ as
\begin{equation}\nonumber
	\mathcal{A}_xu_{i+\frac{1}{2},j}:=\frac{u_{i+1,j}+u_{i,j}}{2},\quad \mathcal{A}_yu_{i,j+\frac{1}{2}}:=\frac{u_{i,j+1}+u_{i,j}}{2},
\end{equation}
and the  difference operators are introduced on the function spaces
\begin{equation}\nonumber
	D_xu_{i+\frac{1}{2},j}:=\frac{1}{h}(u_{i+1,j}-u_{i,j}),\quad D_yu_{i,j+\frac{1}{2}}:=\frac{1}{h}(u_{i,j+1}-u_{i,j}).
\end{equation}
Likewise
\begin{equation}\nonumber
	d_xu_{i,j}:=\frac{1}{h}(u_{i+\frac{1}{2},j}-u_{i-\frac{1}{2},j}),\quad d_yu_{i,j}:=\frac{1}{h}(u_{i,j+\frac{1}{2}}-u_{i,j-\frac{1}{2}}), 
\end{equation}
with $d_x,d_y: X \rightarrow  X $. The discrete gradient $\nabla _h:  X \rightarrow  X $ is defined by
\begin{equation}\nonumber
	\nabla_h u_{i,j}:= (D_xu_{i+\frac{1}{2},j}, D_yu_{i,j+\frac{1}{2}}),
\end{equation}
and the discrete divergence $\nabla_h\cdot:  X ^2\rightarrow  X $ reads that
\begin{equation}\nonumber
	\nabla_h\cdot \vec{f}:=d_xf^x_{i,j}+d_y f^y_{i,j},\quad \text{for~} \vec{f}=(f^x,f^y)\in  X \times X.
\end{equation}
The standard discrete Laplacian, $\Delta _h: X \rightarrow  X $, is given  by
\begin{equation}\nonumber
	\Delta _h u_{i,j} := \nabla_h\cdot (\nabla_h u)_{i,j}=\frac{1}{h^2}(u_{i+1,j}+u_{i-1,j}+u_{i,j+1}+u_{i,j-1}-4u_{i,j}).
\end{equation}
More generally, if $\mathcal{D}$ is a periodic scalar function that is defined at all of the face center points and $\vec{f}\in  X ^2$, by assuming point-wise multiplication, we may define
\begin{equation}\nonumber
	\nabla_h\cdot (\mathcal{D}\vec{f})_{i,j}:= d_x(\mathcal{D}f^x)_{i,j}+d_y(\mathcal{D}f^y)_{i,j}.
\end{equation}
Specifically, if $u\in  X $, then $\nabla_h\cdot (\mathcal{D}\nabla_h): X \rightarrow  X $ is defined pointwisely via
\begin{equation}\nonumber
	\nabla_h\cdot(\mathcal{D}\nabla_h u)_{i,j}:=d_x(\mathcal{D}D_xu)_{i,j}+d_y(\mathcal{D}D_yu)_{i,j}.
\end{equation}

For $u_h,v_h\in  X $, the discrete $L^2$ inner product, the $L^2$ norm and the $L^p(1\leq p\leq\infty)$ norm are given by
\begin{equation}\nonumber
	\langle u_h,v_h\rangle = h^2\sum_{i,j=1}^{N}u_{i,j}v_{i,j},\quad \Vert u_h\Vert_{L^2}^2=\langle u_h,u_h\rangle,\quad \Vert u_h\Vert_{L^p}^p=h^2\sum_{i,j=1}^{N}\vert u_{i,j}\vert^p,
\end{equation}
and the $L^\infty$ norm is defined as $\Vert u\Vert_\infty =\max_{1\leq i,j\leq N}\vert u_{i,j}\vert$.  The discrete $H^1$ norm for $u_h\in  X $ is
\begin{equation}\nonumber
	\Vert \nabla_h u_h\Vert_{L^2}^2=\langle \nabla_h u_h,\nabla_h u_h\rangle,\quad \Vert u_h\Vert_{H^1}^2 = \Vert u_h\Vert_{L^2}^2 + \Vert \nabla u_h\Vert_{L^2}^2.
\end{equation}

Let the vector-valued grid function
$\mathbf{u}_h=(u_{1,h},\dots,u_{m,h})^\top,\mathbf{v}_h=(v_{1,h},\dots,v_{m,h})^\top\in X^m$, $ u_{k,h}=\{u_{k,i,j}\}_{1\le i,j\le N},v_{k,h}=\{v_{k,i,j}\}_{1\le i,j\le N}$, the discrete vector inner product and the $L^2$ norm are then defined by
\begin{equation}\nonumber
	\langle \mathbf{u}_h,\mathbf{v}_h\rangle
	= h^2\sum_{i,j=1}^{N}\sum_{k=1}^m u_{k,i,j}v_{k,i,j},\quad \|\mathbf{u}_h\|_{L^2}^2=\langle \mathbf{u}_h,\mathbf{u}_h\rangle_h
	= h^2\sum_{i,j=1}^{N}\sum_{k=1}^m \vert u_{k,i,j}\vert^2.
\end{equation}
Furthermore, 
\begin{equation}\label{eq:component_sum}
	\|\mathbf{u}_h\|_{L^2}^2=\sum_{k=1}^m\|u_{k,h}\|_{L^2}^2,\quad
	\|\mathbf{u}_h\|_{L^2}=\Big(\sum_{k=1}^m\|u_{k,h}\|_{L^2}^2\Big)^{1/2}
	\le \sum_{k=1}^m\|u_{k,h}\|_{L^2}.
\end{equation}

For any grid functions $u_h,v_h,w_h\in  X $, the summation by parts holds in the discrete sense as follows
\begin{equation}\label{summation_by_parts}
	\langle\nabla_h u_h,\nabla_h v_h\rangle = -\langle \Delta_h u_h, v_h\rangle, \quad
	\langle \nabla_h\cdot (\mathcal{A}_h w_h)\nabla_h u_h,v_h \rangle =-\langle (\mathcal{A}_hw_h)\nabla_h u_h,\nabla_h v_h\rangle.
\end{equation}

\section{SPIMEX  for  chemotaxis modeling}\label{sec:cross}

In this section, we introduce SPIMEX for solving the PDE system \eqref{cross_diffusion_PDE} and provide its error analysis. For the sake of presentation, we focus on the 2D case, i.e. a rectangular box $\Omega = [a_{\min},a_{\max}]\times [b_{\min},b_{\max}]$. The numerical scheme can be easily extended to 1D and 3D cases due to tensor construction and theoretical results remain the same. Without loss of generality, we assume $\gamma=0$ in the subsequent analysis.

\subsection{FDMs with \texorpdfstring{$L^2$-$H^1$}{L2-H1} projection}

Now, we are ready to construct our FDMs based on the prediction-correction strategy. We choose the size of the time steps $\tau>0$, and the time steps are $t_k:=k\tau$ for $k=0,1,\dots$. Let $(\phi_h^k,p_h^k)\in X\times X$ be the numerical approximation of the exact solution $(\phi,p)$ of \eqref{cross_diffusion_PDE} on $\Omega_h$ at time $t_k$.

For $k = 0,1,2,\dots$, we employ a prediction-correction strategy.
\begin{itemize}
	\item[Step 1]  Compute the intermediate solutions $\tilde{\phi}_h^{k+1},\tilde{p}_h^{k+1}$ using the Crank–Nicolson finite difference scheme (CNFD) coupled with an Adams–Bashforth treatment of the nonlinear term:
	\begin{equation}\label{CNFD}
		\begin{split}
			\frac{\tilde{\phi}_h^{k+1}-\phi_h^k}{\tau} = &\frac{D}{8}\Delta_h(\tilde{\phi}_h^{k+1}+\phi_h^k) - \frac{3D}{4}\nabla_h\cdot \left(\frac{\phi_h^k}{p_h^k+p^0}\nabla_h (p_h^k+p^0) \right) \\
			& + \frac{D}{4}\nabla_h\cdot \left(\frac{\phi_h^{k-1}}{p_h^{k-1}+p^0}\nabla_h (p_h^{k-1}+p^0)\right) -\frac{3}{2}(p_h^k+p^0)\phi_h^k+\frac{1}{2}(p_h^{k-1}+p^0)\phi_h^{k-1},\\
			\frac{\tilde{p}_h^{k+1}-p_h^k}{\tau} = &\frac{\eta D}{2}\Delta_h(\tilde{p}_h^{k+1}+p_h^k) +\frac{3}{2}[(p_h^k+p^0)\phi_h^k-p_h^k]  -\frac{1}{2}[(p_h^{k-1}+p^0)\phi_h^{k-1}-p_h^{k-1}].
		\end{split}
	\end{equation}
	
	\item[Step 2]  Correct the  intermediate solutions. In general, $\tilde{\phi}_h^{k+1}$ and $\tilde{p}_h^{k+1}$ computed by the linear semi-implicit schemes may break the positivity preservation \cite{li2020exponential}.  Here, we treat the positivity at the discrete level as the constraints for the numerical solution at $t^{k+1}$ as $( \phi_h^{k+1},p_h^{k+1})$ with $\phi_h^{k+1}, p_h^{k+1}\geq 0$. A natural way to obtain such $( \phi_h^{k+1},p_h^{k+1})$ from $(\tilde{\phi}_h^{k+1}, \tilde{p}_h^{k+1})$ is to project the nodal vectors $\tilde{\phi}_h^{k+1}$ and $\tilde{p}_h^{k+1}$ to the constrained manifold with positivity preservation. We adopt the $L^2$-$H^1$ projection here to enforce the positivity, which reads
	\begin{equation}\label{projection}
		\min_{\phi_h^{k+1},p_h^{k+1}\in  X } \frac{1}{2}(\Vert \phi_h^{k+1} - \tilde{\phi}_h^{k+1}\Vert_{L^2}^2+\Vert p_h^{k+1} - \tilde{p}_h^{k+1}\Vert_{H^1}^2),\quad
		\text{s.t.}~~ \phi_h^{k+1}\geq 0, \quad p_h^{k+1}\geq 0 . 
	\end{equation}
	This is a convex minimization problem with Karush-Kuhn-Tucker (KKT) conditions:
	\begin{equation}\label{KKT_1}
		\begin{aligned}
			&\phi_h^{k+1} =  \tilde{\phi}_h^{k+1}+\lambda_h^{k+1}, \quad (I - \Delta_h) p_h^{k+1} =  (I - \Delta_h) \tilde{p}_h^{k+1}+\xi_h^{k+1},	\\
			&\lambda_h^{k+1}\phi_h^{k+1} = 0,\quad \xi_h^{k+1} p_h^{k+1} = 0,\quad \lambda_h^{k+1}\geq 0,\quad \xi_h^{k+1}\geq 0,
		\end{aligned}
	\end{equation}
	where $\lambda_h^{k+1}, \xi_h^{k+1}\in  X $  are the Lagrange multipliers for the positivity preservation, and $I$ is the identity operator. 
\end{itemize}

Now, \eqref{CNFD}-\eqref{KKT_1} complete SPIMEX scheme.  Since
\eqref{CNFD} is linear and the convex minimization problem \eqref{projection} admits unique solutions, we find SPIMEX is uniquely solvable at each time step.
\begin{remark}
	Since \eqref{CNFD} is a three-level scheme, for the first step $k = 0$, we use the first order scheme instead
	\begin{equation}\label{first scheme}
		\begin{split}
			\frac{\tilde{\phi}_h^{1}-\phi_h^0}{\tau}&=\frac{D}{8}\Delta_h(\tilde{\phi}_h^{1}+\phi_h^0)-\frac{D}{2}\nabla_h\cdot\left(\frac{\phi_h^{0}}{p_h^{0}+p^0}\nabla_h(p_h^{0}+p^0)\right) -(p_h^0+p^0)\phi_h^0,\\
			\frac{\tilde{p}_h^{1}-p_h^0}{\tau}&=\frac{\eta D}{2}\Delta_h\left(\tilde{p}_h^{1}+p_h^0\right) +(p_h^0+p^0)\phi_h^0-p_h^0.
		\end{split}
	\end{equation}
\end{remark}

\begin{remark}
	It is not necessary to explicitly compute the value of $\lambda_h^{n+1}$ and $\xi_h^{n+1}$, since we can use the complementary slackness property in KKT conditions \eqref{KKT_1} to determine the solution.
\end{remark}

\subsection{Error estimation}
Now, we carry out the error analysis for \eqref{CNFD}-\eqref{KKT_1} with \eqref{first scheme}. Let $T>0$ be a fixed time, and $(\phi(\bx,t)\geq0, p(\bx,t)\geq0)$ be the exact solution of \eqref{cross_diffusion_PDE}. Based on the theoretical results, we make the following assumptions,
\begin{equation}\label{assumption}
	\phi(\bx,t)\in C^3([0,T];C_{per}^4(\Omega)),\quad  p(\bx,t)\in C^3([0,T];C_{per}^5(\Omega)),
\end{equation}
where $C_{per}^m(\Omega) = \{u\in C^m(\Omega)|\partial _x^k\partial_y^l u ~  \text{is periodic on}~ \Omega,\forall k,l\geq 0,k+l\leq m\}$. 

We introduce the biased error functions $e_\phi^k,e_p^k,\tilde{e}_\phi^k,\tilde{e}_p^k\in  X  (k\geq 0)$:
\begin{equation}\label{error_denfinition}
	\begin{aligned}
		e_\phi^k=\phi_h(t_k)-\phi_h^k,\quad e_p^k=p_h(t_k)-p_h^k,\\
		\tilde{e}_\phi^k=\phi_h(t_k)-\tilde{\phi}_h^k,\quad \tilde{e}_p^k=p_h(t_k)-\tilde{p}_h^k,
	\end{aligned}
\end{equation}
where $\tilde{\phi}_h^0=\phi_h^0$, $\tilde{p}_h^0=p_h^0$.
The following error bounds can be established.
\begin{theorem}\label{main_theorem}
	Let $(\phi_h^k,p_h^k)\in  X \times X$ be obtained by \eqref{CNFD}-\eqref{KKT_1} with \eqref{first scheme}. Under the assumption \eqref{assumption}, for small $\tau$ and $h$ satisfying a mild CFL type condition $\tau \leq C_0 h ~ (C_0>0)$, the following error estimation holds,
	\begin{equation*}
		\Vert e_\phi^k\Vert_{L^2}+\Vert e_p^k\Vert_{H^1}\leq C(\tau^2+h^2),\quad 0\leq k\leq T/\tau,
	\end{equation*}
	where $C > 0$ is a constant independent of $h$, $\tau$ and $k$.
\end{theorem}

Now define the local truncation errors $R_\phi^k, R_p^k\in  X (k\geq 0)$ as
\begin{equation}\label{R_0_error}
	\begin{split}
		R_\phi^0=&\frac{\phi_h(\tau)-\phi_h(0)}{\tau}-\frac{D}{4}\Delta_h\phi_h(\tau)+\frac{D}{2}\nabla_h\cdot\left(\frac{\phi_h(0)}{p_h(0)+p^0}\nabla_h(p_h(0)+p^0)\right)+(p_h(0)+p^0)\phi_h(0),\\
		R_p^0=&\frac{p_h(\tau)-p_h(0)}{\tau}-\eta D\Delta_h p_h(\tau) - (p_h(0)+p^0)\phi_h(0)+p_h(0),
	\end{split}
\end{equation}
and for $k\geq 1$
\begin{equation}\label{R_error}
	\begin{split}
		R_\phi^k=&\frac{\phi_h(t_{k+1})-\phi_h(t_k)}{\tau}-\frac{D}{8}\Delta_h(\phi_h(t_{k+1})+\phi_h(t_k))\\
		&+\frac{3D}{4}\nabla_h\cdot\left(\frac{\phi_h({t_k})}{p_h({t_k})+p^0}\nabla_h(p_h({t_k})+p^0)\right) - \frac{D}{4}\nabla_h\cdot\left(\frac{\phi_h({t_{k-1}})}{p_h({t_{k-1}})+p^0}\nabla_h(p_h({t_{k-1}})+p^0)\right)\\
		&+ \frac{3}{2}(p_h(t_k)+p^0)\phi_h(t_{k}) - \frac{1}{2} (p_h(t_{k-1})+p^0)\phi_h(t_{k-1}),\\
		R_p^k=&\frac{p_h(t_{k+1})-p_h(t_k)}{\tau}-\frac{\eta D}{2} \Delta_h\left(p_h(t_{k+1})+p_h(t_k)\right)-\frac{3}{2}((p_h(t_k)+p^0)\phi_h(t_{k})-p_h(t_k))\\
		&+\frac{1}{2}((p_h(t_{k-1})+p^0)\phi_h(t_{k-1})-p_h(t_{k-1})).
	\end{split}
\end{equation}
By Taylor's expansion, we can obtain the following estimates of the local errors.
\begin{lemma}\label{bound_R_phi_p}
	For the local truncation errors $R_\phi^k,R_p^k\in  X (k\geq 0)$ defined in \eqref{R_0_error}-\eqref{R_error}, under the assumption \eqref{assumption}, we have the estimates for $0\leq k\leq \frac{T}{\tau}-1$
	\begin{equation}\label{bound_R}
		\Vert R_\phi^{k+1}\Vert_{L^2}+ \Vert R_p^{k+1}\Vert_{H^1}\leq C(\tau^2+h^2),\quad 
		\Vert R_\phi^{0}\Vert_{L^2}+\Vert R_p^{0}\Vert_{H^1}\leq C(\tau+h^2),
	\end{equation}
	where $C$ is independent of $\tau$, $h$ and $k$.
\end{lemma}
\begin{proof}
	By the Taylor expansion approach, it is readily to check the bounds in \eqref{bound_R} when replacing $(\phi_h(t),p_h(t))$ in \eqref{R_0_error}-\eqref{R_error} by $(\phi_h(t),p_h(t))$, and $R_\phi^0=\mathcal{O}(\tau_0+h^2)$. Similarly, we can get $R_\phi^{k+1}=\mathcal{O}(\tau^2+h^2)$ and consequently  have $\Vert R_\phi^0\Vert_{L^2}\leq C (\tau+h^2)$ and 
	$\Vert R_\phi^{k+1}\Vert_{L^2}\leq C (\tau^2+h^2)$. 
	The rest estimation can be derived similarly.
\end{proof}

The following lemma characterizes the relation between the biased error functions \eqref{error_denfinition} under the $L^2$-$H^1$ projection.
\begin{lemma}[Lemma 4.3 in \cite{TongFenghua2024Positivity}]\label{error_projection}
	For the errors defined in \eqref{error_denfinition}, it holds that
	\begin{equation*}
		\Vert e_\phi^k\Vert_{L^2}^2+\Vert e_\phi^k-\tilde{e}_\phi^k\Vert_{L^2}^2\leq \Vert \tilde{e}_\phi^k\Vert_{L^2}^2,\quad 
		\Vert e_p^k\Vert_{H^1}^2+\Vert e_p^k-\tilde{e}_p^k\Vert_{H^1}^2\leq \Vert \tilde{e}_p^k\Vert_{H^1}^2,\quad 0\leq k\leq\frac{T}{\tau}.
	\end{equation*}
\end{lemma}
\begin{proof}
	$k=0$ is trivial. For $k\geq 1$, we have from \eqref{KKT_1} 
	\begin{equation*}
		e_\phi^k = \tilde{e}_\phi^k - \lambda_h^k,\quad
		(I-\Delta_h)e_p^k = (I-\Delta_h)\tilde{e}_p^k  -\xi_h^k.
	\end{equation*} 
	Taking the inner product of both sides with $e_\phi^k$ and $e_p^k$, respectively, we have
	\begin{equation*}
		\begin{aligned}
			&\frac{1}{2}(\Vert e_\phi^k\Vert_{L^2}^2+\Vert e_\phi^k-\tilde{e}_\phi^k\Vert_{L^2}^2-\Vert \tilde{e}_\phi^k\Vert_{L^2}^2)=-\langle \lambda_h^k,e_\phi^k\rangle,\\
			&\frac{1}{2}(\Vert e_p^k\Vert_{H^1}^2+\Vert e_p^k-\tilde{e}_p^k\Vert_{H^1}^2-\Vert \tilde{e}_p^k\Vert_{H^1}^2)=-\langle \xi_h^k,e_p^k\rangle.
		\end{aligned}
	\end{equation*}
	Using the KKT conditions and $\lambda_h^k\geq 0$, $\xi_h^k\geq 0$, we have $-\langle\lambda_h^k,e_\phi^k\rangle=-\langle\lambda_h^k,\phi_h(t_k)\rangle\leq 0$ $ (\phi_h(t)\geq 0)$ as well as  the estimate on $e_\phi^k$. The case for $e_p^k$ is similar and thus omitted here for brevity.
\end{proof}

By subtracting \eqref{CNFD} from \eqref{R_error}, we obtain the error equations for $k\geq 1$ as 
\begin{equation}\label{error_equation}
	\begin{aligned}
		\frac{\tilde{e}_\phi^{k+1}-e_\phi^k}{\tau}&=\frac{D}{8}\Delta_h(\tilde{e}_\phi^{k+1}+e_{\phi}^{k}) - \frac{3D}{4}T_1^k + \frac{D}{4}T_1^{k-1} - \frac{3}{2}U_1^k + \frac{1}{2}U_1^{k-1} + R_\phi^k,\\
		\frac{\tilde{e}_p^{k+1}-e_p^k}{\tau}&=\frac{\eta D}{2} \Delta_h(\tilde{e}_p^{k+1}+e_{p}^{k})+ \frac{3}{2}U_2^k -\frac{1}{2}U_2^{k-1}+R_p^k,
	\end{aligned}
\end{equation}
where $T_1^m, U_1^m,U_2^m\in  X  (m = k,k-1, k \geq 1)$ are defined as 
\begin{equation}\label{T_denfinition}
	T_1^m=\nabla_h\cdot\biggr(\frac{\phi_h({t_m})}{p_h(t_{m})+p^0}\nabla_h(p_h(t_{m})+p^0) - \frac{\phi_h^{m}}{p_h^{m}+p^0}\nabla_h(p_h^{m}+p^0)\biggr),
\end{equation}
and
\begin{equation}\label{U_denfinition}
	\begin{aligned}
		U_1^m&=p_h(t_m)e_\phi^m+e_p^m\phi_h(t_m)+p^0 e_\phi^m,\\
		U_2^m&=p_h(t_m)e_\phi^m+e_p^m\phi_h(t_m)+p^0 e_\phi^m - e_p^m.
	\end{aligned}
\end{equation}
For the nonlinear part \eqref{T_denfinition}, we denote
\begin{equation*}
	B=\max_{0\leq m\leq\frac{T}{\tau}}\{\Vert\phi_h(t_m)\Vert_{L^\infty}+\Vert p_h(t_m)\Vert_{L^\infty}+\Vert \nabla p_h(t_m)\Vert_{L^\infty}\},
\end{equation*}
where $B > 0$ is well defined under the assumption \eqref{assumption} and sufficiently small $h>0$.

The difficulty lies in the singular sensitivity in $T_1^m$. Fortunately, the Lagrange multiplier correction provides a posterior estimation on $e_p^m = p_h(t_m) - p_h^m$ and ensures positivity of the numerical solution $p_h^m$. 
\begin{lemma}\label{about T}
	Assuming $\Vert \phi_h^m\Vert_{L^\infty}+\Vert p_h^m\Vert_{L^\infty}+\Vert \nabla_h p_h^m\Vert_{L^\infty}\leq B+1(m=k,k-1,k\geq 1)$,  under the assumption \eqref{assumption}, for $T_1^m$ given in \eqref{T_denfinition} and any $f_h\in  X $, we have
	\begin{equation*}
		\begin{aligned}
			\vert\langle T_1^m,f_h\rangle\vert \leq C_{B,m_p,\tilde{C}}(\Vert e_\phi ^ {m} \Vert_{L^2} + \Vert e_p^{m}\Vert_{H^1} ) \Vert \nabla_h f_h \Vert_{L^2},
		\end{aligned}
	\end{equation*}
	where $C_{B,m_p,\tilde{C}}$ is a constant depending on $B$ ,$m_p$ and $\tilde{C}$.
\end{lemma}
\begin{proof}
	Let $a:=p_h(t_m)+p^0\geq m_p>0$. Note that \(a-e_p^m=p_h^m+p^0\geq m_p>0\). By definition \eqref{T_denfinition} we have
	\begin{equation}\label{T_e_1}
		\begin{aligned}
			T_1^m&=\nabla_h \cdot \Big(\frac{\phi_h(t_m)}{a}\nabla_h a-\frac{\phi_h^m}{a-e_p^m}(\nabla_h a-\nabla_h e_p^m)\Big)\\    
			&= \nabla_h \cdot \Bigg(\frac{e_\phi^m\,a-\phi_h(t_m)\,e_p^m}{a(a-e_p^m)}\nabla_h a +\frac{\phi_h^m}{a-e_p^m}\nabla_h e_p^m \Bigg)\\
			&= \nabla_h \cdot \Big(\underbrace{ \frac{e_\phi^m}{a-e_p^m}\nabla_h a}_{V_1} \Big)+ \nabla_h\cdot\Big(\underbrace{-\frac{\phi_h(t_m)\,e_p^m}{a(a-e_p^m)}\nabla_h a }_{V_2}\Big)+ \nabla_h\cdot\Big(\underbrace{\frac{\phi_h^m}{a-e_p^m}\nabla_h e_p^m}_{V_3}\Big).
		\end{aligned}
	\end{equation}
	Under assumptions in Lemma~\ref{about T}, since $p^0$ is sufficiently smooth and $\Vert \nabla_h p^0\Vert_{L^\infty}\leq \tilde{C}$, we have
	\begin{equation}\label{T_e_2}
		\begin{aligned}
			&\|V_1\|_{L^2}\le \|\tfrac{\nabla_h a}{a-e_p^m} \|_{L^\infty}\,\|e_\phi^m\|_{L^2}=\|\tfrac{\nabla_h a}{p_h^m+p^0} \|_{L^\infty}\,\|e_\phi^m\|_{L^2} \le \frac{B+\tilde{C}}{m_p} \|e_\phi^m\|_{L^2},\\
			&\Vert V_2\Vert_{L^2}\le\Vert \frac{\phi_h(t_m)}{a(a-e_p^m)}\nabla_h a\Vert_{L^{\infty}}\Vert e_p^m\Vert_{L^2}\leq \frac{B(B+\tilde{C})}{m_p^2}\Vert e_p^m\Vert_{L^2},\\
			& \|V_3\|_{L^2}\le \Vert\frac{\phi_h^m}{a-e_p^m} \Vert_{L^{\infty}}\Vert\nabla_h e_p^m\Vert_{L^2}\le \frac{B+1}{m_p}\,\|\nabla_h e_p^m\|_{L^2}.
		\end{aligned}
	\end{equation}
	Combining \eqref{T_e_1}, \eqref{T_e_2} and applying discrete summation by parts \eqref{summation_by_parts} to the inner product with $f_h$, we obtain
	\[
	|\langle T_1^m,f_h\rangle|
	\le C_{B,m_p,\tilde{C}}(\|e_\phi^m\|_{L^2}+\|e_p^m\|_{H^1})\|\nabla_h f_h\|_{L^2},
	\]
	which completes the proof.
\end{proof}

Now, we proceed to prove the main theorem.

\begin{proof}[Proof of Theorem \ref{main_theorem}] We shall prove by induction that for sufficiently small $h$ and $\tau$ satisfying $\tau\leq C_0h$ ($C_0>0$ is a constant),
	\begin{equation}\label{bound}
		\Vert \phi_h^k\Vert_{L^\infty}+\Vert p_h^k\Vert_{L^\infty}+\Vert \nabla _h p_h^k\Vert_{L^\infty}\leq B+1,\quad 0 \le k \le \frac{T}{\tau},	
	\end{equation} 
	and 
	\begin{equation}\label{subject}
		\Vert e_\phi^k\Vert_{L^2}^2+\Vert e_p^k\Vert_{H^1}^2+\Vert\tilde{e}_\phi^k\Vert_{L^2}^2+\Vert\tilde{e}_p^k\Vert_{H^1}^2\leq Ce^{Ck\tau}(\tau^2+h^2)^2,\quad 0\leq k\leq\frac{T}{\tau},
	\end{equation}
	where $C$ is constant (to be determined later) which is independent of $\tau$, $h$ and $k$. 
	
	\noindent \textbf{Step 1.} For $k = 0$, due to the error definition \eqref{error_denfinition} , we have $e_\phi^k = e_p^k = \tilde{e}_\phi^k = \tilde{e}_p^k =0$. Then for
	sufficiently small $h > 0$, we have \eqref{bound} and \eqref{subject}.\\
	\noindent \textbf{Step 2.} For $k = 1$, by subtracting \eqref{first scheme} from \eqref{R_0_error}, we have the error equations
	\begin{equation}\label{k1}
		\begin{aligned}
			&\frac{\tilde{e}_\phi^1-e_\phi^0}{\tau} = \frac{D}{4}\Delta_h\tilde{e}_\phi^1 - \frac{D}{2}T_1^0 - U_1^0 +R_\phi^0,\quad \frac{\tilde{e}_p^1-e_p^0}{\tau} = \eta D\Delta_h\tilde{e}_p^1 + U_2^0 + R_p^0,
		\end{aligned}
	\end{equation}
	where we have used the fact that $\phi_h^0=\phi_h(0)$ and $p_h^0=p_h(0)$. Under the assumption \eqref{assumption}, we take the $L^2$ inner
	product of \eqref{k1} with $2\tau\tilde{e}_\phi^1$ and $2\tau(I-\Delta_h)\tilde{e}_p^1$, respectively, and use the Cauchy inequality, then it yields
	\begin{equation*}
		\begin{aligned}
			&\Vert \tilde{e}_\phi^1\Vert_{L^2}^2 - \Vert e_\phi^0\Vert_{L^2} ^2+ \Vert \tilde{e}_\phi^1 - e_\phi^0\Vert_{L^2}^2 + \frac{D\tau}{2}\Vert \nabla_h \tilde{e}_\phi^1\Vert_{L^2}^2 = 2\tau \langle \tilde{e}_\phi^1,R_\phi^0\rangle - D\tau\langle \tilde{e}_\phi^1,T_1^0\rangle - 2\tau\langle \tilde{e}_\phi^1,U_1^0\rangle\\
			&\leq \frac{1}{2}\Vert \tilde{e}_\phi^1\Vert_{L^2}^2 + \frac{D\tau}{4} \Vert\nabla_h\tilde{e}_\phi^1\Vert_{L^2}^2 + 4\tau^2\Vert R_\phi^0\Vert_{L^2}^2+ (4\tau^2C_1+D\tau C_{B,p^0})(\Vert e_\phi^0\Vert_{L^2}^2+\Vert e_p^0\Vert_{H^1}^2),
		\end{aligned}
	\end{equation*}
	and
	\begin{equation*}
		\begin{aligned}
			&\Vert \tilde{e}_p^1\Vert_{H^1}^2 - \Vert e_p^0\Vert_{H^1} ^2+ \Vert \tilde{e}_p^1 - e_p^0\Vert_{H^1}^2 + 2\eta D\tau\Vert \nabla_h \tilde{e}_p^1\Vert_{H^1}^2 = 2\tau \langle (I-\Delta_h)\tilde{e}_p^1,R_p^0\rangle + 2\tau\langle (I-\Delta_h)\tilde{e}_p^1,U_2^0\rangle\\
			&\leq \frac{1}{2}\Vert \tilde{e}_p^1\Vert_{H^1}^2+4\tau^2\Vert R_p^0\Vert_{H^1}^2 + (4\tau^2C_1+1)(\Vert e_\phi^0\Vert_{L^2}^2+\Vert e_p^0\Vert_{H^1}^2),
		\end{aligned}
	\end{equation*}  
	where $C_1>0$ is a constant depending on $\phi(\bx,t), p(\bx,t), M_p$.
	From Lemma \ref{bound_R_phi_p}, we get
	$\Vert \tilde{e}_\phi^1\Vert_{L^2}^2 + \Vert \tilde{e}_p^1\Vert_{H^1}^2 = \mathcal{O}(\tau^2(\tau+h^2)^2 )$, which leads to $\Vert \tilde{e}_\phi^1\Vert_{L^2}+\Vert \tilde{e}_p^1\Vert_{H^1} = \mathcal{O}(\tau^2+h^2)$. Lemma \ref{error_projection} implies that for some constant $C_2>0$
	\begin{equation}\label{phi p 1}
		\Vert e_\phi^1\Vert_{L^2}^2+\Vert e_p^1\Vert_{H^1}^2+\Vert\tilde{e}_\phi^1\Vert_{L^2}^2+\Vert\tilde{e}_p^1\Vert_{H^1}^2\leq C_2(\tau^2+h^2)^2,
	\end{equation}
	i.e., Eq.~\eqref{subject} holds. Again, for any grid functions $u_h\in  X $, an application of $2$-D inverse inequality implies that $\Vert u_h\Vert _{\infty}\leq Ch^{-1} \Vert u_h\Vert_2$
	for $\tau\leq C_0h$ and $h>0$ ($h$ is sufficiently small). Thus we have
	\eqref{bound} for $k = 1$.
	
	\noindent \textbf{Step 3.} Now, assume \eqref{bound} and \eqref{subject} hold for $k \leq m (1 \leq m \leq T/\tau -1)$. We are going to prove the case $k = m+1$.
	
	Taking the inner products of \eqref{error_equation} with $\tilde{e}_\phi^{k+1}+e_\phi^k\in  X $ and $(I-\Delta_h)(\tilde{e}_p^{k+1}+e_p^k)\in  X $, respectively, applying Sobolev inequality, we have the $L^2$ error equation for $k\geq 1$,
	\begin{equation}\label{total}
		\begin{split}
			\frac{1}{\tau}(\Vert \tilde{e}_\phi^{k+1}\Vert_{L^2}^2-\Vert e_\phi^k\Vert_{L^2}^2)&+\frac{D}{8}\Vert\nabla_h(\tilde{e}_\phi^{k+1}+e_\phi^k)\Vert_{L^2}^2 = \underbrace{\langle -\frac{3D}{4}T_1^k+\frac{D}{4}T_1^{k-1},\tilde{e}_\phi^{k+1}+e_\phi^k\rangle}_{E_1} \\
			&+ \underbrace{\langle -\frac{3}{2}U_1^k+\frac{1}{2}U_1^{k-1},\tilde{e}_\phi^{k+1}+e_\phi^k\rangle}_{E_2} + \underbrace{\langle R_\phi^k,\tilde{e}_\phi^{k+1}+e_\phi^k\rangle}_{E_4},\\
			\frac{1}{\tau}(\Vert \tilde{e}_p^{k+1}\Vert_{H^1}^2-\Vert e_p^k\Vert_{H^1}^2)& + \frac{\eta D}{2}\Vert\nabla_h(\tilde{e}_p^{k+1}+e_p^k)\Vert_{L^2}^2 + \frac{\eta D}{2}\Vert\Delta_h(\tilde{e}_p^{k+1}+e_p^k)\Vert_{L^2}^2\\
			&= \underbrace{\langle \frac{3}{2}U_2^k-\frac{1}{2}U_2^{k-1},(I-\Delta_h)(\tilde{e}_p^{k+1}+e_p^k)\rangle}_{E_3} + \underbrace{\langle R_p^k,(I-\Delta_h)(\tilde{e}_p^{k+1}+e_p^k)\rangle}_{E_5}.
		\end{split}
	\end{equation}
	
	In view of Lemma \ref{about T}, under the induction hypothesis and applying the Cauchy inequality, we obtain 
	\begin{equation}\label{e2}
		\begin{aligned}
			E_1 &\leq C_{B,m_p,D,\tilde{C}}\Vert \nabla_h (\tilde{e}_\phi^{k+1}+e_\phi^k)\Vert_{L^2}\sum_{l=k,k-1}(\Vert e_\phi^l\Vert_{L^2}+\Vert e_p^l\Vert_{H^1})\\
			&\leq \frac{D}{16}\Vert \nabla_h(\tilde{e}_\phi^{k+1}+e_\phi^k)\Vert_{L^2}^2 + C_{B,m_p,D,\tilde{C},\Omega}\sum_{l=k,k-1}(\Vert e_\phi^l\Vert_{L^2}^2+\Vert e_p^l\Vert_{H^1}^2),
		\end{aligned}
	\end{equation}
	where $C_{B,m_p,D,\tilde{C},\Omega}$ is a constant depending on $B$, $p^0$, $D$, $\tilde{C}$ and $\Omega$. 
	
	
	An application of the Cauchy inequality yields the following series of estimates
	\begin{equation}\label{e3}
		\begin{aligned}
			E_2 \leq  &\frac{1}{16}\Vert \tilde{e}_\phi^{k+1}\Vert_{L^2} ^2+4 C_{B,M_p,\Omega} \sum_{l=k,k-1}(\Vert e_\phi^l\Vert_{L^2}^2 + \Vert e_p^l\Vert_{H^1}^2),\\
			E_3 \leq &\frac{1}{16}\Vert \tilde{e}_p^{k+1}\Vert_{H^1}^2+\frac{\eta D}{8}\Vert \Delta_h(\tilde{e}_p^{k+1}+e_p^k)\Vert_{L^2}^2+(8+\frac{2}{\eta D}) C_{B,M_p,\Omega} \sum_{l=k,k-1}(\Vert e_\phi^l\Vert_{L^2}^2 + \Vert e_p^l\Vert_{H^1}^2),\\
			E4 \leq &\frac{1}{16}\Vert \tilde{e}_\phi^{k+1}\Vert_{L^2}^2+ \frac{1}{16}\Vert e_\phi^{k}\Vert_{L^2}^2+ 8\Vert R_\phi^k\Vert_{L^2}^2,\\
			E5\leq &\frac{1}{16}\Vert \tilde{e}_p^{k+1}\Vert_{H^1}^2 +\frac{1}{16} \Vert e_p^{k}\Vert_{H^1}^2 + (8 + \frac{2}{\eta D})\Vert R_p^k\Vert_{H^1}^2 +\frac{\eta D}{8}\Vert\Delta_h(\tilde{e}_p^{k+1}+e_p^k)\Vert_{L^2}^2.
		\end{aligned}
	\end{equation}
	
	Substituting \eqref{e2} and \eqref{e3} from \eqref{total}, we obtain the $L^2$-$H^1$ error estimation
	\begin{equation}\label{e4}
		\begin{aligned}
			(\frac{1}{\tau} - \frac{1}{8})(\Vert \tilde{e}_\phi^{k+1}\Vert_{L^2}^2 & -\Vert e_\phi^k\Vert_{L^2}^2)+\frac{D}{16}\Vert\nabla_h(\tilde{e}_\phi^{k+1}+e_\phi^k)\Vert_{L^2}^2 \leq C_1(\Vert R_\phi^k\Vert_{L^2}^2+\sum_{l=k-1}^{k}(\Vert e_\phi^l\Vert_{L^2}^2+\Vert e_p^l\Vert_{H^1}^2)),\\
			(\frac{1}{\tau} - \frac{1}{8})(\Vert \tilde{e}_p^{k+1}\Vert_{H^1}^2& -\Vert e_p^k\Vert_{H^1}^2)+\frac{\eta D}{4}\Vert\nabla_h(\tilde{e}_p^{k+1}+e_p^k)\Vert_{H^1}^2
			\leq  C_2(\Vert R_p^k\Vert_{H^1}^2+\sum_{l=k-1}^{k}(\Vert e_\phi^l\Vert_{L^2}^2+\Vert e_p^l\Vert_{H^1}^2)),
		\end{aligned}
	\end{equation}
	where $C_1 = \max\{C_{B,m_p,D,\tilde{C},\Omega},4C_{B,M_p,\Omega} + 1\}$, $C_2 = \max\{C_{B,M_p,\Omega},(8 + \frac{2}{\eta D})C_{B,M_p,\Omega}+1\}$. 
	
	Denote by $S^k = \Vert \tilde{e}_\phi^{k}\Vert_{L^2}^2+\Vert \tilde{e}_p^{k}\Vert_{H^1}^2+\Vert e_\phi^{k}\Vert_{L^2}^2+\Vert e_p^{k}\Vert_{H^1}^2(k\geq 1)$. In view of \eqref{e4}, recalling Lemma \ref{error_projection} where $\Vert e_\phi^l\Vert_{L^2}^2\leq \Vert \tilde{e}_\phi^l\Vert_{L^2}^2$ and $\Vert e_p^l\Vert_{L^2}^2\leq \Vert \tilde{e}_p^l\Vert_{L^2}^2  (l = k,k-1)$, with sufficiently small $\tau$ and $h$, we have for $k\geq 1$
	\begin{equation}\label{e6}
		\begin{aligned}
			&\frac{1}{2\tau}(S^{k+1}-S^k)+\frac{D}{8}\Vert\nabla_h(\tilde{e}_\phi^{k+1}+e_\phi^k)\Vert_{L^2}^2+\frac{\eta D}{2}\Vert\nabla_h(\tilde{e}_p^{k+1}+e_p^k)\Vert_{H^1}^2\\
			&\leq 4C_3(\Vert R_\phi^k\Vert_{L^2}^2+\Vert R_p^k\Vert_{H^1}^2+\sum_{l=k,k-1}S^l),
		\end{aligned}
	\end{equation}
	where $C_3 = \max\{C_1,C_2\}$. Summing \eqref{e6} together for $1,2,\dots,k$, and using the local error in Lemma \ref{bound_R_phi_p} and the estimates \eqref{phi p 1} at the first step, for $1 \leq k \leq m$, we arrive at
	\begin{equation*}
		\begin{aligned}
			&S^{k+1}+\frac{D\tau}{4}\Vert\nabla_h(\tilde{e}_\phi^{k+1}+e_\phi^k)\Vert_{L^2}^2+\eta D\tau \Vert\nabla_h(\tilde{e}_p^{k+1}+e_p^k)\Vert_{H^1}^2\\
			&\leq S^1+8C_3\tau\sum_{l=1}^{k}(\Vert R_\phi^k\Vert_{L^2}^2+\Vert R_p^k\Vert_{H^1}^2+S^l) \leq 8C_3\tau\sum_{l=1}^{k}S^l+C_4(\tau^2+h^2)^2,
		\end{aligned}
	\end{equation*}
	where $C_4 > 0$ is a constant independent of $\tau$, $h$ and $k$. Using the discrete Gronwall inequality, it yields for some $ \tau >0 $,
	\begin{equation*}
		S^{k+1}\leq \exp(8C_3(k+1)\tau)C_4(\tau^2+h^2)^2,\quad 1\leq k\leq m.
	\end{equation*}
	So we have \eqref{subject} holds at $k = m + 1$, if we set $C = \max\{8C_3, C_4\}$. It is easy to check that the constant $C$ is independent of $h$, $\tau$ and $k$. The remaining is \eqref{bound} for $k = m+1$, which can be derived similarly as the $k = 0$ case by the inverse inequality $\Vert u_h\Vert _{\infty}\leq \tilde{C}h^{-1} \Vert u_h\Vert_2$. More precisely, for $k = m + 1$, \eqref{bound} implies for $\tau\leq C_0h$ with $h >  0$
	\begin{equation*}
		\begin{aligned}
			\Vert e_\phi^{k+1}\Vert_\infty + \Vert e_p^{k+1}\Vert_\infty &\leq h^{-1}(\Vert e_\phi^{k+1}\Vert_{L^2} + \Vert e_p^{k+1}\Vert_{H^1}) \leq 3\sqrt{Ce^{CT}}(h^{-1}\tau^2+h)\leq 1.
		\end{aligned}
	\end{equation*}
	The triangular inequality implies \eqref{bound} at $k = m + 1$. By the induction process, it completes the proof of Theorem \ref{main_theorem}.
\end{proof}

\section{SPIMEX for epidemic modeling}\label{sec:epi}

SPIMEX scheme can be extended to epidemic modeling involving multiple agents for investigating the impact of heterogeneous human behavior factors on the dynamics of infectious diseases. Here we consider a real example in \cite{XiongWangZhang2024}, in which the population density is divided into eight compartments: the susceptible agent $ {S}$, the exposed agent  ${E}$, the infectious and pre-symptomatic agent ${P}$, the asymptomatic agent $ {A}$, the mildly infectious symptomatic agent $I^-$, the infectious and symptomatic agent $ {I}^+$, the hospitalized agent $ {H}$ and the recovered agent $ {R}$.
\begin{equation}
	\left\{
	\begin{aligned}
		&\partial_t  {S} =   \mathcal{R} {S} 
		-   { \lambda } ( \beta (    {P}+  {A} ) +  I^-+  I^+ )  {S  } + \delta _{R}   {R}, \\
		&\partial_t  {E}=     \mathcal{R} {E}+
		\lambda   ( \beta (  {P}  +  {A } ) +    I^-+  I^+ ) { S     }  -\alpha   {E}  ,\\
		&\partial_t  {P} =     \mathcal{R} {P} + \alpha   {E} -\eta^\prime  {P }, \\
		&\partial_t  {A} =    \mathcal{R} {A} + \eta ^\prime (1- \rho)  {P} - \delta_A   {A }, \\
		&\partial_t  I^-=     \mathcal{R} I^- +  \eta ^\prime \rho (1- p_H) { P }-\delta_I ^{-}  I^-, \\
		&\partial_t  I^+ =     \mathcal{R} { I } ^{+} + \eta ^\prime \rho  p_H  { P }-\delta_I^{+}   I^+, \\
		&\partial_t  H =    \delta_I^{+}   I^+ - \delta _{H}  {H}, \\
		&\partial_t  R =     \mathcal{R}    R+ \delta _A   {A} + \delta_I^{-}   I^-  + \delta_H  { H} - \delta _{R}  {R },\\ 
		& \partial_t p =  \eta D \Delta p  + \bar{\delta}^{+}_{{\mathcal P}} (S+E + P +A+R) - \bar{\delta}^{-}_{{\mathcal P}}  p.
	\end{aligned}
	\right.
\end{equation} 
The chemotaxis-type operator $\mathcal{R}$ is defined by Eq.~\eqref{chemotaxis_operator}. All the event types, the parameter values and their corresponding physical meanings are put in Tables \ref{tab:parameter1233_1} and \ref{biase} in \ref{app:ODE_PDE}.

The above PDE system exhibits typical features of epidemic dynamics, involving (possibly degenerate) diffusion, aggregation, and quadratic interaction. Now denote by ${\bm{\Phi}} (\bx, t):= (\phi_{i})_{7\times 1} (\bx, t) =(S(\bx, t),  E(\bx, t) , P(\bx, t) , A(\bx, t) ,  I^+(\bx, t) ,  I^-(\bx, t) , R(\bx, t))^{\top}$.
The general form of the chemotaxis epidemic PDE system with singular sensitivity reads that
\begin{equation}\label{Epidemic_PDE}
	\left\{
	\begin{aligned}
		&\frac{\partial {{\bm{\Phi}}}}{\partial t}   =   \mathcal{R} {{\bm{\Phi}}} + \mathcal{L} {{\bm{\Phi}}} + \mathcal{N} ({{\bm{\Phi}}},H), & \mathbf{x}\in \Omega,t>0,\\
		&\frac{\partial H}{\partial t} = \delta_I^{+}I^{+}  -\delta_H H, & \mathbf{x}\in \Omega,t>0,\\
		&\frac{\partial p }{\partial t} = \eta D \Delta p + \bar{\delta}^{+}_{{\mathcal P}} (S+E + P +A+R) - \bar{\delta}^{-}_{{\mathcal P}}  p, & \mathbf{x}\in \Omega,t>0, 
	\end{aligned}
	\right.
\end{equation}
and the initial conditions are given as
\begin{equation}
	{{\bm{\Phi}}}(\mathbf{x},0) = {{\bm{\Phi}}}_0(\mathbf{x})\geq 0,~ p(\mathbf{x},0) = p_0(\mathbf{x})\geq 0,~ H(\mathbf{x},0) = H_0(\mathbf{x})\geq 0,\quad \mathbf{x}\in \Omega,
\end{equation}
where
\begin{equation}
	\mathcal{L} =
	\begin{pmatrix}
		0 	&	0	&	0 					&	0		&	0			&	0			&	&  	\delta_R \\
		0 	&	-\alpha	&	0 					&	0		&	0			&	0			&	&  	0 \\
		0 	&	\alpha	&	-\eta^{\prime} 			&	0		&	0			&	0			&	&  	0 \\
		0 	&	0	&	\eta^{\prime}(1-\rho) 		&	-\delta_A	&	0			&	0			&	&  	0 \\
		0 	&	0	&	\eta^{\prime}\rho(1-p_H) 	&	0		&	-\delta_I^{-}	&	0			&	&  	0 \\
		0 	&	0	&	\eta^{\prime}\rho p_H 	&	0		&	0			&	-\delta_I^{+}	&	&  	0 \\
		0 	&	0	&	0				 	&	\delta_A	&	\delta_I^{-}	&	0			&	&	-\delta_R	 \\
	\end{pmatrix},
\end{equation}
and 
\begin{equation}
	\mathcal{N} ({{\bm{\Phi}}},H) = (-   { \lambda }( \beta (    {P}+  {A} ) +  {  I  }^{-}+  {  I  }^{+})S,  ~{ \lambda } ( \beta (    {P}+  {A} ) +  {  I  }^{-}+  {  I  }^{+})S, ~0, ~0, ~0, ~0, ~\delta_H H)^{\top}.
\end{equation} 
In the absence of the attractiveness field $p(\bx, t)$, Eq.~\eqref{Epidemic_PDE} features the epidemic PDE system involving degenerated diffusion in $H(\bx, t)$ as the hospitalized agents are usually immobile.

Under the periodic boundary condition, it satisfies the positivity constraints and the mass conservation
\begin{equation}\label{mass_conservation}
	\psi_i(\mathbf{x},t) \geq 0, \quad p(\mathbf{x},t)\geq 0,
	\quad\frac{\mathrm{d}}{\mathrm{d} t}\int_{\Omega } \sum_{i=1}^{8}\psi_i(\mathbf{x},t) \mathrm{d} \bx = 0, 
\end{equation}
where $\bm{\Psi}:=(\psi_i)_{8\times1} = ((\phi_i)_{7\times1},H)^{\top}= (S,E,P,A,I^+,I^-,R, H)^{\top}$.

\subsection{FDMs with \texorpdfstring{$L^2$-$H^1$}{L2-H1} projection}

Now we introduce the FDMs based on the prediction-correction strategy for solving the epidemic model and give the error estimation. Similarly, we focus on the 2D case with $\Omega = [a_{\min},a_{\max}]\times [b_{\min},b_{\max}]$, and choose the time step size $\tau>0$, and the time steps are $t_k:=k\tau$ for $k=0,1,\dots$. Let $(\bm{\Psi}_h^k, p_h^k)\in X^8\times X$ be the numerical approximation to the exact solution $(\bm{\Psi},p)$ of \eqref{Epidemic_PDE} on $\Omega_h$ at time $t_k$. 

For $k = 0,1,2,\dots$, the prediction-correction strategy is given as follows.		
\begin{itemize}
	\item[Step 1.]	Intermediate solutions $\tilde{\bm{\Psi}}_h^{k+1} = 
	\left(
	\begin{array}{c}
		\tilde{{\bm{\Phi}}}_h^{k+1} \\
		\tilde{H}_h^{k+1}
	\end{array}
	\right)\in X^8$, $\tilde{p}_h^{k+1}\in  X $ are obtained via the Crank-Nicolson FDM with the Adams-Bashforth strategy for the nonlinear part: 
	\begin{equation}\label{epidemic_pre_step}
		\begin{aligned}
			\frac{\tilde{{{\bm{\Phi}}}}_h^{k+1}-{{\bm{\Phi}}}_h^k}{\tau} &= \frac{D}{8}\Delta_h(\tilde{{{\bm{\Phi}}}}_h^{k+1}+{{\bm{\Phi}}}_h^k) - \frac{3D}{4}\nabla_h\cdot (\frac{{{\bm{\Phi}}}_h^k}{p_h^k+p^0}\nabla_h (p_h^k+p^0)) \\
			&+ \frac{D}{4}\nabla_h\cdot \left(\frac{{{\bm{\Phi}}}_h^{k-1}}{p_h^{k-1}+p^0}\nabla_h (p_h^{k-1}+p^0)\right) + \frac{3}{2}\mathcal{L} {{\bm{\Phi}}}_h^k -\frac{1}{2}\mathcal{L}{{\bm{\Phi}}}_h^{k-1}\\
			&+ \frac{3}{2}\mathcal{N}({{\bm{\Phi}}}_h^k,H_h^k) - \frac{1}{2} \mathcal{N}({{\bm{\Phi}}}_h^{k-1},H_h^{k-1}),\\
			\frac{\tilde{H}_h^{k+1}-H_h^k}{\tau} &= \frac{3}{2} \delta_I^+(I^+)_h^k - \frac{1}{2} \delta_I^+(I^+)_h^{k-1} - \frac{3}{2}\delta_HH_h^k + \frac{1}{2}\delta_HH_h^{k-1},\\
			\frac{\tilde{p}_h^{k+1}-p_h^k}{\tau} &= \frac{\eta D}{2}\Delta_h(\tilde{p}_h^{k+1}+p_h^k) +\frac{3}{2}[ \bar{\delta}^{+}_{{\mathcal P}} (S_h^k +E_h^k + P_h^k + A_h^k + R_h^k) - \bar{\delta}^{-}_{{\mathcal P}}  p_h^{k}]\\
			&~~~ -\frac{1}{2}[ \bar{\delta}^{+}_{{\mathcal P}} (S_h^{k-1} +E_h^{k-1} + P_h^{k-1} + A_h^{k-1} + R_h^{k-1}) - \bar{\delta}^{-}_{{\mathcal P}}  p_h^{k-1}].
		\end{aligned}
	\end{equation}
	\item[Step 2.] Correct  the intermediate numerical solutions $\tilde{\bm{\Psi}}_{h}^{k+1}$ and $\tilde{p}_h^{k+1}$. 
	Here we treat the positivity for $\psi_{i,h}^{k+1}$ ($i=1, \dots, 8$), $p_h^{k+1}$ at the discrete level as the constraints for the numerical solution at $t_{k+1}$ with $\psi_{i,h}^{k+1}, p_h^{k+1}\geq 0$. At the same time, we also treat the mass conservation $\langle \sum_{i=1}^{8}\psi_{i,h}^{k+1},1\rangle = \langle \sum_{i=1}^{8}\psi_{i,h}^{0},1\rangle$. A natural approach to obtain $(\bm{\Psi}_{h}^{k+1},  p_h^{k+1})$ from the intermediate solutions $(\tilde{\bm{\Psi}}_{h}^{k+1}, \tilde{p}_h^{k+1})$ is to project the nodal vectors $\tilde{\bm{\Psi}}_{h}^{k+1}$ and $\tilde{p}_h^{k+1}$ onto a constrained manifold. The $L^2$ projection is applied to ensure both the positivity and mass conservation of $\bm{\Psi}_{h}^{k+1}$, while the $H^1$ projection is used to enforce the positivity of $p_h^{k+1}$, which reads
	\begin{equation}\label{min_pro_epi}
		\begin{aligned}
			&\min_{\psi_{i,h}^{k+1},p_h^{k+1}\in  X } \frac{1}{2}(\sum_{i=1}^{8}\Vert \psi_{i,h}^{k+1} - \tilde{\psi}_{i,h}^{k+1}\Vert_{L^2}^2+\Vert p_h^{k+1} - \tilde{p}_h^{k+1}\Vert_{H^1}^2),\\
			&\quad
			\text{s.t.} \quad \psi_{i,h}^{k+1}\geq 0 ~ (i=1,\dots,8), \quad  p_h^{k+1}\geq 0,\quad\langle \sum_{i=1}^{8}\psi_{i,h}^{k+1},1\rangle = \langle \sum_{i=1}^{8}\psi_{i,h}^{0},1\rangle. 
		\end{aligned}
	\end{equation}
	This is a convex minimization problem with the following KKT conditions, and can be solved efficiently by a simple semi-smooth Newton solver. 
	\begin{equation}\label{KKT_8}
		\begin{aligned}
			&{\bm{\Psi}}_{h}^{k+1} = \tilde{{\bm{\Psi}}}_{h}^{k+1}+\bm{\lambda}_h^{k+1}-\xi^{k+1}\mathbf{I}_{8\times1},\quad (I - \Delta_h) p_h^{k+1} =  (I - \Delta_h) \tilde{p}_h^{k+1}+\zeta_h^{k+1},\\
			& \lambda_{i,h}^{k+1}\psi_{i,h}^{k+1} = 0,\quad  \zeta_h^{k+1} p_h^{k+1} = 0,\quad  \lambda_{i,h}^{k+1}\geq 0,\quad \zeta_h^{k+1}\geq 0, \quad i = 1,\dots, 8,\\
			&\langle \sum_{i=1}^{8}\psi_{i,h}^{k+1},1\rangle = \langle \sum_{i=1}^{8}\psi_{i,h}^{0},1\rangle,
		\end{aligned}
	\end{equation}
	where $\zeta_h^{k+1}\in  X $, $\bm{\lambda}_h^{k+1} = (\lambda_{1,h}^{k+1},\dots,\lambda_{8,h}^{k+1})^\top\in  X^8 $ are the Lagrange multipliers for the positivity preservation, and $\xi^{k+1} \in \mathbb{R}$ is the Lagrange multiplier for the mass conservation. $I$ is the identity operator. 
\end{itemize}

Now, \eqref{epidemic_pre_step}--\eqref{KKT_8} constitute the SPIMEX scheme for epidemic dynamics. Since
\eqref{epidemic_pre_step} is linear and the convex minimization problem \eqref{min_pro_epi} admits unique solutions, it is still uniquely solvable at each time step; see \ref{app:Lagrange_method}.
\begin{remark}
	Since \eqref{epidemic_pre_step} is a three-level scheme, for the first step $k = 0$, we use the first order scheme instead
	\begin{equation}\label{epi_first_order}
		\begin{aligned}
			\frac{\tilde{{{\bm{\Phi}}}}_h^{1}-{{\bm{\Phi}}}_h^0}{\tau}&=\frac{D}{8}\Delta_h(\tilde{{{\bm{\Phi}}}}_h^{1}+{{\bm{\Phi}}}_h^0)-\frac{D}{2}\nabla_h\cdot\left(\frac{{{\bm{\Phi}}}_h^{0}}{p_h^{0}+p^0}\nabla_h(p_h^{0}+p^0)\right) +\mathcal{L}{{\bm{\Phi}}}_h^0 +\mathcal{N}({{\bm{\Phi}}}_h^0,H_h^0),\\
			\frac{\tilde{p}_h^{1}-p_h^0}{\tau}&=\frac{\eta D}{2}\Delta_h\left(\tilde{p}_h^{1}+p_h^0\right) + \bar{\delta}^{+}_{{\mathcal P}} (S_h^0 + E_h^0 + P_h^0 +A_h^0+R_h^0) - \bar{\delta}^{-}_{{\mathcal P}}  p_h^0,\\
			\frac{\tilde{H}_h^{1}-H_h^0}{\tau}& = \delta_I^+(I^+)_h^0  - \delta_HH_h^0 .
		\end{aligned}
	\end{equation}
\end{remark}

	Regarding the projection part \eqref{epi_first_order}, the following proposition states the positivity of the Lagrangian multiplier $\xi^{k+1}$.  In other words, the role of $\xi^{k+1}$ is to cancel out the increments of total mass induced by the Lagrangian multiplies $\bm{\lambda}_h^{k+1}$ which smooth out negative parts of ${{\bm{\tilde{\Psi}}}}_h^{k+1}$.
	\begin{proposition}\label{mass_phi}
		For \eqref{epidemic_pre_step}-\eqref{KKT_8}, the Lagrange multipliers $\xi^{k+1}\in \mathbb{R}$ satisfy
		\begin{equation}\nonumber
			\xi^{k+1} \geq 0,\quad k =0, 1, 2, \dots.
		\end{equation}
	\end{proposition}
	\begin{proof}
		For $k\geq 1$, combining \eqref{epidemic_pre_step} and \eqref{KKT_8}, we have
		\begin{equation}\label{mass_mult}
			\begin{split}
				\frac{{{\bm{\Phi}}}_{h}^{k+1}-{{\bm{\Phi}}}_{h}^k}{\tau} &= \frac{D}{8}\Delta_h(\tilde{{{\bm{\Phi}}}}_{h}^{k+1}+{{\bm{\Phi}}}_{h}^k) - \frac{3D}{4}\nabla_h\cdot (\frac{{{\bm{\Phi}}}_{h}^k}{p_h^k+p^0}\nabla_h (p_h^k+p^0))\\
				& + \frac{D}{4}\nabla_h\cdot (\frac{{{\bm{\Phi}}}_{h}^{k-1}}{p_h^{k-1}+p^0}\nabla_h (p_h^{k-1}+p^0)) + \frac{3}{2}\mathcal{L}{{\bm{\Phi}}}_{h}^k - \frac{1}{2}\mathcal{L}{{\bm{\Phi}}}_{h}^{k-1} \\
				& + \frac{3}{2}\mathcal{N}({{\bm{\Phi}}}_{h}^k,H_h^k)  - \frac{1}{2}\mathcal{N}({{\bm{\Phi}}}_{h}^{k-1},H_h^{k-1})+ \frac{\bm{\lambda}_{h}^{k+1}-\xi^{k+1}\mathbf{I}_{7\times 1}}{\tau},\\
				\frac{H_h^{k+1}-H_h^k}{\tau} &= \frac{3}{2} \delta_I^+(I^+)_h^k - \frac{1}{2} \delta_I^+(I^+)_h^{k-1} - \frac{3}{2}\delta_HH_h^k + \frac{1}{2}\delta_HH_h^{k-1} + \frac{\lambda_{8,h}^{k+1}-\xi^{k+1}}{\tau}.
			\end{split}
		\end{equation}
		Taking inner products of \eqref{mass_mult} with $\mathbf{1}$ on both sides and summing over all grid points, noticing $\lambda_{i,h}^{k+1}\geq 0$, we get $ \langle \sum_{i=1}^{8}(\lambda_{i,h}^{k+1}-\xi^{k+1}), 1\rangle = 0$, implying that $\xi^{k+1} = \frac{\langle \sum_{i=1}^{8}\lambda_{i,h}^{k+1},1\rangle}{\vert \Omega\vert}\geq 0$.  The case of $k=0$ is similar and omitted.         
	\end{proof}

	\subsection{Error estimation}
	Now we carry out the error analysis for \eqref{epidemic_pre_step}-\eqref{KKT_8} with \eqref{epi_first_order}.		
	
	Let $T>0$ be a fixed time, and $(\psi_i(\mathbf{x},t)\geq 0,p(\mathbf{x},t)\geq 0)$ be the exact solution of \eqref{Epidemic_PDE}. Based on the theoretical results, we make the following assumptions
	\begin{equation}\label{assumption_epidemic}
		\psi_i(\bx,t)\in C^3([0,T];C_{per}^4(\Omega)),\quad  p(\bx,t)\in C^3([0,T];C_{per}^5(\Omega)), \quad i=1,\dots,8,
	\end{equation}
	where $C_{per}^m(\Omega) = \{u\in C^m(\Omega)|\partial _x^k\partial_y^l u ~  \text{is periodic on}~ \Omega,\forall k,l\geq 0,k+l\leq m\}$. 
	
	We introduce the biased error functions $e_{\psi_{i}}^k$, $e_p^k$, $\tilde{e}_{\psi_{i}}^k$, $\tilde{e}_p^k\in  X (k\geq 0)$ ($i=1,\dots,8$):
	\begin{equation}\label{epi_error_definition}
		\begin{aligned}
			&e_{\psi_{i}}^k = \psi_{i,h}(t_k)-\psi_{i}^k,\quad e_{p}^k = p_h(t_k)-p_h^k,\quad  \tilde{e}_{\psi_{i}}^k = \psi_{i,h}(t_k)-\tilde{\psi}_{i}^k,\quad \tilde{e}_{p}^k = p_h(t_k)-\tilde{p}_h^k,
		\end{aligned}
	\end{equation}
	where $\tilde{e}_{\psi_{i}}^0 = e_{\psi_{i}}^0$, $\tilde{e}_p ^0 = e_p^0$ and  $e_{\bm{\Psi}}^k = (e_{\psi_1}^k,\dots,e_{\psi_8}^k)^\top$.
	The following error bounds can be established.
	
	\begin{theorem}\label{main_theorem_2}
		Let $(\bm{\Psi}_{h}^k,p_h^k)\in X$ be obtained by CNFDP \eqref{epidemic_pre_step}--\eqref{KKT_8} subject to \eqref{epi_first_order}.  
		Suppose $\mathcal{L}$ is a bounded operator, $\mathcal{N}(\bm{\Phi}, H)$ is Lipschitz continuous, and assumption \eqref{assumption_epidemic} holds.  
		Then for sufficiently small $\tau$ and $h$ satisfying the mild CFL-type condition $\tau \le C_0 h\ (C_0>0)$, the following error estimation holds
		\begin{equation*}
			\Vert e_{\bm{\Psi}}^k\Vert_{L^2}+\Vert e_{p}^k\Vert_{H^1}\leq C(\tau^2+h^2),
		\end{equation*}
		where $C > 0$ is a constant independent of $h$, $\tau$ and $k$.
	\end{theorem}
	
	Before proving the above theorem, we establish the relationship between the predicted errors $\tilde{e}_{\psi_i}^k~(i=1,\dots,8),~\tilde{e}_p^k$, and the corrected errors $e_{\psi_i}^k~(i=1,\dots,8), ~ e_p^k$.

	\begin{lemma}\label{error_psi}
		For the errors defined in \eqref{epi_error_definition}, it holds for $0\leq k\leq\frac{T}{\tau}$,
		\begin{equation*}
			\begin{aligned}
				&\sum_{i=1}^{8}\Vert e_{\psi_{i}} ^k\Vert_{L^2}^2+\sum_{i=1}^{8}\Vert e_{\psi_{i}}^k-\tilde{e}_{\psi_{i}}^k\Vert_{L^2}^2\leq \sum_{i=1}^{8}\Vert \tilde{e}_{\psi_{i}}^k\Vert_{L^2}^2,\\
				&\Vert e_p^k\Vert_{H^1}^2+\Vert e_p^k-\tilde{e}_p^k\Vert_{H^1}^2\leq \Vert \tilde{e}_p^k\Vert_{H^1}^2.
			\end{aligned}
		\end{equation*}
	\end{lemma}
	\begin{proof}
		$k=0$ is trivial. For $k\geq 1$, from \eqref{KKT_8}, we have  
		\begin{equation*}
			\begin{aligned}
				&e_{\psi_{i}}^k = \tilde{e}_{\psi_{i}}^k - \lambda_{i,h}^k + \xi^k,\quad
				(I-\Delta_h)e_p^k = (I-\Delta_h)\tilde{e}_p^k  -\zeta_h^k,\quad i=1,\dots,8.
			\end{aligned}
		\end{equation*} 
		Taking the inner product of both sides with $e_{\psi_{i}}^k$ and $e_p^k$, respectively, we have
		\begin{equation*}
			\begin{aligned}
				&\frac{1}{2}(\Vert e_{\psi_{i}}^k\Vert_{L^2}^2+\Vert e_{\psi_{i}}^k-\tilde{e}_{\psi_{i}}^k\Vert_{L^2}^2-\Vert \tilde{e}_{\psi_{i}}^k\Vert_{L^2}^2)=-\langle \lambda_{i,h}^k,e_{\psi_{i}}^k\rangle + \langle \xi^k,e_{\psi_{i}}^k\rangle ,\quad i=1,\dots,8,\\
				&\frac{1}{2}(\Vert e_p^k\Vert_{H^1}^2+\Vert e_p^k-\tilde{e}_p^k\Vert_{H^1}^2-\Vert \tilde{e}_p^k\Vert_{H^1}^2)=-\langle \zeta_h^k,e_p^k\rangle,
			\end{aligned}
		\end{equation*}
		where we use the fact that $\langle \sum_{i=1}^{8}e_{\psi_{i}}^k,\xi^k\rangle =0 $ due to the mass conservation.  Using the KKT conditions, we have $\lambda_{i,h}^k, \zeta_h^k\geq 0$ and $\lambda_{i,h}^{k+1}\psi_{i,h}^{k+1} = 0$.  Then it further yields $-\langle\lambda_{i,h}^k,e_{\psi_{i}}^k\rangle=-\langle\lambda_i^k,\psi_{i,h}(t_k)\rangle\leq 0$ ($\psi_{i,h}(t)\geq 0$) and the estimate on $e_{\psi_{i}}^k$ follows. The case of $e_p^k$ is similar and is omitted for brevity.
	\end{proof}

	Now, we proceed to complete the error estimation.
	\begin{proof}
		[Proof of Theorem \ref{main_theorem_2}]
		First, define the local truncation errors $R_{\bm{\Phi}}^k, R_H^k, R_p^k\in  X (k\geq 0)$ as 
		\begin{equation}\label{epi_R_k}
			\begin{split}
				R_{\bm{\Phi}}^k=&\frac{\bm{\Phi}_h(t_{k+1})-\bm{\Phi}_h(t_k)}{\tau}-\frac{D}{8}\Delta_h(\bm{\Phi}_h(t_{k+1})+\bm{\Phi}_h(t_k))\\
				&+\frac{3D}{4}\nabla_h\cdot\left(\frac{\bm{\Phi}_h({t_k})}{p_h({t_k})+p^0}\nabla_h(p_h({t_k})+p^0)\right) - \frac{D}{4}\nabla_h\cdot\left(\frac{\bm{\Phi}_h({t_{k-1}})}{p_h({t_{k-1}})+p^0}\nabla_h(p_h({t_{k-1}})+p^0)\right)\\
				&- \frac{3}{2}\mathcal{L}\bm{\Phi}_h(t_{k}) + \frac{1}{2} \mathcal{L}\bm{\Phi}_h(t_{k-1})- \frac{3}{2}\mathcal{N}(\bm{\Phi}_h(t_{k}),H_h(t_k)) + \frac{1}{2} \mathcal{N}(\bm{\Phi}_h(t_{k-1}),H_h(t_{k-1})),\\
				R_H^k = &\frac{H_h(t_{k+1})-H_h(t_k)}{\tau} - \frac{3}{2} \delta_I^+(I^+)_h(t_k) + \frac{1}{2} \delta_I^+(I^+)_h(t_{k-1}) + \frac{3}{2}\delta_HH_h(t_k) - \frac{1}{2}\delta_HH_h(t_{k-1}),\\
				R_p^k=&\frac{p_h(t_{k+1})-p_h(t_k)}{\tau}-\frac{\eta D}{2} \Delta_h\left(p_h(t_{k+1})+p_h(t_k)\right)\\
				&-\frac{3}{2}(\bar{\delta}^{+}_{{\mathcal P}} (S_h(t_k) +E_h(t_k) + P_h(t_k) + A_h(t_k) + R_h(t_k)) - \bar{\delta}^{-}_{{\mathcal P}}  p_h(t_k))\\
				&+\frac{1}{2}(\bar{\delta}^{+}_{{\mathcal P}} (S_h(t_{k-1}) +E_h(t_{k-1}) + P_h(t_{k-1}) + A_h(t_{k-1}) + R_h(t_{k-1})) - \bar{\delta}^{-}_{{\mathcal P}}  p_h(t_{k-1})),
			\end{split}
		\end{equation}
		for $k\geq 1$, and in particular, $k=0$,
		\begin{equation}\label{epi_R_0}
			\begin{split}
				R_{\bm{\Phi}}^0=&\frac{\bm{\Phi}_h(\tau)-\bm{\Phi}_h(0)}{\tau}-\frac{D}{4}\Delta_h\bm{\Phi}_h(\tau)+\frac{D}{2}\nabla_h\cdot\left(\frac{\bm{\Phi}_h(0)}{p_h(0)+p^0}\nabla_h(p_h(0)+p^0)\right)\\
				&-\mathcal{L}\bm{\Phi}_h(0)-\mathcal{N}(\bm{\Phi}_h(0),H_h(0)),\\
				R_H^0 =& \frac{H_h(\tau)-H_h(0)}{\tau}- \delta_I^+(I^+)_h(0) + \delta_HH_h(0),\\
				R_p^0=&\frac{p_h(\tau)-p_h(0)}{\tau}-\eta D\Delta_h p_h(\tau) - \bar{\delta}^{+}_{{\mathcal P}} (S_h(0) +E_h(0) + P_h(0) + A_h(0) + R_h(0)) + \bar{\delta}^{-}_{{\mathcal P}}  p_h(0).
			\end{split}
		\end{equation}
		Similar to Lemma \ref{bound_R_phi_p}, under the assumption \eqref{assumption_epidemic}, using Taylor's expansion, we have the estimates for $0\leq k\leq \frac{T}{\tau}-1$,
		\begin{equation}\label{epi_local}
			\begin{split}
				&
				\sum_{i=1}^{7}\Vert R_{\phi_i}^{k+1}\Vert_{L^2}+ \Vert R_p^{k+1}\Vert_{H^1}\leq C(\tau^2+h^2),\quad \Vert R_H^{k+1}\Vert_{L^2}\leq C\tau^2,
				\\
				&\sum_{i=1}^{7}\Vert R_{\phi_i}^{0}\Vert_{L^2}+\Vert R_p^{0}\Vert_{H^1}\leq C(\tau+h^2),\quad \Vert R_H^{0}\Vert_{L^2}\leq C\tau,
			\end{split}
		\end{equation}
		where $C$ is independent of $\tau$, $h$ and $k$.
		
		By subtracting \eqref{epidemic_pre_step} from  \eqref{epi_R_k} , we obtain the error equations for $k\geq 1$ as 
		\begin{equation*}
			\begin{aligned}
				\frac{{\tilde{e}_{\bm{\Phi}_h}}^{k+1}-e_{\bm{\Phi}_{h}}^k}{\tau} = &\frac{D}{8}\Delta_h(\tilde{e}_{\bm{\Phi}_h}^{k+1}+e_{\bm{\Phi}_h}^k) - \frac{3D}{4}T_1^k +\frac{D}{4} T_1^{k-1} + \frac{3}{2}Q_1^k- \frac{1}{2}Q_1^{k-1}+ R_{\bm{\Phi}}^k,\\
				\frac{\tilde{e}_{H_h}^{k+1}-e_{H_h}^k}{\tau} = &\frac{3}{2} \delta_I^+e_{I_h^+}^k - \frac{1}{2} \delta_I^+e_{I_h^+}^{k-1} - \frac{3}{2}\delta_He_{H_h}^k + \frac{1}{2}\delta_He_{H_h}^{k-1} + R_H^k,\\
				\frac{\tilde{e}_{p_h}^{k+1}-e_{p_h}^k}{\tau} = &\frac{\eta D}{2}\Delta_h(\tilde{e}_{p_h}^{k+1}+e_{p_h}^k) +\frac{3}{2}M_1^k  -\frac{1}{2} M_1^{k-1} +R_{p}^k,
			\end{aligned}
		\end{equation*}
		where $T_1^m, Q_1^m, M_1^m \in  X  (m = k,k-1, k \geq 1)$ are defined as 
		\begin{equation*}
			\begin{aligned}
				T_1^m&=\nabla_h\cdot\biggr(\frac{\bm{\Phi}_h(t_m)}{p_h(t_m)+p^0}\nabla_h (p_h(t_m)+p^0) - \frac{{{\bm{\Phi}}}_{h}^m}{p_h^m+p^0}\nabla_h (p_h^m+p^0) \biggr),\\
				Q_1^m&=\mathcal{L}({{\bm{\Phi}}}_{h}(t_m)- {{\bm{\Phi}}}_{h}^m)+ \mathcal{N}\biggr(({{\bm{\Phi}}}_{h}(t_m),H_h(t_m))-( {{\bm{\Phi}}}_{h}^m,H_h^m)\biggr),\\
				M_1^m &= \bar{\delta}^{+}_{{\mathcal P}} (e_{S_h}^{m} + e_{E_h}^{m} + e_{P_h}^{m} + e_{A_h}^{m} + e_{R_h}^{m}) - \bar{\delta}^{-}_{{\mathcal P}}  e_{p_h}^{m}.
			\end{aligned}
		\end{equation*}
		
		Following Theorem \ref{main_theorem},  together with Lemma \ref{about T} for $T_1^m$ and  Lemma \ref{error_psi} for the corrected error, we obtain that if $\mathcal{L}$ is bounded, $\mathcal{N}(\bm{\Phi},H)$ is Lipschitz continuous, and the mild CFL-type condition $\tau \le C_0 h\ (C_0>0)$ holds, then
		\begin{equation}\label{epi_fin_err}
			\sum_{i=1}^8 \Vert e_{\psi_i}^k\Vert_{L^2}^2+\Vert e_p^k\Vert_{H^1}^2 \leq \hat{C}(\tau^2+h^2)^2,  \quad 0\leq k\leq\frac{T}{\tau},
		\end{equation}
		where $\hat{C}$ is a constant independent of $\tau$, $h$, and $k$.
		
		Finally, combining this estimate \eqref{epi_fin_err} with \eqref{eq:component_sum} completes the proof of Theorem \ref{main_theorem_2}.
	\end{proof}
	
	\section{Numerical experiments}\label{sec:numerical}
	In this section, we first perform a series of numerical experiments to validate our theoretical error bounds.  To evaluate numerical errors, the discrete version of $L^2$ and $H^1$ norms are adopted as the metrics
	\begin{equation*}
		\|e\|_{L^2}
		= h\Bigl(\sum_{i=1}^N \sum_{j=1}^N e_{i,j}^2\Bigr)^{\!1/2}\,,
		\quad
		\|e\|_{H^1}
		= \Bigl(\Vert e\Vert_{L^2}^2 + \|\nabla_h e\|_{L^2}^2\Bigr)^{\!1/2},
	\end{equation*}
	where $e_{i,j}$ is the difference of numerical and reference solutions at node $(x_i,y_j)$. The reference solutions are produced under a fine grid mesh with spacings $\tau_{ref}$ and $h_{ref}$.  
	
	
	
	After that, we  investigate the competition of aggregation and diffusion  in both crime and epidemic models via the proposed SPIMEX scheme. Both the generation of stable hotspots and phase transition are observed and coincide with the early observations in \cite{ShortBrantinghamBertozziTita2010,ShortDorsognaPasourTitaBrantinghamBertozziChayes2008,XiongWangZhang2024,HaoMilyQuainiZhong2026}. In addition,  we provide an example to demonstrate the possible numerical blow-up when the positivity of solutions is absent, which confirms the necessity of structure-preserving property.

The theoretical convergence order is verified by simulating the chemotaxis system under a smooth initial condition and typical parameters. The first‐order SPIMEX scheme with much smaller time step is adopted to obtain the missing starting points. 

\begin{example}\label{ex:crime-convergence}
	Consider the chemotaxis model \eqref{cross_diffusion_PDE} in the domain $\Omega = [0,1]^2$ with periodic boundary conditions.  Let
	\begin{equation}\label{crime_initial_data}
		\begin{aligned}
			&  \phi_0(x,y) = 0.1\sin(2\pi (x-0.5))\cos(2\pi (y-0.5)) + 0.2,\quad\\
			&  p_0(x,y)    = 0.5\sin(2\pi (x-0.5))\cos(2\pi (y-0.5)) + 1,
		\end{aligned}
	\end{equation}
	and choose the model parameters $p^0 = 1/30,  \eta = 1, D   = 0.1.$ 
\end{example}
For the spatial error analysis, we employ a very small time step $\tau = 10^{-4}$ and terminate at $T = 0.1$, and let $(\phi_{h_{\mathrm{ref}}},p_{h_{\mathrm{ref}}})$ be the reference solutions on the finest mesh with $N=512$, $h_{ref} = 1/512$. Table \ref{space L2 phi} reports the $L^2$–$H^1$ errors and the convergence order.
\begin{table}[htbp]
	\centering
	\caption{\small Spatial error analysis with initial value \eqref{crime_initial_data} and $h_{ref}=1/512$.}
	\label{space L2 phi}
	\begin{tabular}{c|c|c|c|c|c}
		\hline\hline
		$h$  & $\Vert e_{\phi}\Vert_{L^2}$&  Order & $\Vert e_p\Vert_{H^1}$ & Order & Convergence \\ \hline
		$1/8$  & 2.0727$\times 10^{-4}$ & - & 3.5359$\times 10^{-2}$& - & \multirow{4}{*}{\centering\includegraphics[width=0.2\textwidth]{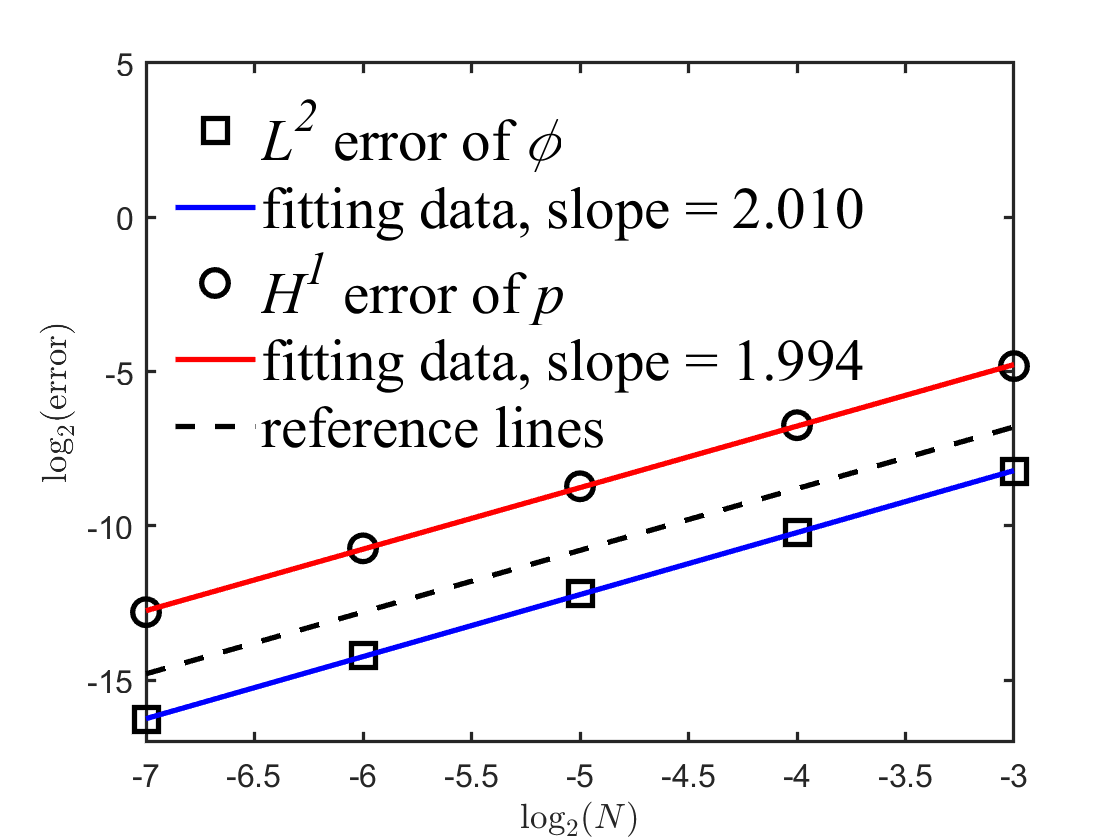}} \\
		$1/16$ & 5.3142$\times 10^{-5}$ &1.9636 &9.3838$\times 10^{-3}$ &1.9138 & \\ 
		$1/32$ & 1.3330$\times 10^{-5}$ & 1.9952 &2.3795$\times 10^{-3}$ &1.9795& \\ 
		$1/64$ & 3.2985$\times 10^{-6}$ & 2.0148 & 5.9054$\times 10^{-4}$ & 2.0105 & \\ 
		$1/128$ & 7.8566$\times 10^{-7}$ & 2.0698 & 1.4077$\times 10^{-4}$ & 2.0687 & \\ 
		\hline\hline
	\end{tabular}
\end{table}

For the temporal error analysis, we fix the spatial mesh at $h = 1/128$. Denote by $\phi_{\tau}$ and $p_{\tau}$ the numerical solutions at time step $\tau$, and use the solutions at $\tau_{\mathrm{ref}} = 10^{-5}$ as the references $\phi_{\tau_{\mathrm{ref}}}$ and $p_{\tau_{\mathrm{ref}}}$.  The final time is $T = 0.2$. Table~\ref{time L2 H1} summarizes the  $L^2$–$H^1$ errors and convergence rates. From the results, the second-order convergence in both spatial and temporal directions is verified.
\begin{table}[htbp]
	\centering
	\caption{\small Temporal error analysis with initial value \eqref{crime_initial_data} and $\tau_{ref}=10^{-4}$.}
	\label{time L2 H1}
	\begin{tabular}{c|c|c|c|c|c}
		\hline\hline
		$\tau$  & $\Vert e_{\phi}\Vert_{L^2}$ &  Order & $\Vert e_p\Vert_{ H^1}$ & Order & Convergence \\ \hline
		$2\times10^{-3}$   &  $4.5715\times10^{-7}$  & - &  $1.2283\times 10^{-5}$ & - & \multirow{4}{*}{\centering\includegraphics[width=0.2\textwidth]{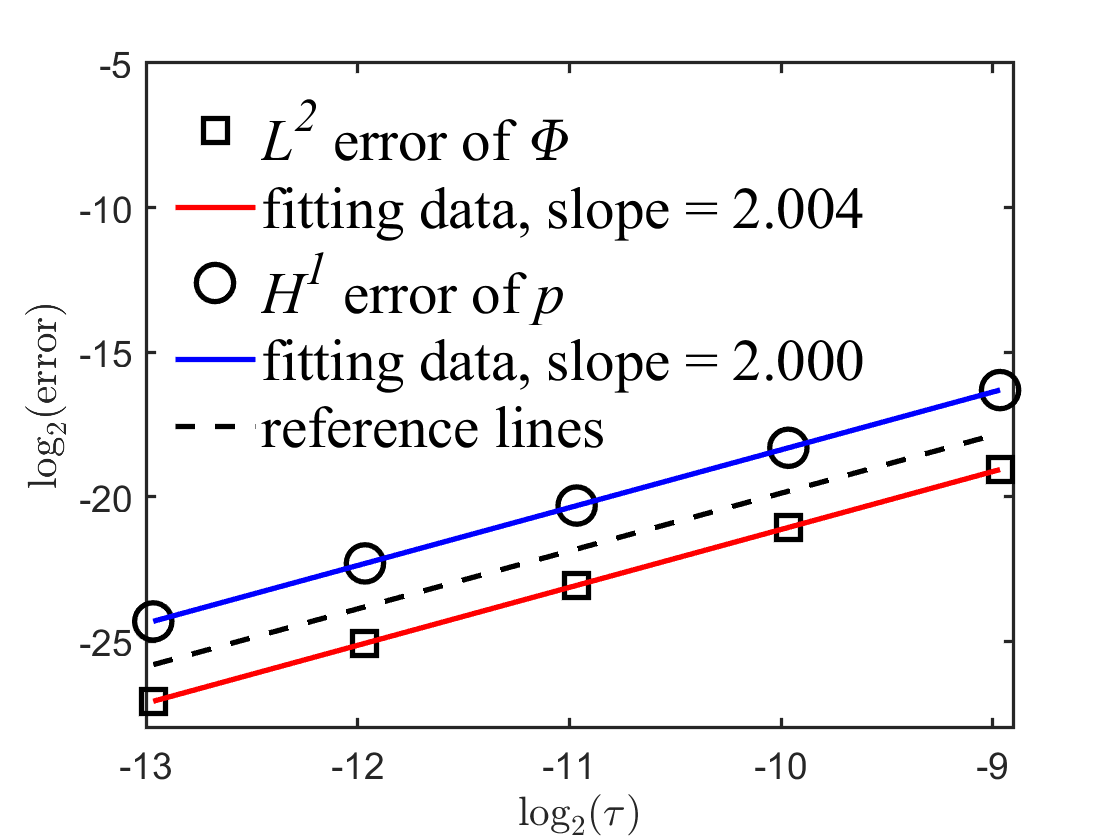}} \\
		$1\times10^{-3}$   & $1.1399\times10^{-7}$ &2.0038 & $3.0802\times10^{-6}$ &  1.9956& \\ 
		$5\times10^{-4}$   & $2.8453\times10^{-8}$ &2.0023  & $7.7103\times10^{-7}$ & 1.9982 & \\ 
		$2.5\times10^{-4}$  & $7.1002\times10^{-9}$ & 2.0026 & $1.9268\times10^{-7}$& 2.0006 & \\ 
		$1.25\times10^{-4}$  & $1.7660\times10^{-9}$ & 2.0074 & $4.7958\times10^{-8}$ & 2.0064 & \\ 
		\hline\hline
	\end{tabular}
\end{table}

\subsection{Pattern formation and phase transition in crime modeling}
		
		We consider the chemotaxis--type crime model \cite{ShortDorsognaPasourTitaBrantinghamBertozziChayes2008}
		\begin{equation}\label{crime_model}
			\left\{
			\begin{aligned}
				\rho_t
				&= D\nabla\!\cdot\!\Big(\nabla\rho-\frac{2\rho}{A+A^0}\nabla(A+A^0)\Big)-\rho(A+A^0)+\gamma,\\
				A_t
				&= \eta D\Delta A-\omega A+\kappa\rho(A+A^0),
			\end{aligned}
			\right.
		\end{equation}
		where $\rho(\bx,t)$ and $A(\bx,t)$ denote the criminal density and attractiveness,
		and $D,\eta,\omega,\kappa,\gamma,A^0>0$.

		The spatially homogeneous equilibrium of \eqref{crime_model} is
		\[
		\bar A=\frac{\kappa\gamma}{\omega},\qquad
		\bar Q=\bar A+A^0,\qquad
		\bar\rho=\frac{\gamma}{\bar Q}.
		\]
		Let $\rho=\bar\rho+u$, $A=\bar A+v$ and $U = (\rho, A)^T$.
		Seeking normal modes $(u,v)^T=Ue^{\sigma t+i k\cdot x}$ with $s=|k|^2$
		yields 
		\[ \sigma U=M(s)U,\quad \text{where}~~
		M(s)=
		\begin{pmatrix}
			-(Ds+\bar Q) & 2D\frac{\bar\rho}{\bar Q}s-\bar\rho\\
			\kappa\bar Q & \kappa\bar\rho-\omega-\eta Ds
		\end{pmatrix}.
		\]
		Hence
		\[
		\sigma^2-\mathrm{tr}(M(s))\,\sigma+\det(M(s))=0
		\]
		with
		\[
		\mathrm{tr}(M(s))=-(1+\eta)Ds-(\bar Q+\omega-\kappa\bar\rho),\quad
		\det(M(s))=\eta D^2s^2+D(\eta\bar Q+\omega-3\kappa\bar\rho)s+\omega\bar Q.
		\]
		
		Since $\mathrm{tr}(M(s))<0$ for all $s\ge0$ when
		$\bar Q+\omega-\kappa\bar\rho>0$, instability occurs if
		$\det(M(s))<0$ for some $s>0$, which is equivalent to
		\begin{equation}\label{condition}
			3\kappa\bar\rho-\eta\bar Q-\omega>2\sqrt{\eta\omega\bar Q}.
		\end{equation}
		In this case the unstable wave-number band is
		\[
		|k|^2\in(s_-,s_+),\qquad
		s_\pm=
		\frac{
			3\kappa\bar\rho-\eta\bar Q-\omega
			\pm
			\sqrt{(3\kappa\bar\rho-\eta\bar Q-\omega)^2-4\eta\omega\bar Q}
		}{2\eta D}.
		\]

\begin{figure}[h!]
			\centering
			\includegraphics[width=0.65\textwidth]{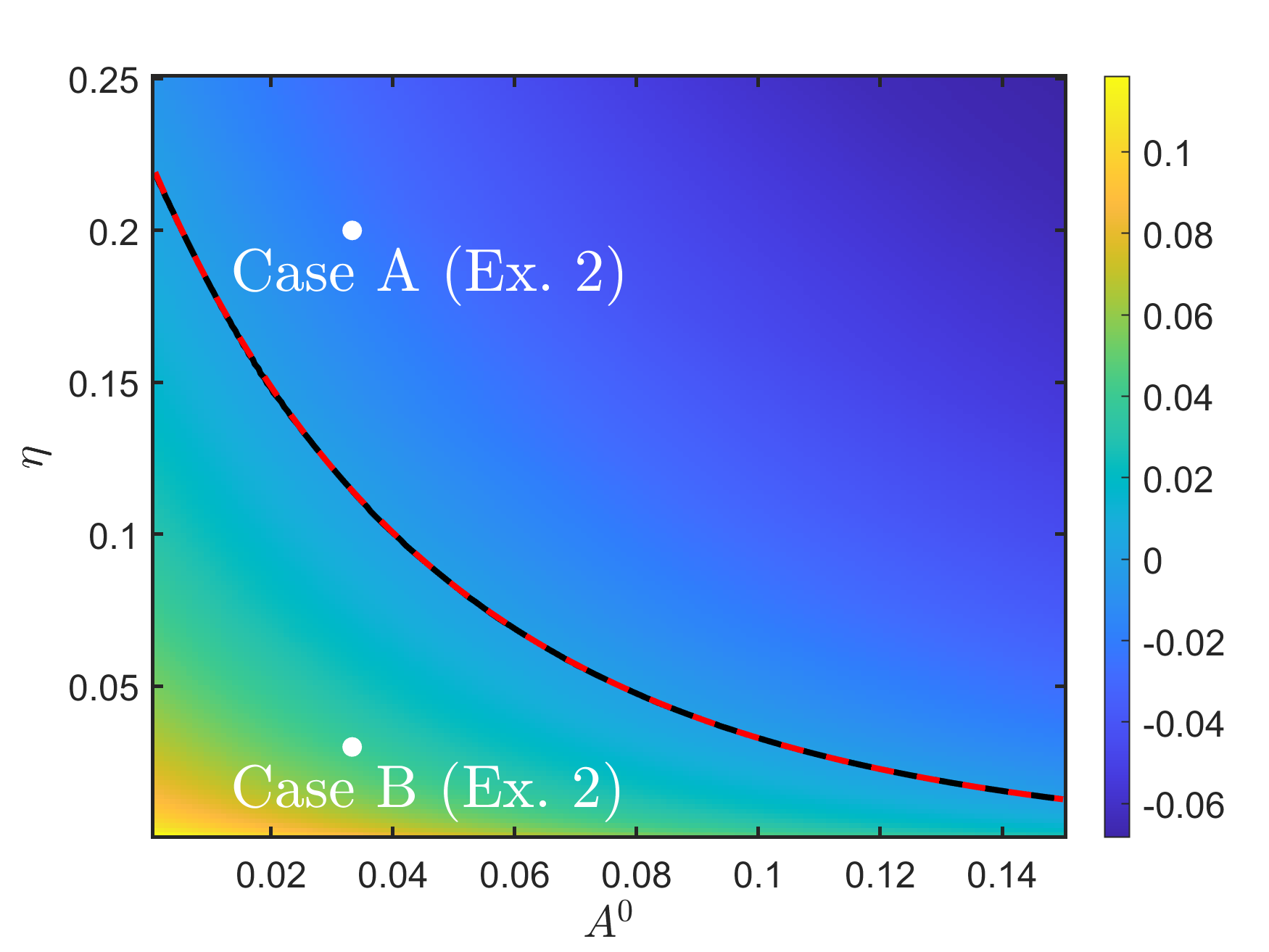}
			\vspace{-10pt}
			\caption{
			Phase diagram of $\sigma_{\max}(A^0, \eta)$ in the crime model.
			The stability boundary $\sigma_{\max}=0$ is shown by the black solid curve
			(numerical computation), which almost coincides with  analytical instability condition \eqref{condition} in
			the red dashed curve.
			This boundary separates the stable region ($\sigma_{\max}<0$) from the
			unstable region ($\sigma_{\max}>0$).
			In the unstable regime, spatial perturbations with wave numbers
			$|k|^2\in(s_-,s_+)$ grow exponentially, leading to hotspot formation.
			The white markers indicate the parameter choices $(A^0,\eta)$ used in
			Example~\ref{exa:crime} (Case~A and Case~B).
			}
			\label{crime_phase}
		\end{figure}

To illustrate the instability mechanism,
we resort to the maximal growth rate
\[
\sigma_{\max}(\eta,A^0)
=\max_{s>0}\,\Re\lambda_{\max}(M(s))
\]
over the parameter plane $(A^0,\eta)$.
The remaining parameters are fixed as
$\omega = 1/15$, $\kappa = 0.56$, and $\gamma = 0.019$.
The resulting phase diagram of $\sigma_{\max}(A^0, \eta)$ is shown in Figure~\ref{crime_phase}.
The  boundary between the stable and unstable regimes defined by $\sigma_{\max}=0$ is shown by the
black solid curve obtained from numerical computation, and the
red dashed curve presents the analytical result by \eqref{condition}.
The two curves almost coincide, confirming the theoretical prediction.

The diagram clearly separates the parameter region where the
spatially homogeneous equilibrium is stable ($\sigma_{\max}<0$)
from the region where it becomes unstable ($\sigma_{\max}>0$).
Being consistent with the analytical condition \eqref{condition},
instability occurs when $\eta$ and 
baseline attractiveness $A^0$ are sufficiently small.
As either $\eta$ or $A^0$ increases, the homogeneous equilibrium
becomes stable and spatial pattern formation is suppressed.
		
		To further verify the predictions of the linear stability analysis
		and the phase diagram in Figure \ref{crime_phase},
		we perform numerical simulations of the PDE system
		\eqref{crime_model}.
		In particular, we consider perturbations of the spatially homogeneous
	equilibrium and examine the subsequent evolution of the crime density
	and attractiveness fields.
				
		\begin{example}\label{exa:crime}
			The initial condition is set by adding small perturbations on the equilibrium
			\begin{equation}\label{init_condition}
				\rho(x, y, 0) =  \bar{\rho}~+ \delta (x, y), \quad
				A(x, y, 0) = \bar{A}~+ \delta (x, y),
			\end{equation}
			where the equilibrium of Eq.~\eqref{crime_model}, and the perturbation $\delta(x, y)$ is composed of $30$ independent Gaussian functions with randomly chosen centers $(x_i, y_i)$, heights $h_i$ and widths $\sigma_i$
			\begin{equation}\label{delta_set}
				\delta (x, y)  \approx \sum_{i=1}^{30}  \sum_{j= -L}^{L} \sum_{k = -L}^{L}  h_i \exp\left( \frac{-(x - x_i +j L_x)^2 - (y - y_i +k L_y)^2}{\sigma_i}\right),
			\end{equation}
			where $h_i = 0.02 r_i^{(1)}, \sigma_i = 0.005 r_i^{(2)}, x_i = x_{\min} + (x_{\max} - x_{\min})r_i^{(3)}, y_i = y_{\min} + (y_{\max} - y_{\min})r_i^{(4)}$, and $r_i^{(1)}$ to $ r_i^{(4)}$ are independent uniform random number in $[0, 1)$. The periodic images should be added to ensure the periodic boundary condition of all agent and field variables. With the same parameters as in Figure \ref{crime_phase}, we choose $\omega = {1}/{15}$, $\kappa = 0.56$, $\gamma = 0.019$, $\Omega = [0,2\pi]\times[0,2\pi]$, $h=2\pi/256$, $\Delta t=0.01$. 
		\end{example}
		Two representative parameter sets are selected according to the
		phase diagram in Figure \ref{crime_phase},
		\begin{itemize}
			\item[(i)] \text{Case A:} $\eta=0.2,\ A^0=1/30$;
			\item[(ii)] \text{Case B:} $\eta=0.03,\ A^0=1/30$.
		\end{itemize}
		
\begin{figure}[h!]
			\centering
		\subfigure[PDE system (Case A): Crime density $\rho$ at $T = 800$
		for $D = 0.001, 0.01, 0.1, 1$ ($\eta = 0.2$).\label{crime_plot_1}]
			{
				\includegraphics[width=0.25\textwidth]{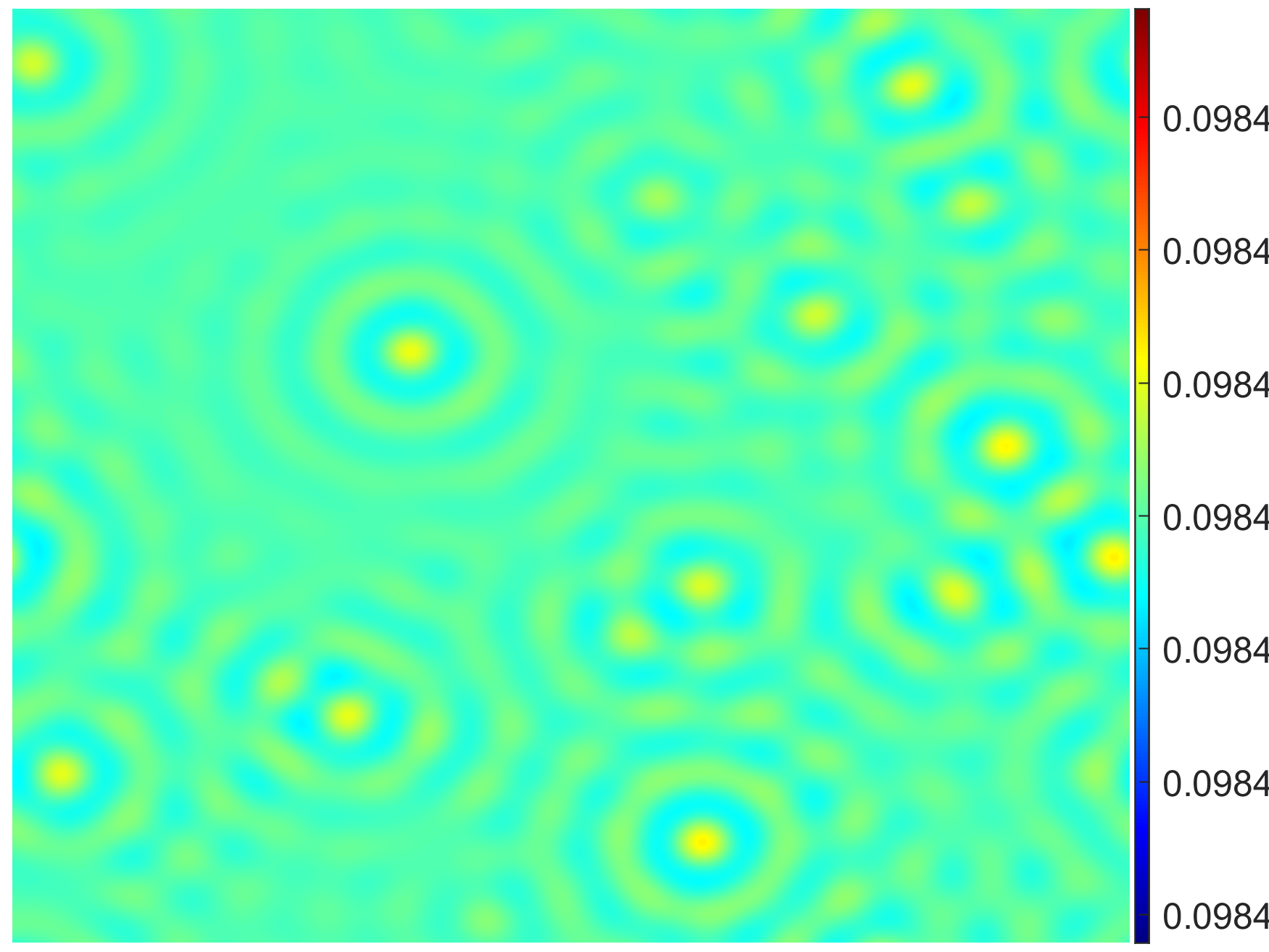} 
				\includegraphics[width=0.25\textwidth]{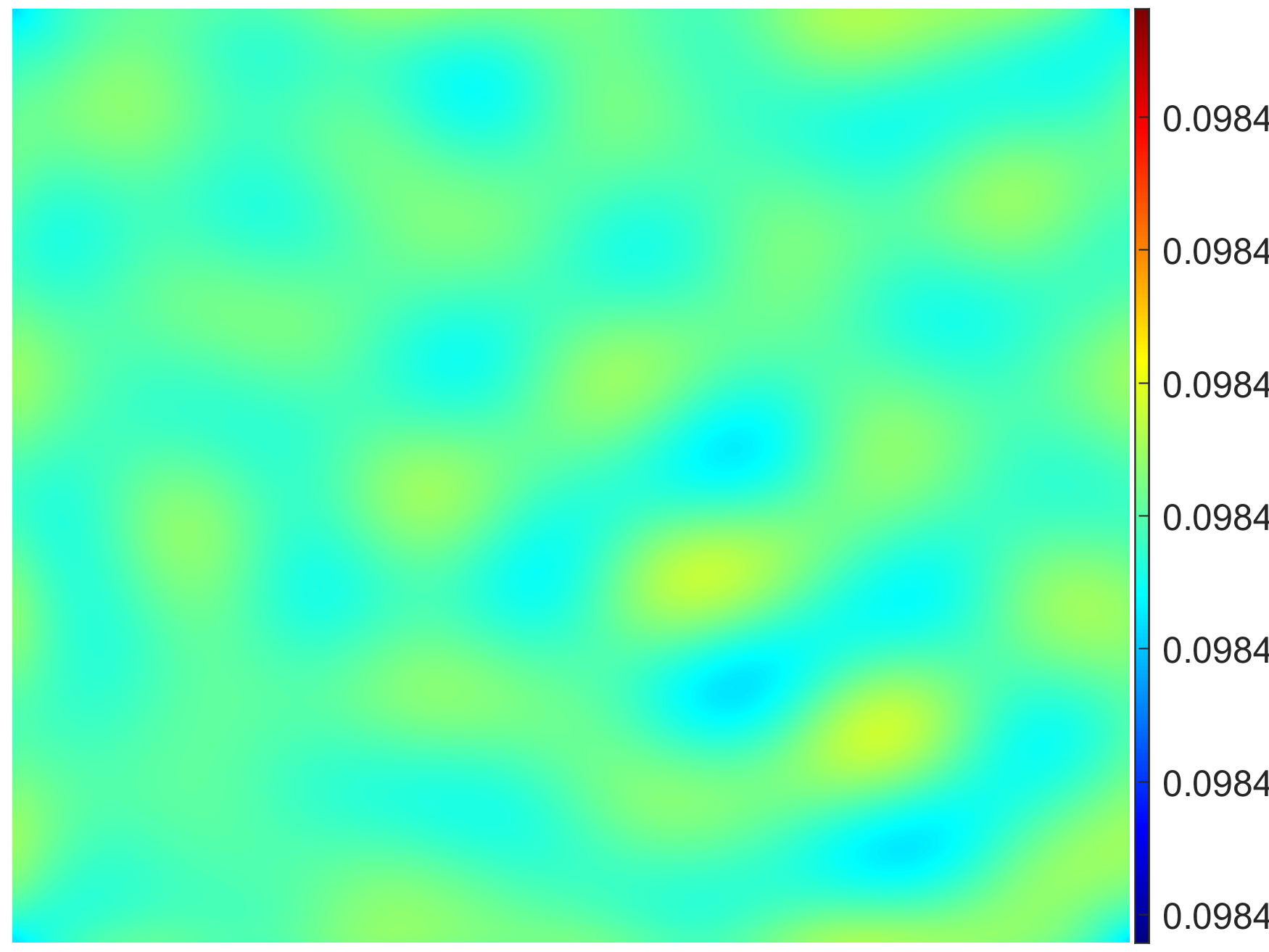} 
				\includegraphics[width=0.25\textwidth]{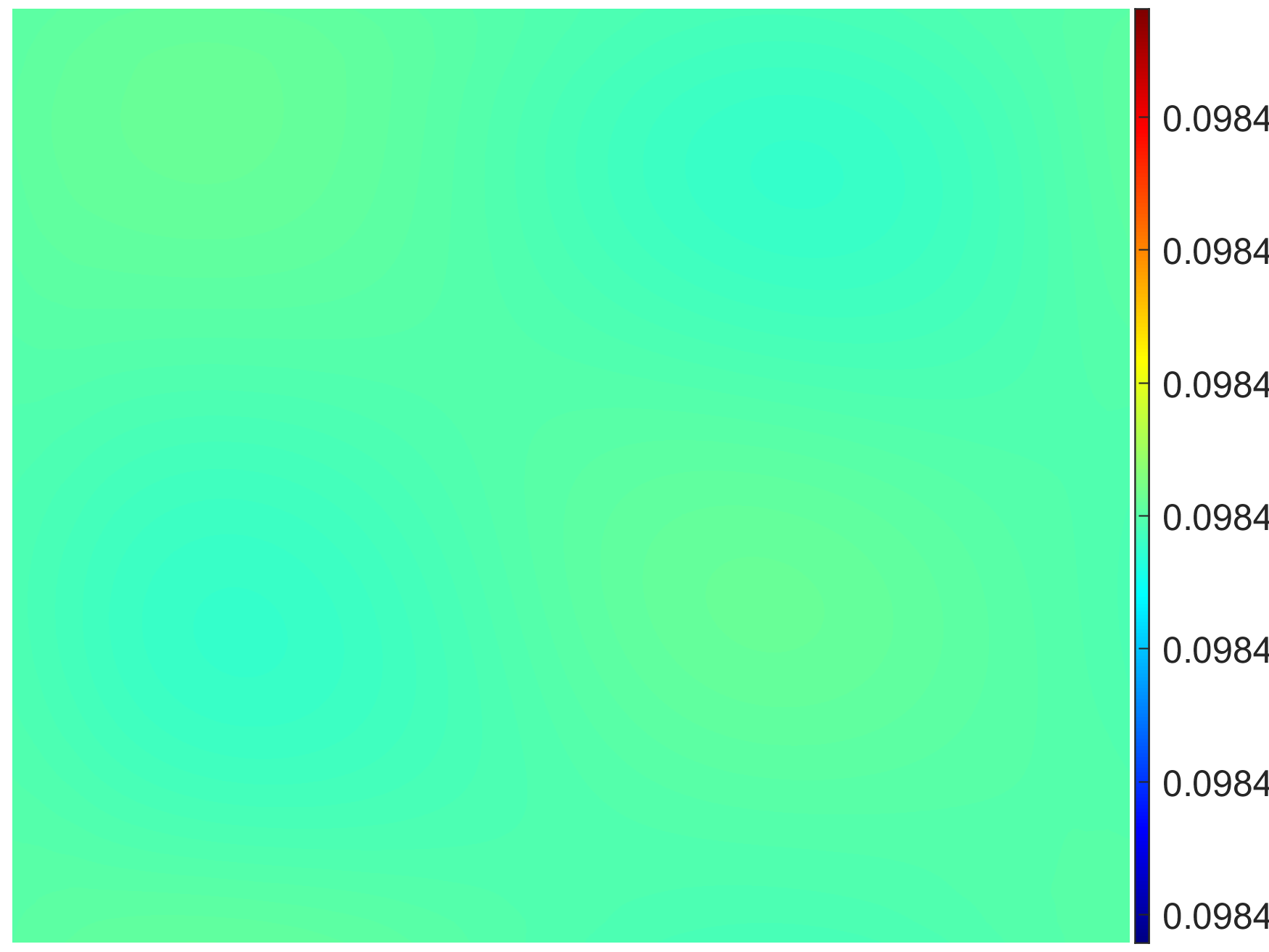}
				\includegraphics[width=0.25\textwidth]{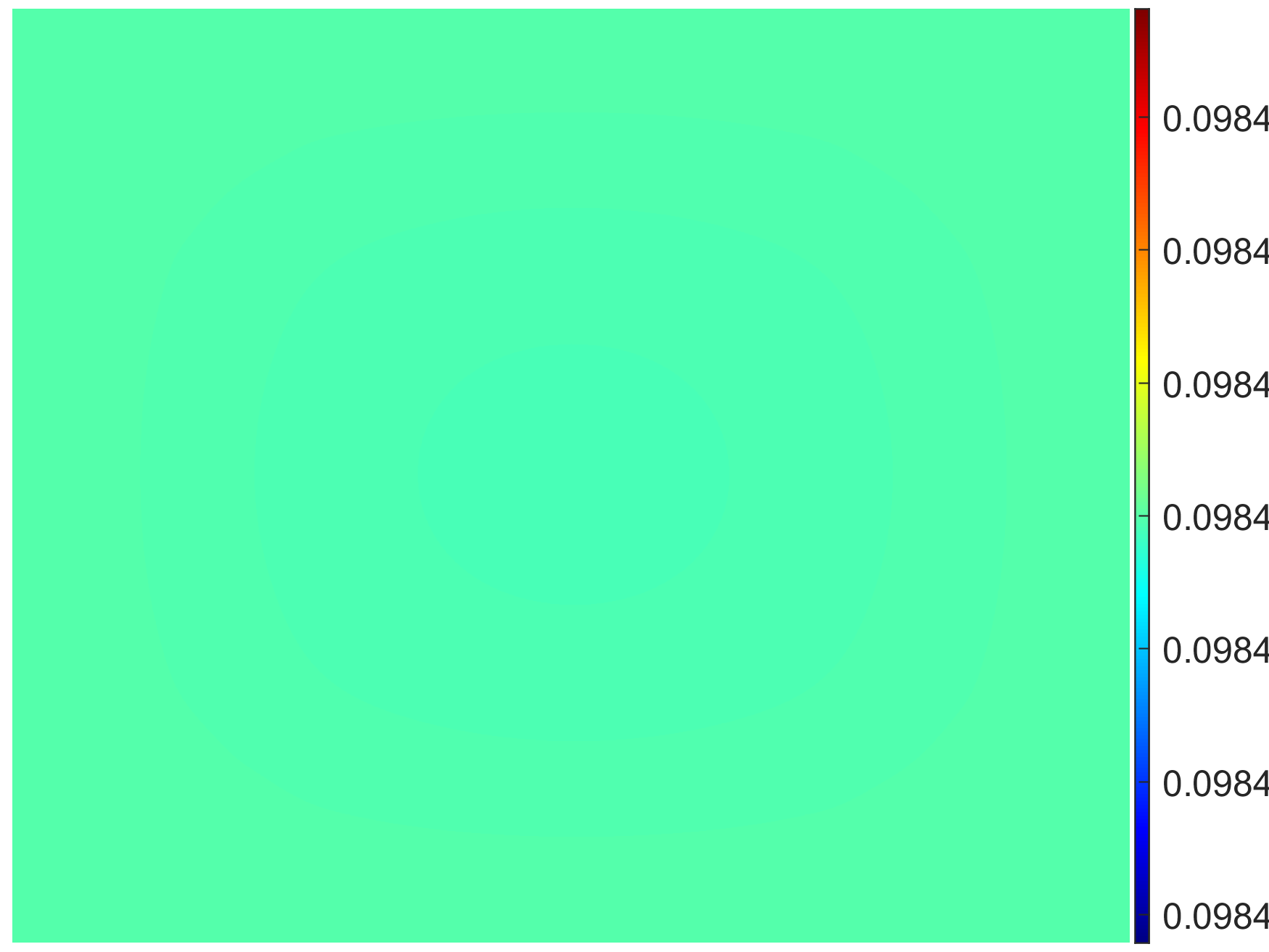} 
			}
			\\
		\subfigure[PDE system (Case B): Crime density $\rho$ at $T = 800$
		for $D = 0.001, 0.01, 0.1, 1$ ($\eta = 0.03$).\label{crime_plot_2}]
			{
				\includegraphics[width=0.25\textwidth]{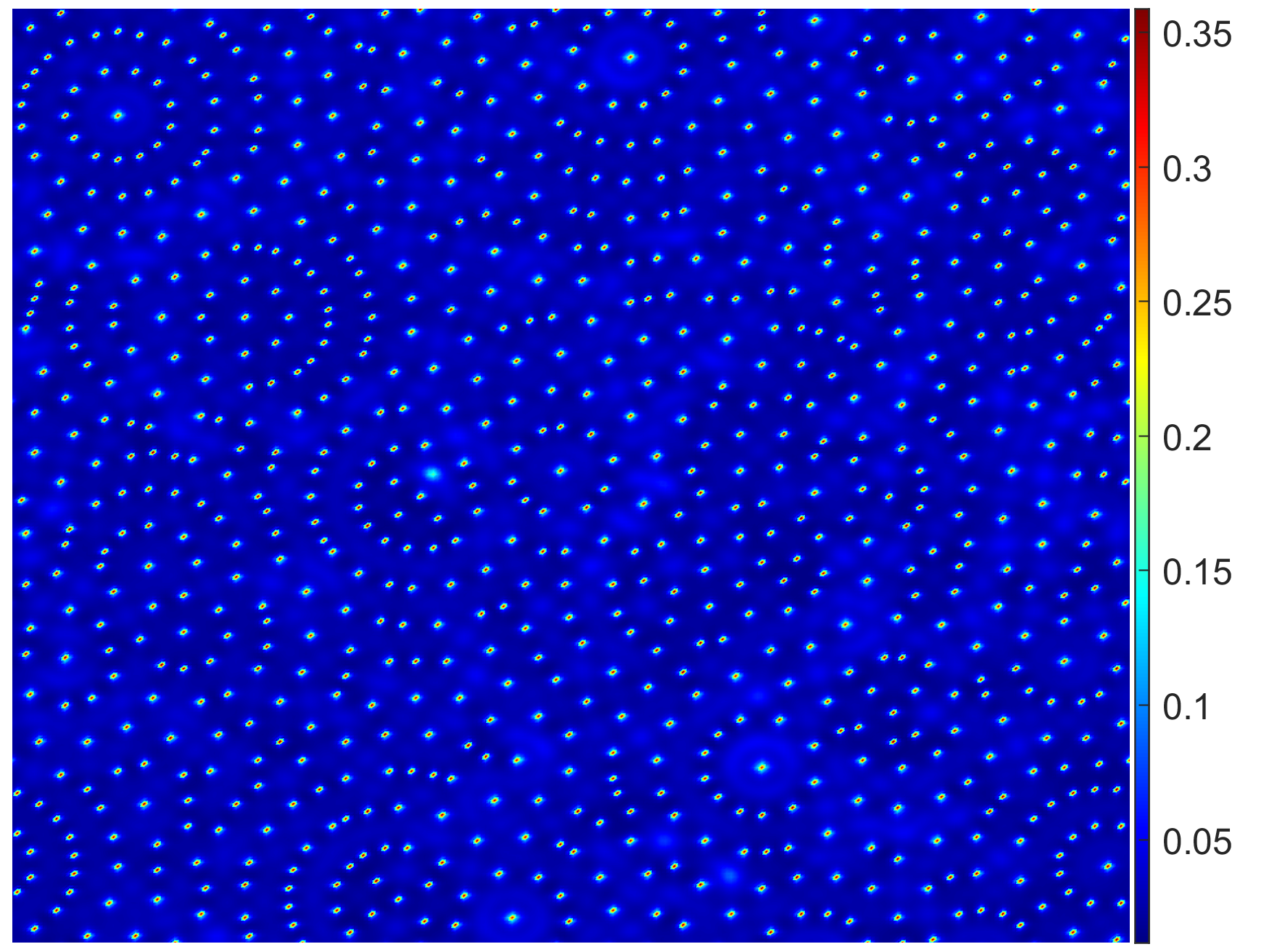}
				\includegraphics[width=0.25\textwidth]{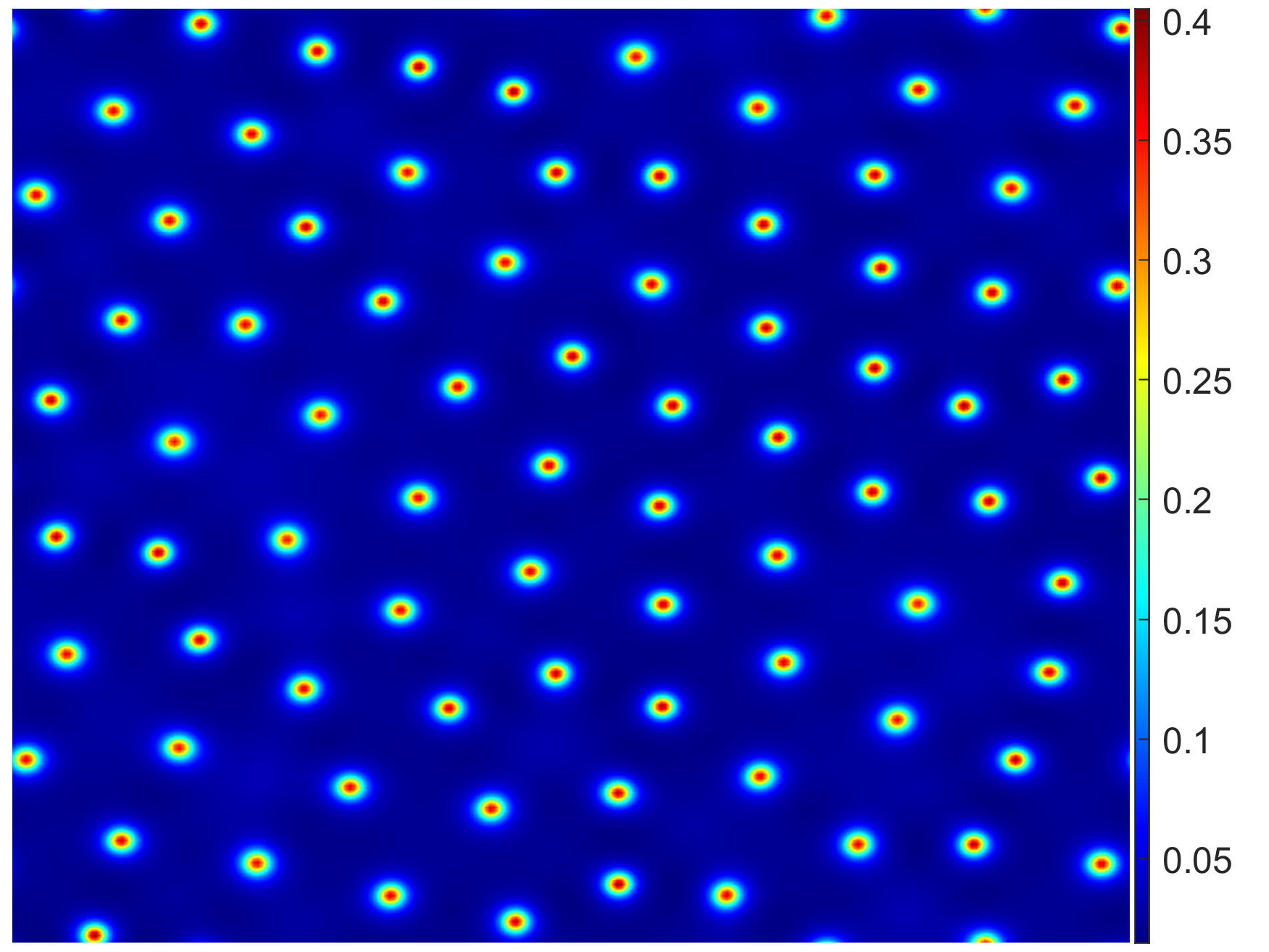} 
				\includegraphics[width=0.25\textwidth]{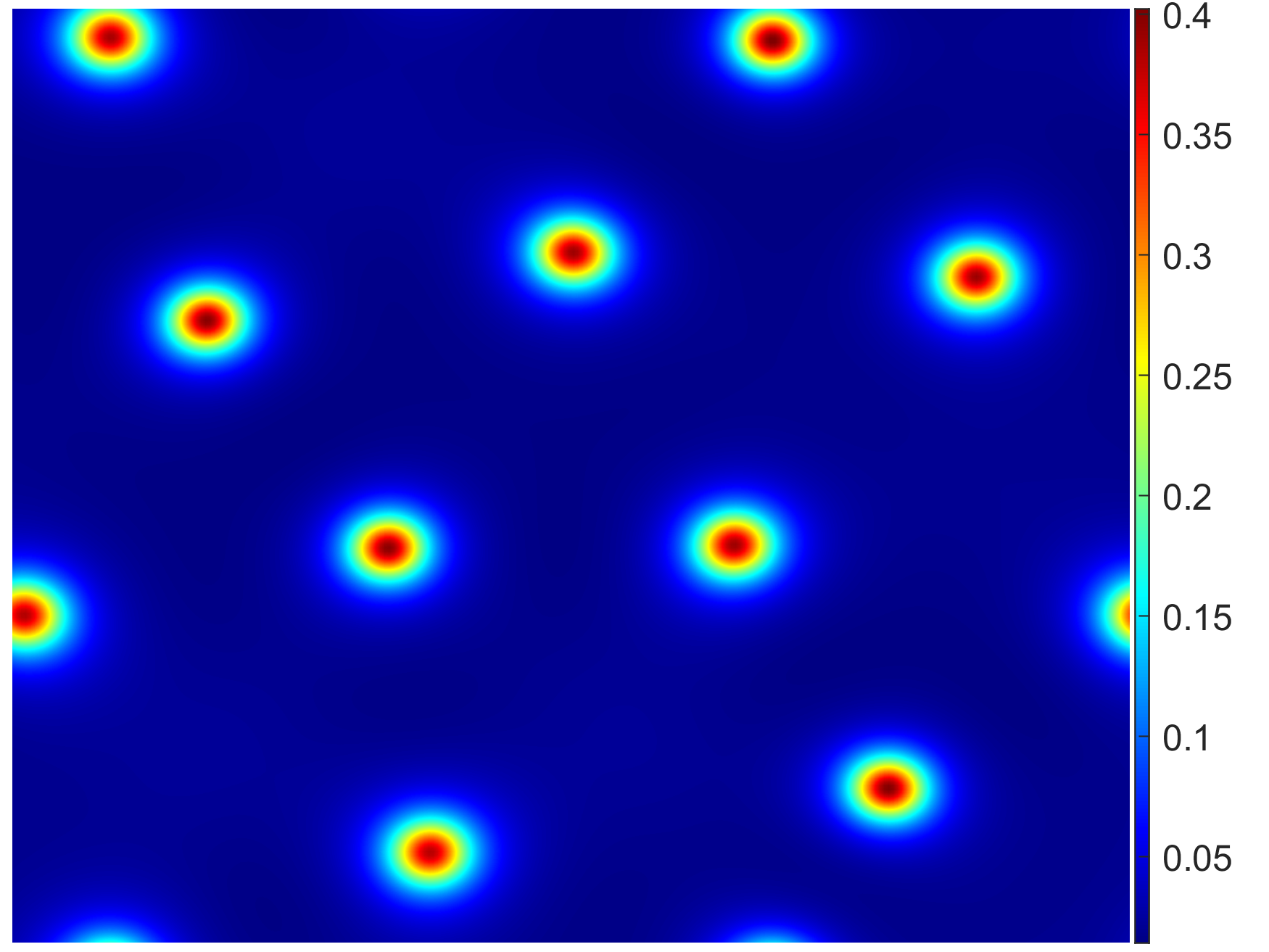}
				\includegraphics[width=0.25\textwidth]{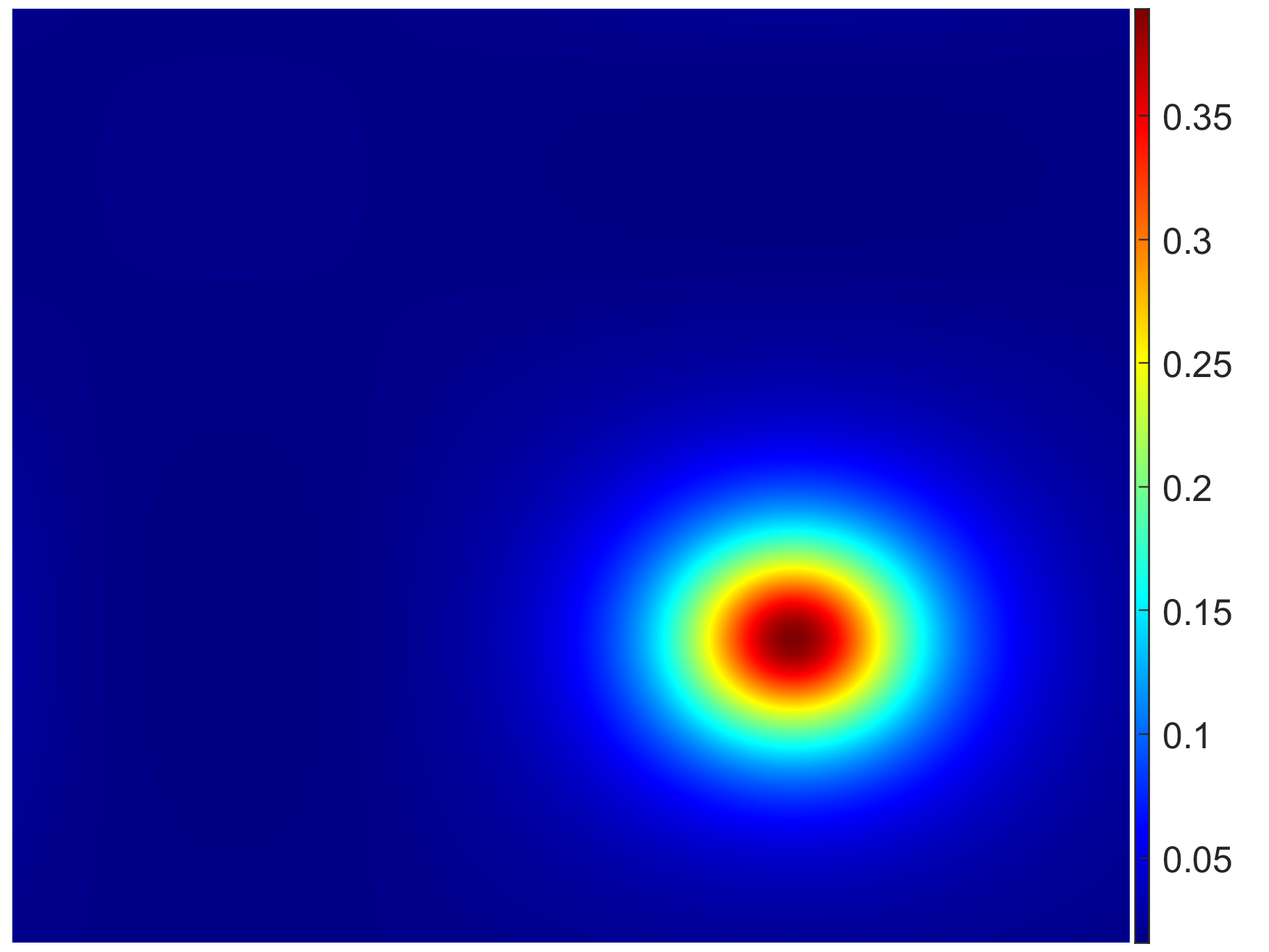}
			}
			\\
			\centering
		\subfigure[Particle system (Case A): Agent distribution at $T = 800$
		for $D = 0.001, 0.01, 0.05, 0.1$ ($\eta = 0.2$).\label{crime_plot_3}]
			{
				\includegraphics[width=0.25\textwidth]{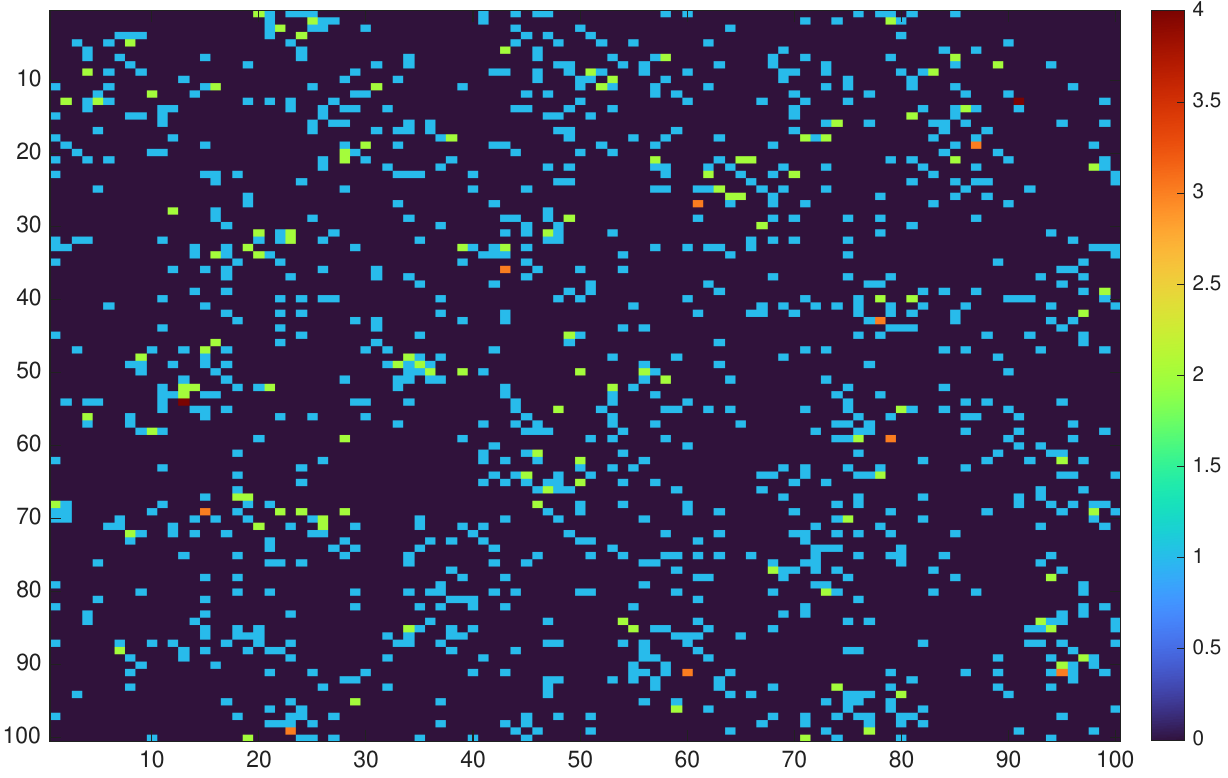} 
				\includegraphics[width=0.25\textwidth]{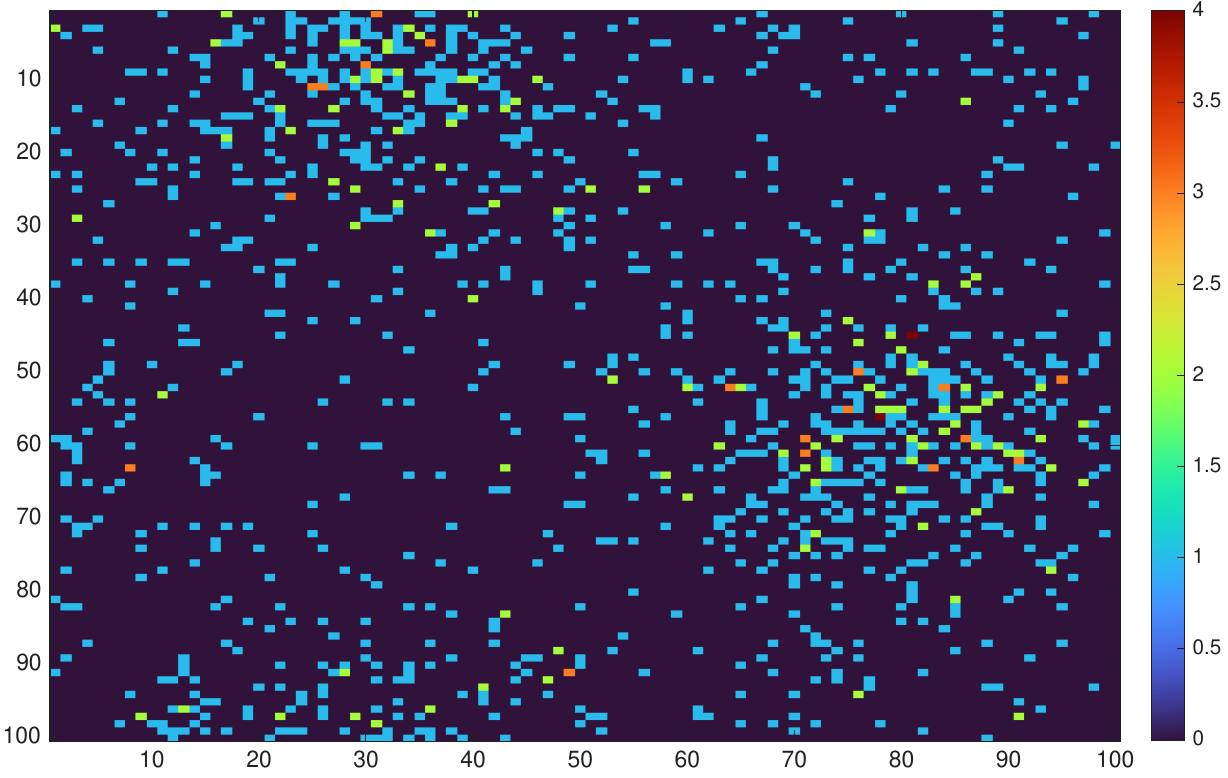} 
				\includegraphics[width=0.25\textwidth]{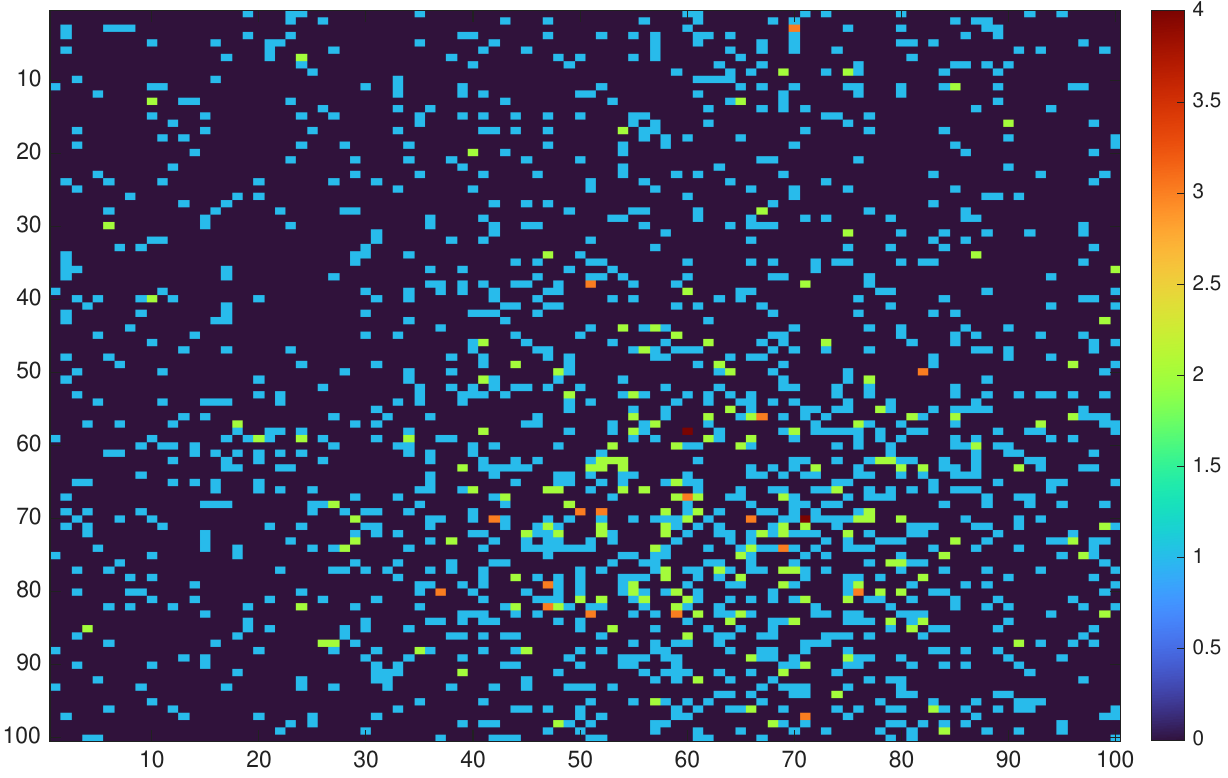}
				\includegraphics[width=0.25\textwidth]{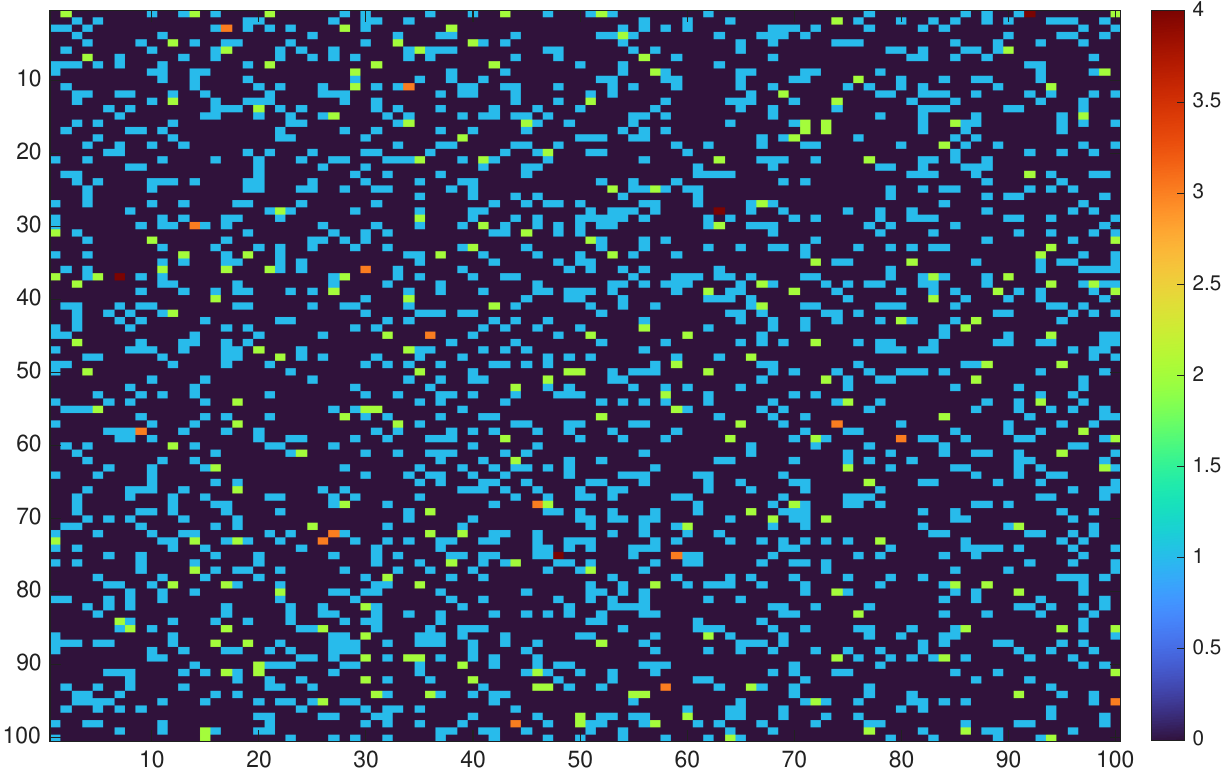} 
			}
			\\
			\centering
		\subfigure[Particle system (Case B): Agent distribution at $T = 800$
		for $D = 0.001, 0.01, 0.05, 0.1$ ($\eta = 0.03$).\label{crime_plot_4}]
			{
				\includegraphics[width=0.25\textwidth]{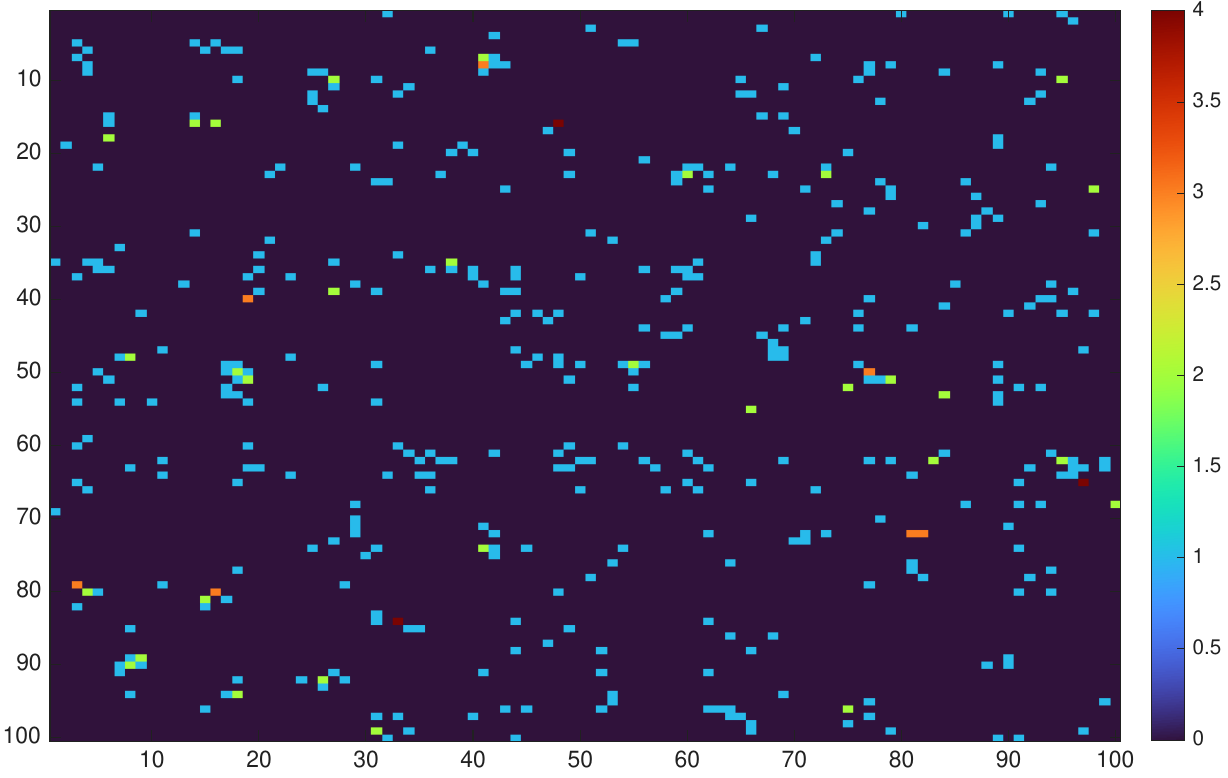} 
				\includegraphics[width=0.25\textwidth]{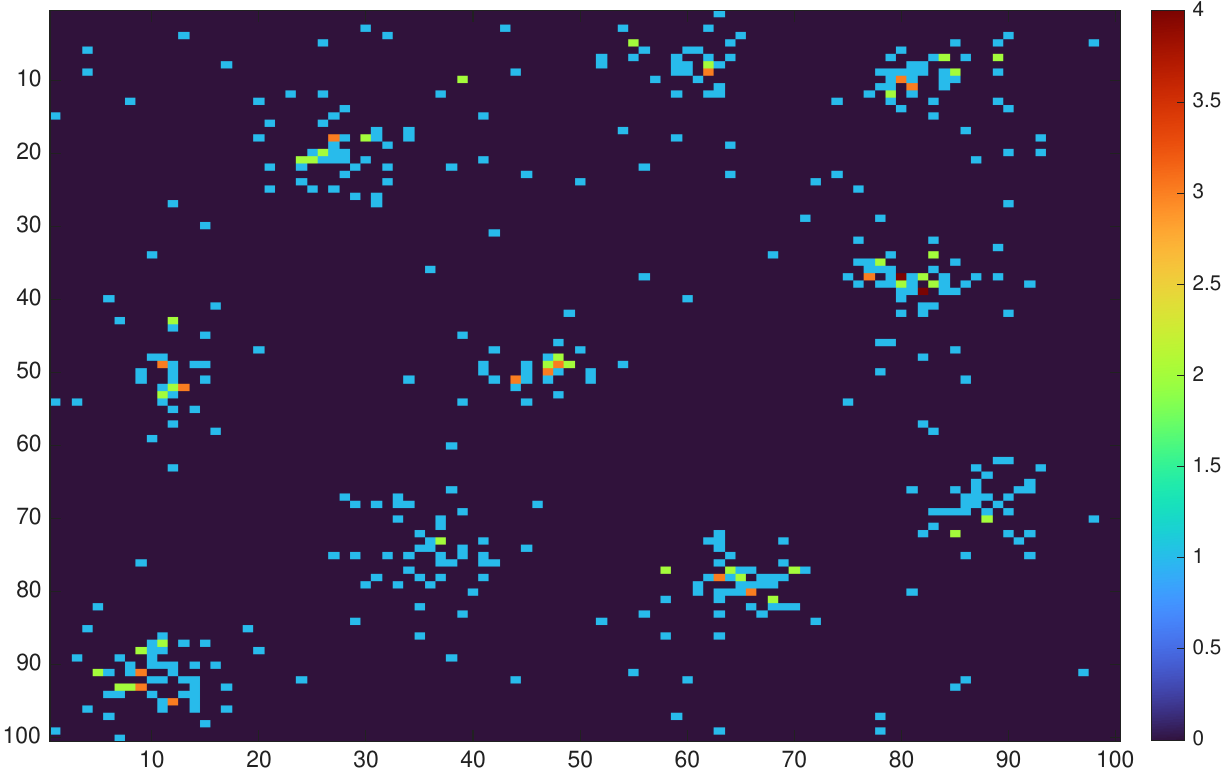} 
				\includegraphics[width=0.25\textwidth]{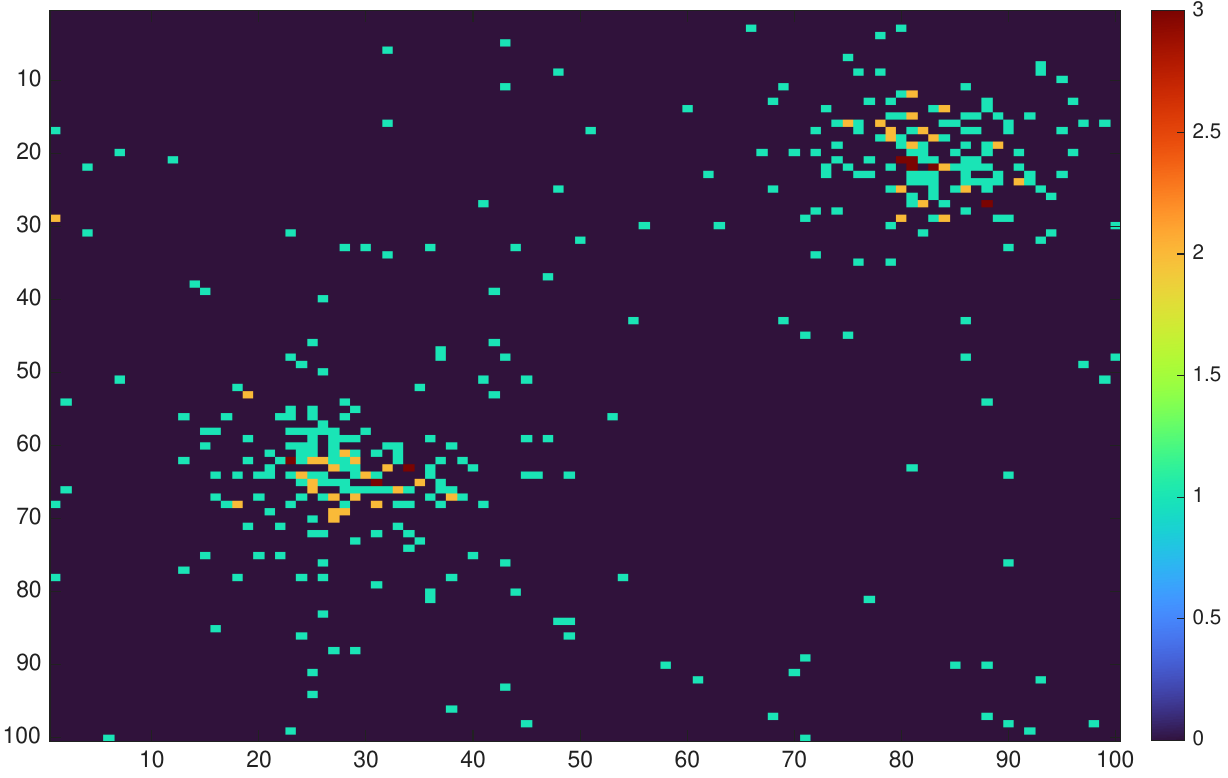}
				\includegraphics[width=0.25\textwidth]{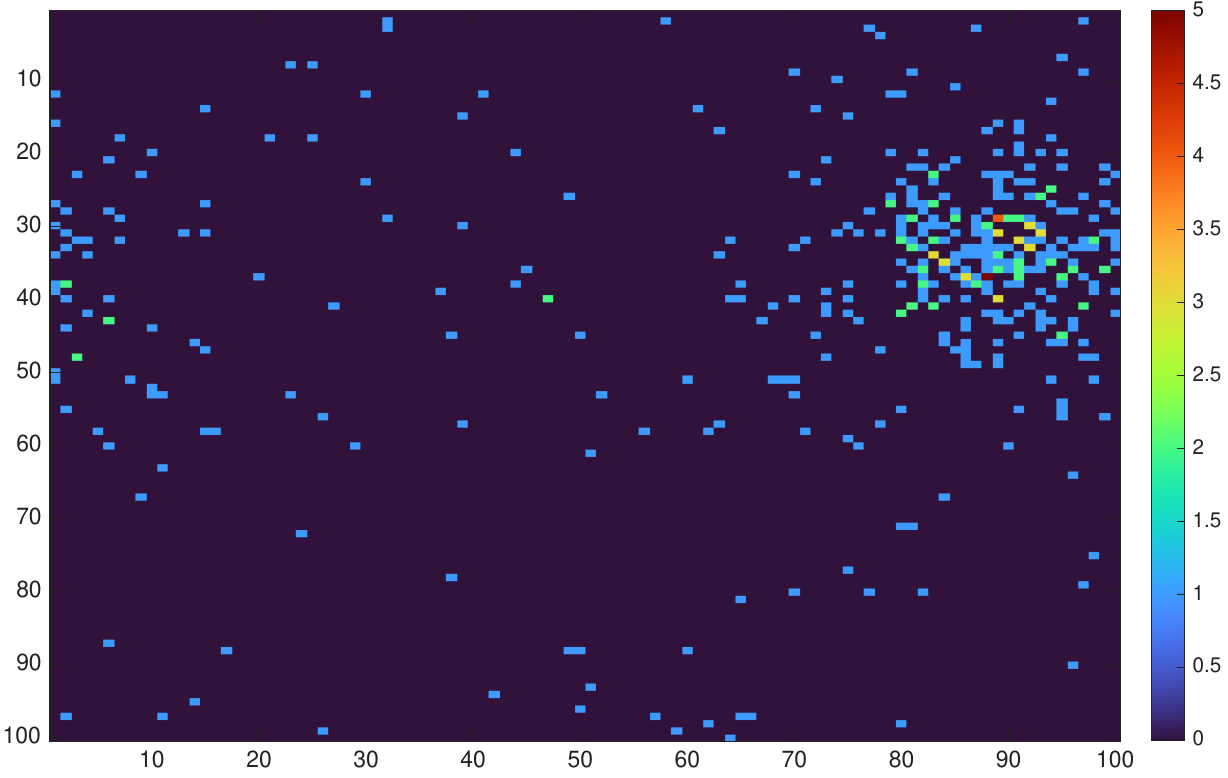} 
			}
			\vspace{-10pt}
			\caption{Population density at $T = 800$ day in PDE-based and particle-based crime modelings: A similar pattern is observed in both particle system and PDE system. From the horizontal view, the hotspots become more concentrated as the diffusion coefficient $D$ increases. When the value of $D$ is small, the density of the crime population shows a local clustering pattern such as "points" or "rings". In contrast, when the value of $D$ is large, the density distribution exhibits several large spikes. From the vertical view,  when the diffusion coefficient $D$ is fixed, the parameter $\eta$ has a qualitative impact on the pattern of distribution of crime. When $\eta$ decreases (from $\eta = 0.2$ to $\eta = 0.03$), the clustered 'points' become more regular, and the phase diagram transits from diffusion-dominated to aggregation-dominated.}
			\label{crime_plot}
		\end{figure}

		According to the phase diagram in Figure \ref{crime_phase}, the first case
		lies in the stable regime, while the second lies in the unstable regime.
		Figures \ref{crime_plot_1} and \ref{crime_plot_2} show the evolution of the
		PDE crime model at $T=800$, while
		Figures \ref{crime_plot_3} and \ref{crime_plot_4} display the corresponding
		particle simulations from
		\cite{ShortDorsognaPasourTitaBrantinghamBertozziChayes2008}.
		The two models exhibit qualitatively similar patterns.
		
		When the diffusion coefficient $D$ increases, the hotspots become more
		concentrated and may merge into larger spikes.
		For fixed $D$, decreasing $\eta$ (from $0.2$ to $0.03$)
		produces more regular hotspot patterns, indicating a transition from
		diffusion-dominated behavior to aggregation-dominated behavior.
		These observations are consistent with the theoretical prediction in
		Figure \ref{crime_phase}.

	Now we  investigate how model parameters
	influence the spatial characteristics of hotspot patterns.
	In particular, we examine how the neighbourhood effect parameter $\eta$
	influences the characteristic scale of hotspot patterns,
	following the approach in \cite{HaoMilyQuainiZhong2026}.
			
\begin{figure}[h!]
			\centering
			\includegraphics[width=0.65\textwidth]{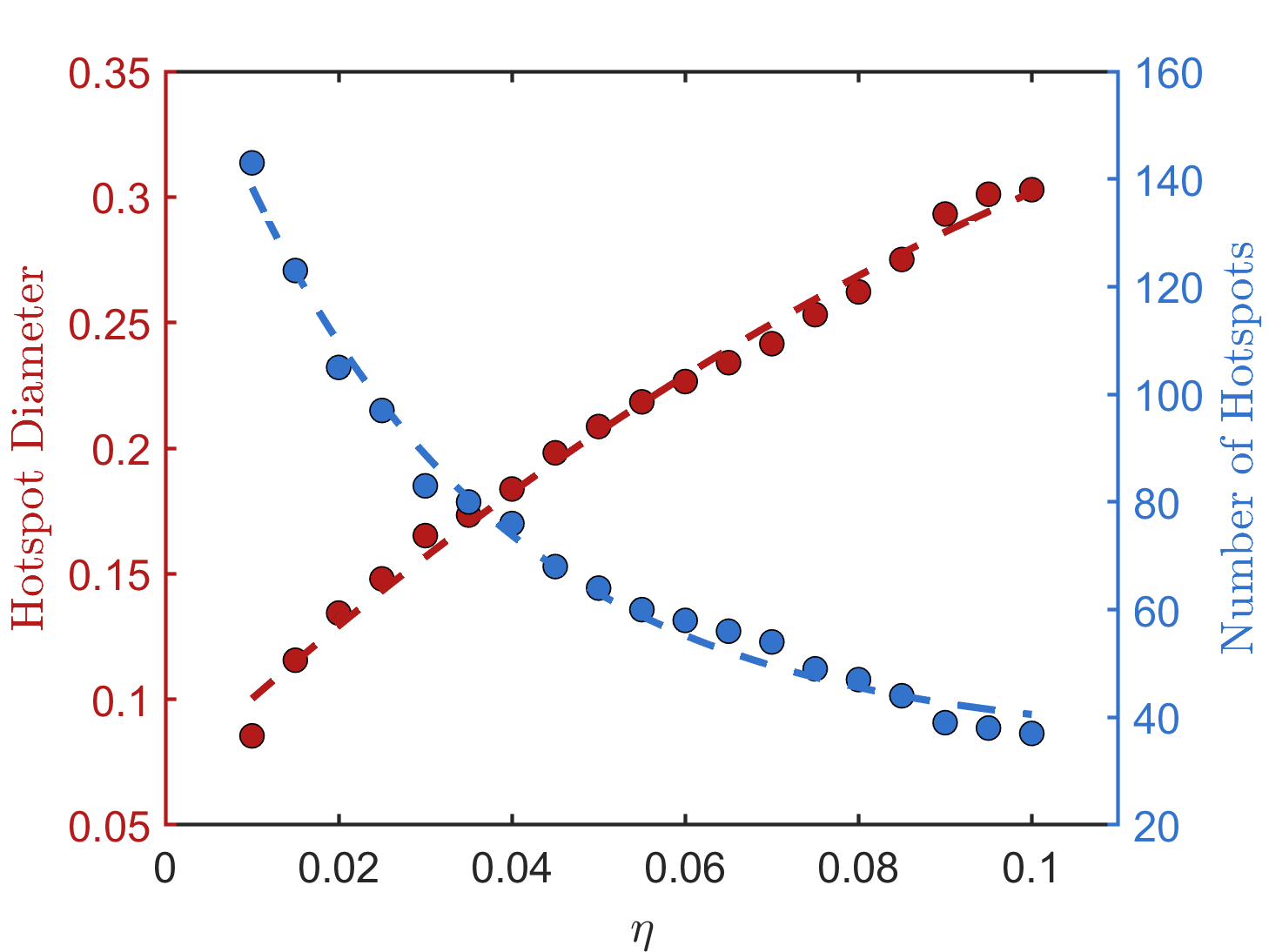}
			\vspace{-10pt}
			\caption{(Color online) Hotspot diameter (left axis, red markers)
				and number of hotspots (right axis, blue markers)
				as the neighbourhood effect parameter $\eta$ varies in
				$[0.01,0.1]$. The dashed curves show the corresponding
				empirical fitting results.}
			\label{fig:hotspots_fit}
		\end{figure}

		Figure~\ref{fig:hotspots_fit} shows the hotspot diameter
		(left axis) and the number of hotspots (right axis)
		for $\eta\in[0.01,0.1]$.
		The results reveal a clear dependence of the spatial pattern on $\eta$.
		As $\eta$ decreases, the diffusion of the attractiveness field becomes
		weaker, which enhances the aggregation mechanism and leads to smaller
		and more numerous hotspots.
		To quantify this trend, we perform an empirical curve fitting
		based on the computed data.
		The number of hotspots is well approximated by the exponential relation
		$143.40\exp(-32.73\eta)+35.10$,
		while the hotspot diameter follows approximately the quadratic law
		$-8.30\eta^2+3.15\eta+0.07$.
		These relations provide a phenomenological description of the numerical
		results, although a rigorous analytical derivation remains further investigations.
	
\begin{figure}[h!]
			\centering
			\subfigure
			{
				\includegraphics[width=0.25\textwidth]{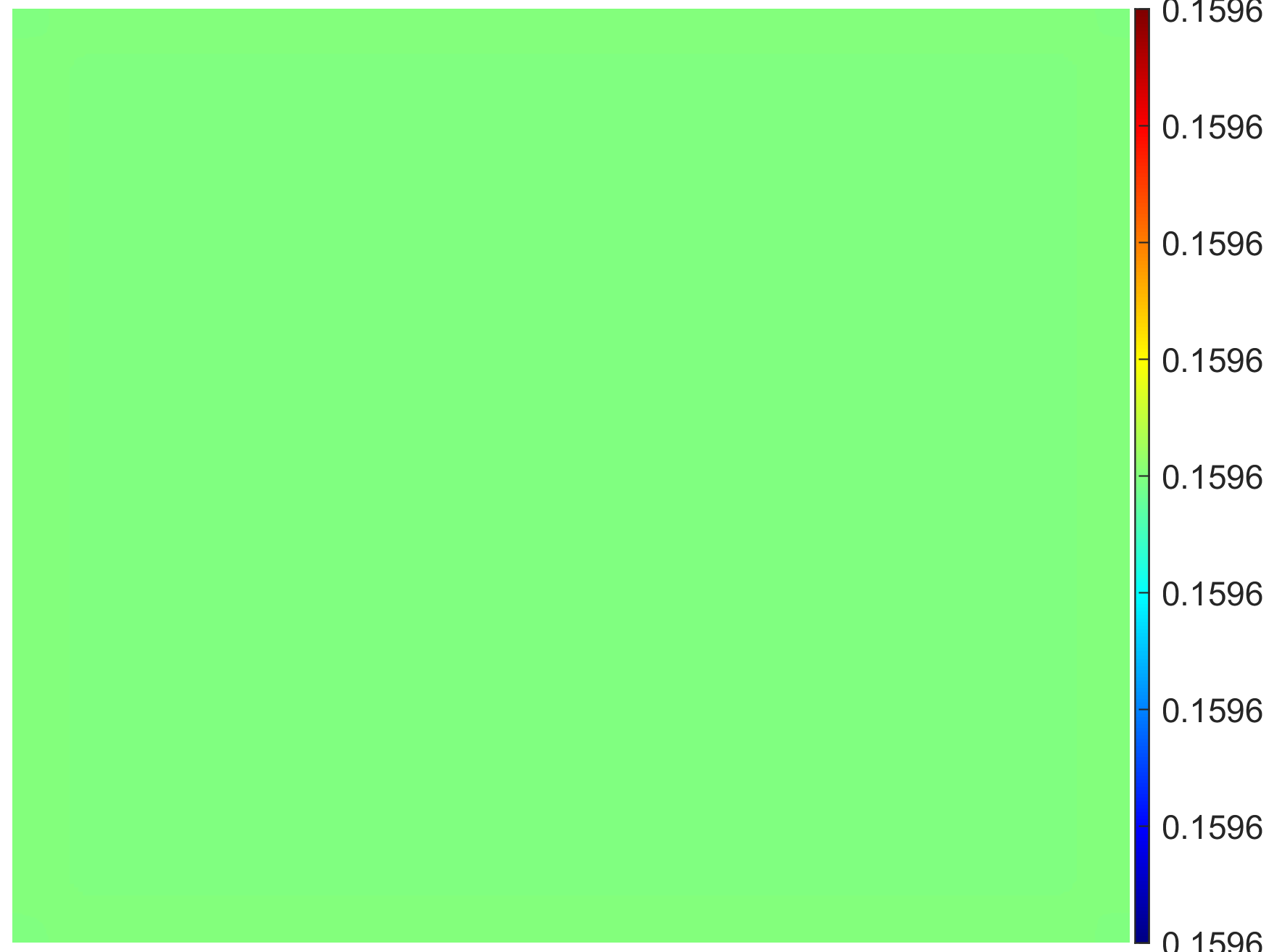}
				\includegraphics[width=0.25\textwidth]{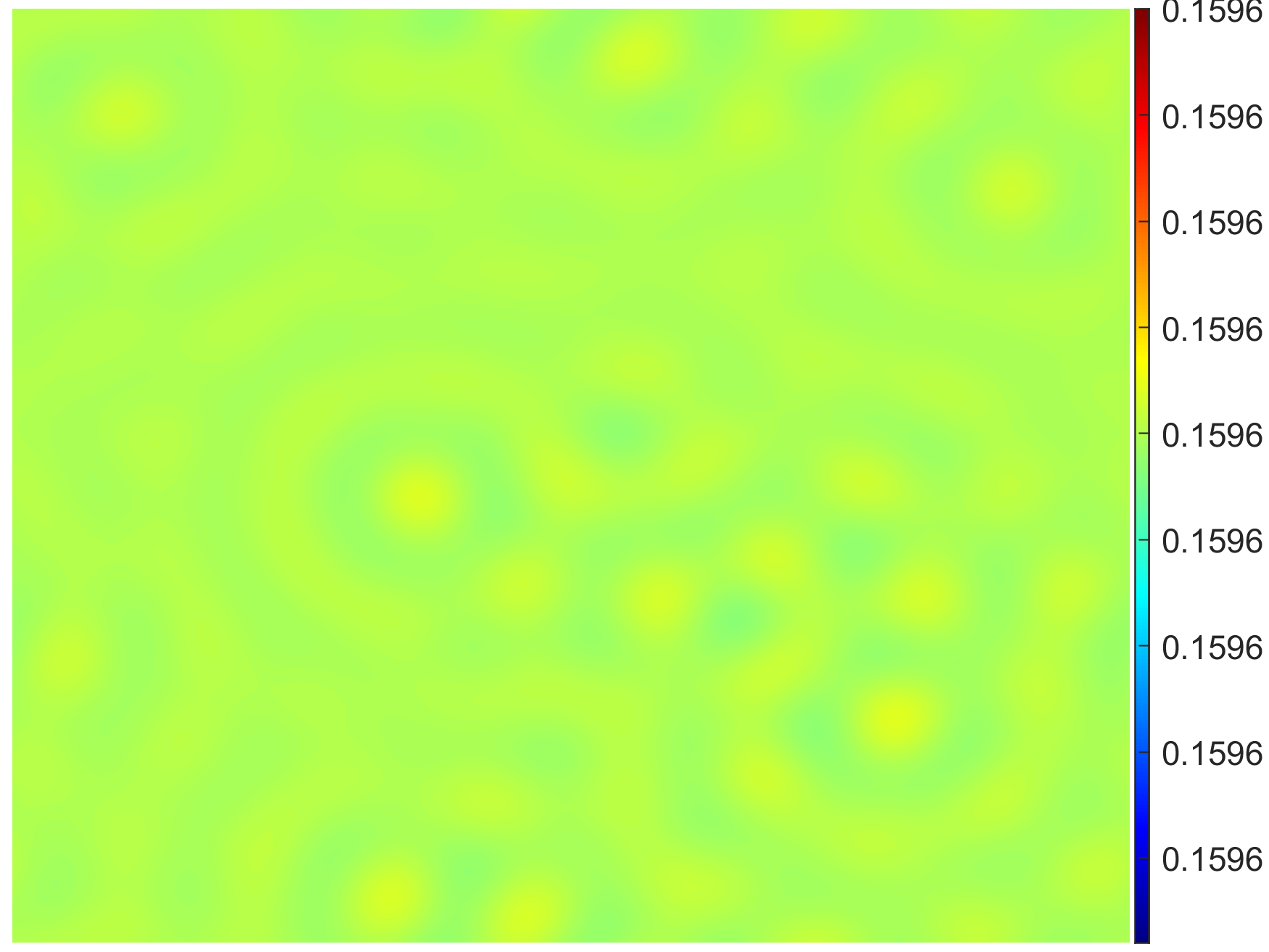} 
				\includegraphics[width=0.25\textwidth]{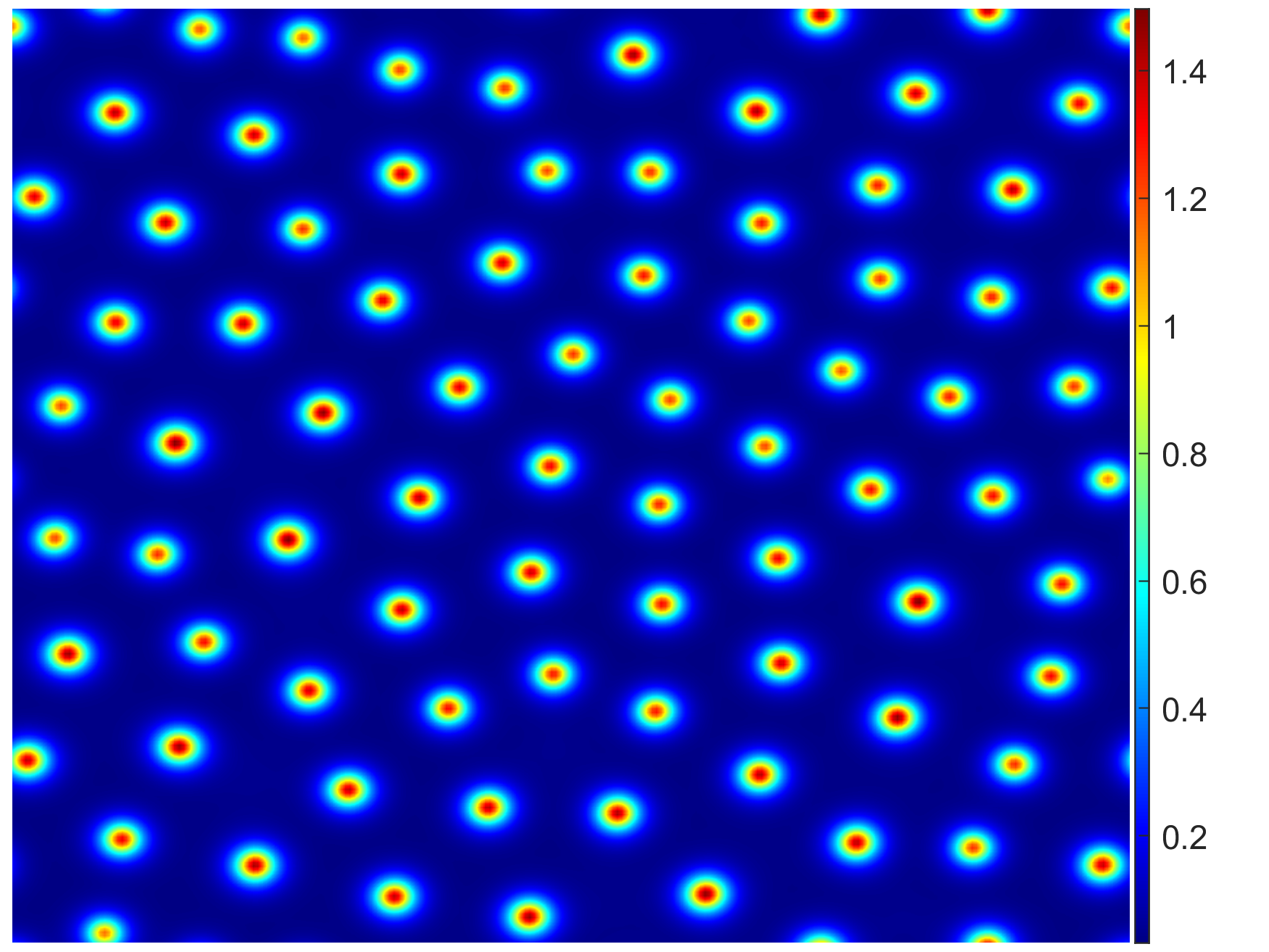}
				\includegraphics[width=0.25\textwidth]{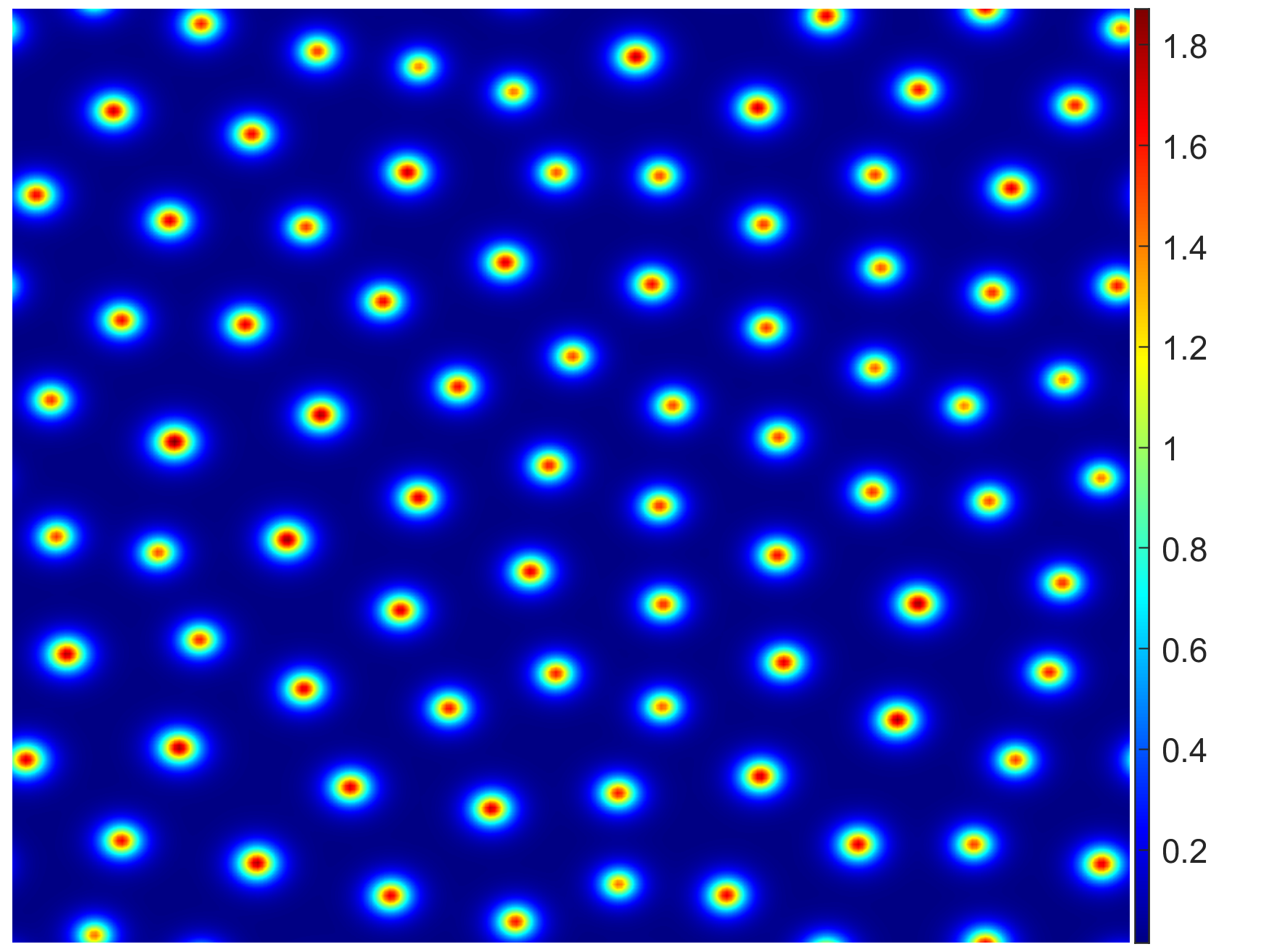}
			}\\
			\vspace{-5pt}
			\subfigure
			{
				\includegraphics[width=0.25\textwidth]{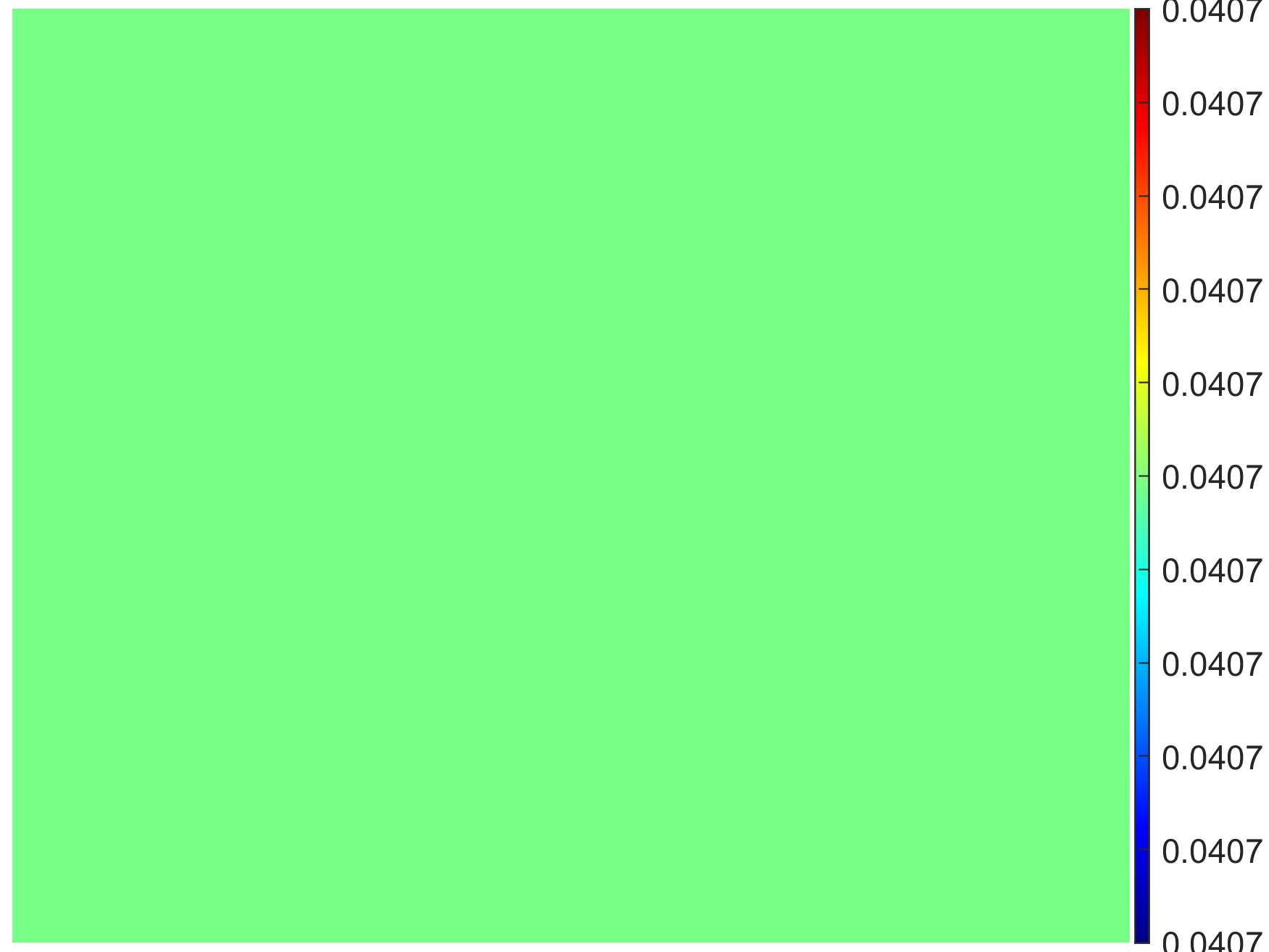} 
				\includegraphics[width=0.25\textwidth]{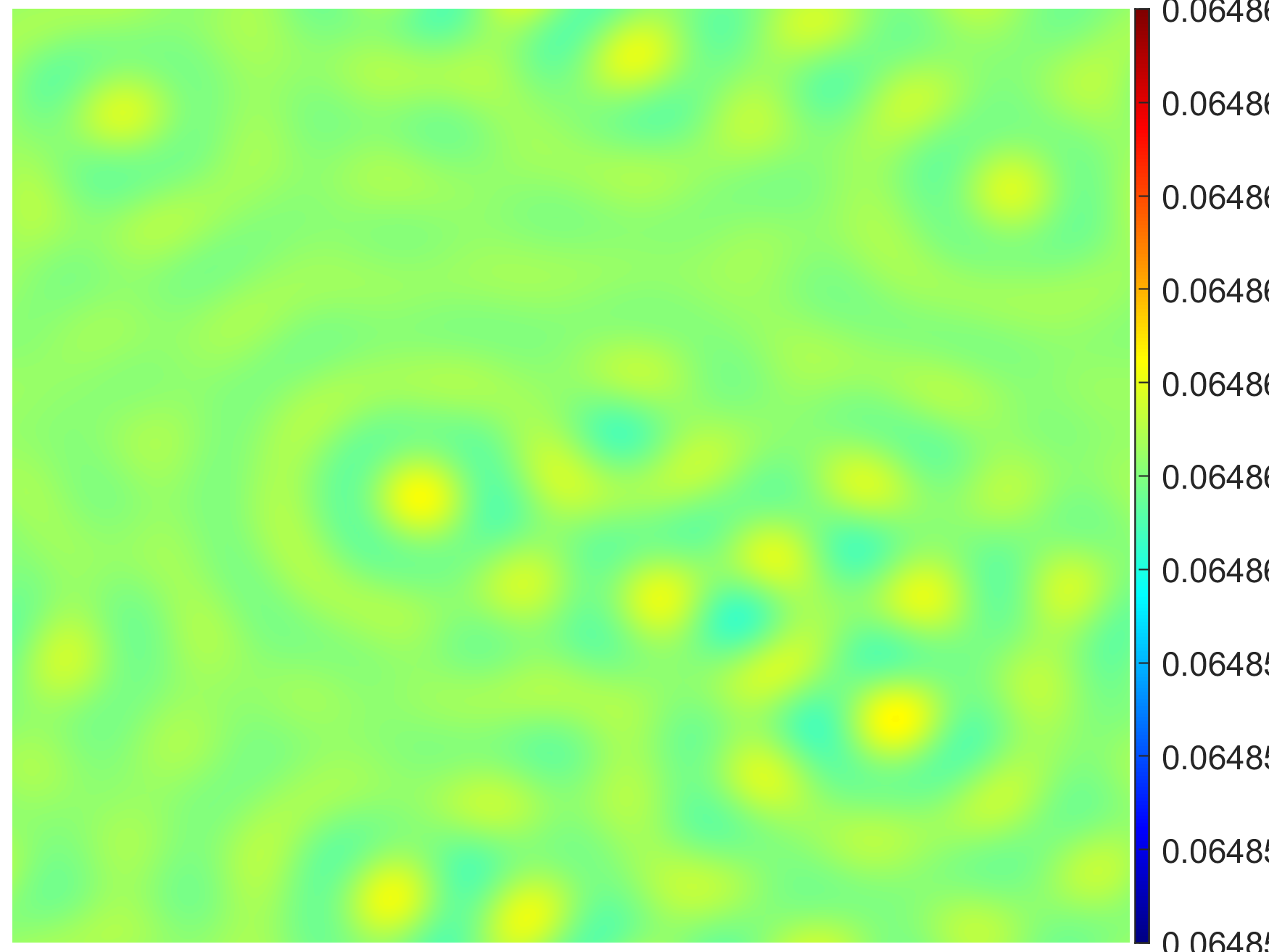} 
				\includegraphics[width=0.25\textwidth]{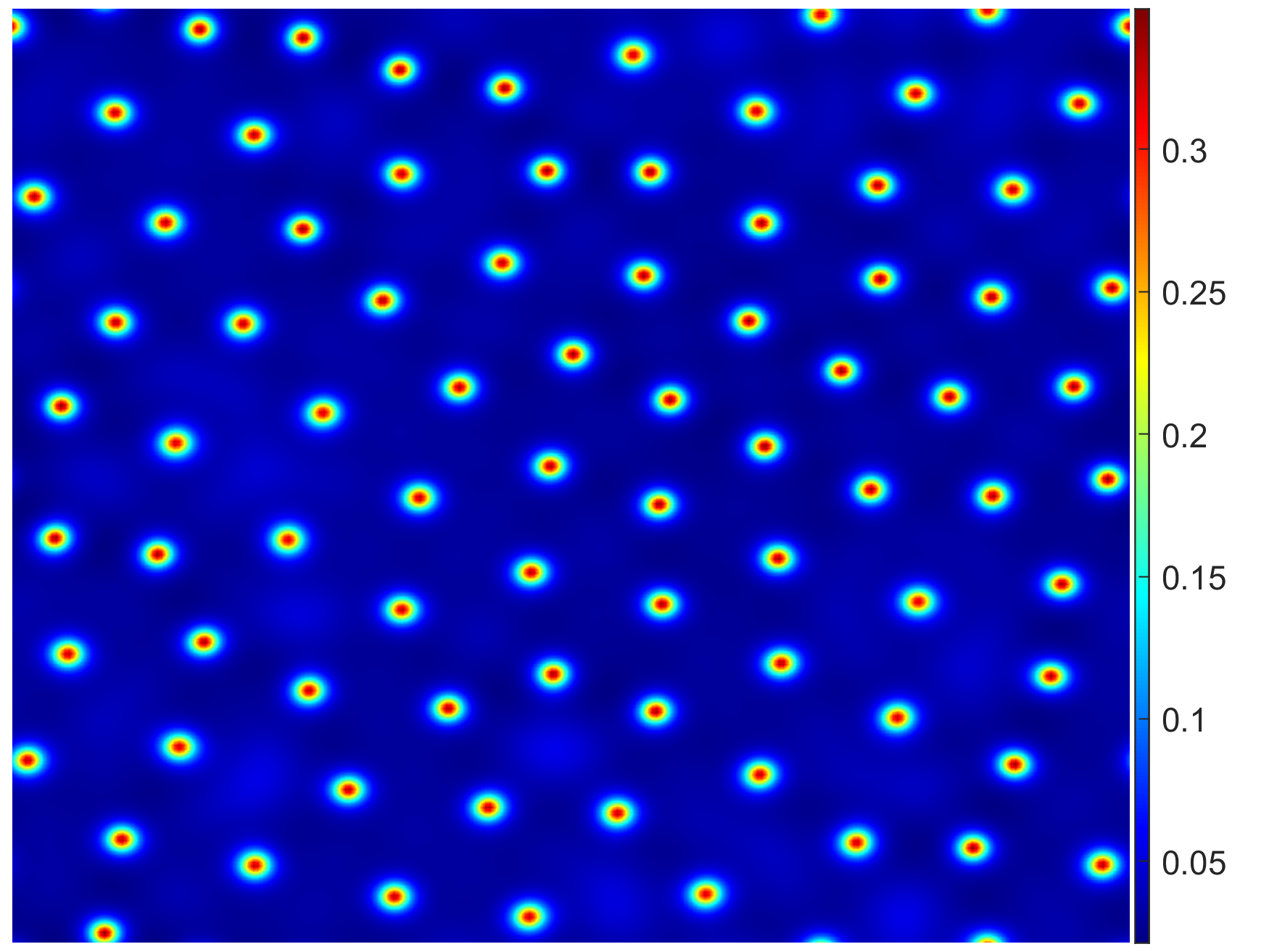}
				\includegraphics[width=0.25\textwidth]{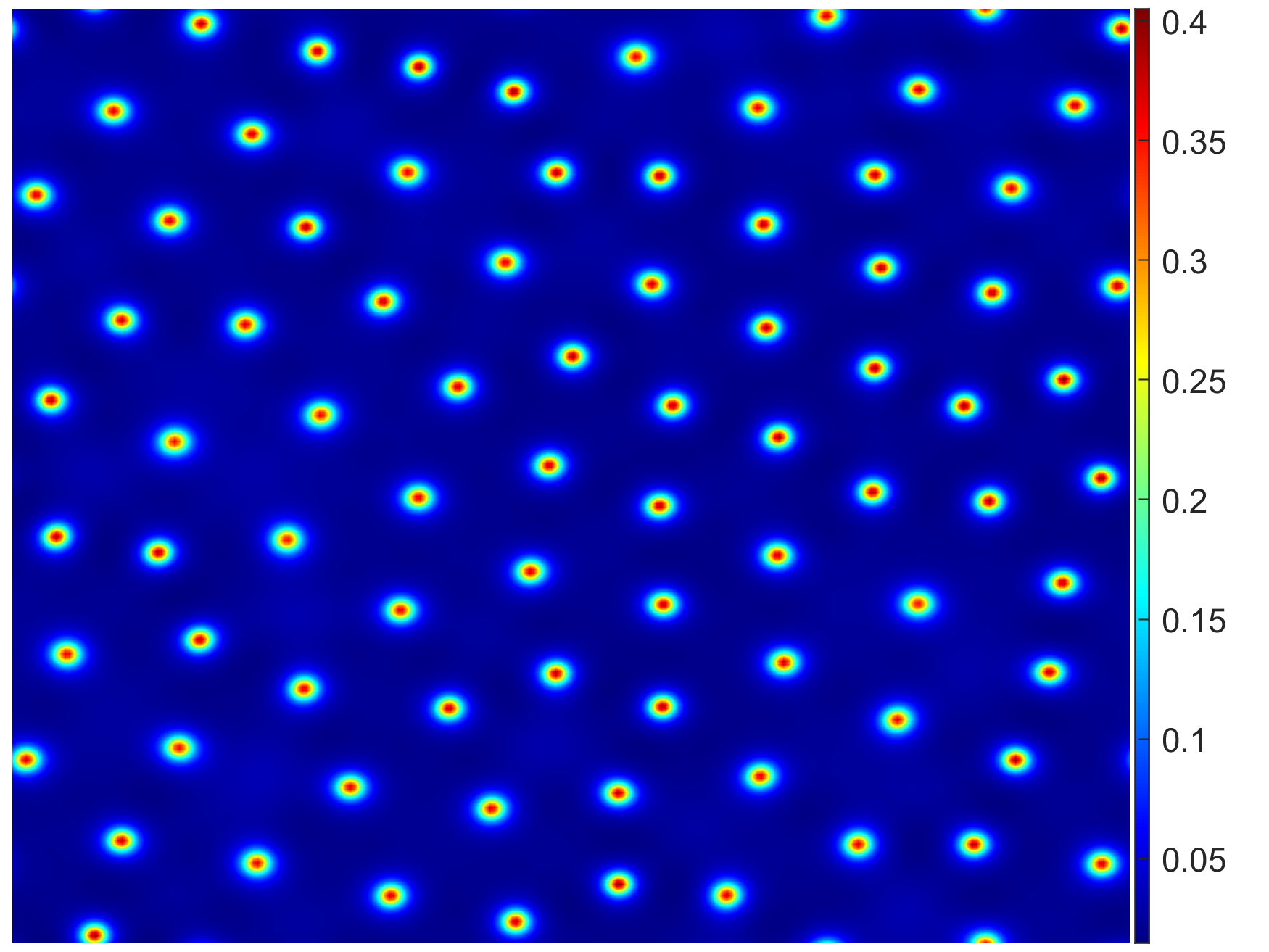} 
			}
			\vspace{-15pt}
			\caption{Snapshots of the solutions for different values of the baseline
				attractiveness $A_0$: evolution of the attractiveness field $A$ (top)
				and criminal density $\rho$ (bottom) at $T=800$.
				From left to right,
				$A_0 = 1/3.75,\;1/7.5,\;1/15,\;1/30$.
			}
			\label{fig_A0}
		\end{figure}

		Figure~\ref{fig_A0} shows the final states of the system
		at time $T=800$ under different values of $A_0$.
		When $A_0$ is relatively large ($A_0=1/3.75$ and $1/7.5$),
		the system remains close to a spatially homogeneous state and
		no hotspot structures are observed.
		As $A_0$ decreases ($A_0=1/15$ and $1/30$),
		localized peaks begin to emerge and eventually form
		regular hotspot patterns.
		This indicates that a smaller baseline attractiveness
		facilitates the self-organization of crime hotspots.
	
		Next, we investigate the cases where the neighborhood effect varies
		spatially, which may reflect the self-aggregation in a multi-scale level.

		\begin{example}
		We consider a spatially heterogeneous neighborhood effect
		by prescribing $\eta=\eta(x, y)$ piecewisely in the $x$-direction:
		\begin{equation}\label{eta_change_equation}
			\eta(x, y)=
			\left\{
			\begin{aligned}
				&0.03 \quad \textup{(as in Case B)}, &&1.5\pi \le x \le 2\pi,\\
				&0.06, &&\pi \le x < 1.5\pi,\\
				&0.1,  &&0.5\pi \le x < \pi,\\
				&0.2 \quad \textup{(as in Case A)}, &&0 \le x < 0.5\pi,
			\end{aligned}
			\right.
		\end{equation}
		with the same initial condition in Eq.~\eqref{init_condition}.
		This setting divides the domain into four vertical subregions
		with different values of $\eta$, allowing us to examine how spatial
		heterogeneity in the neighborhood effect influences the scale
		of hotspot formation. 
		\end{example}
		
\begin{figure}[h!]
			\centering
			\subfigure
			{
				\includegraphics[width=0.25\textwidth]{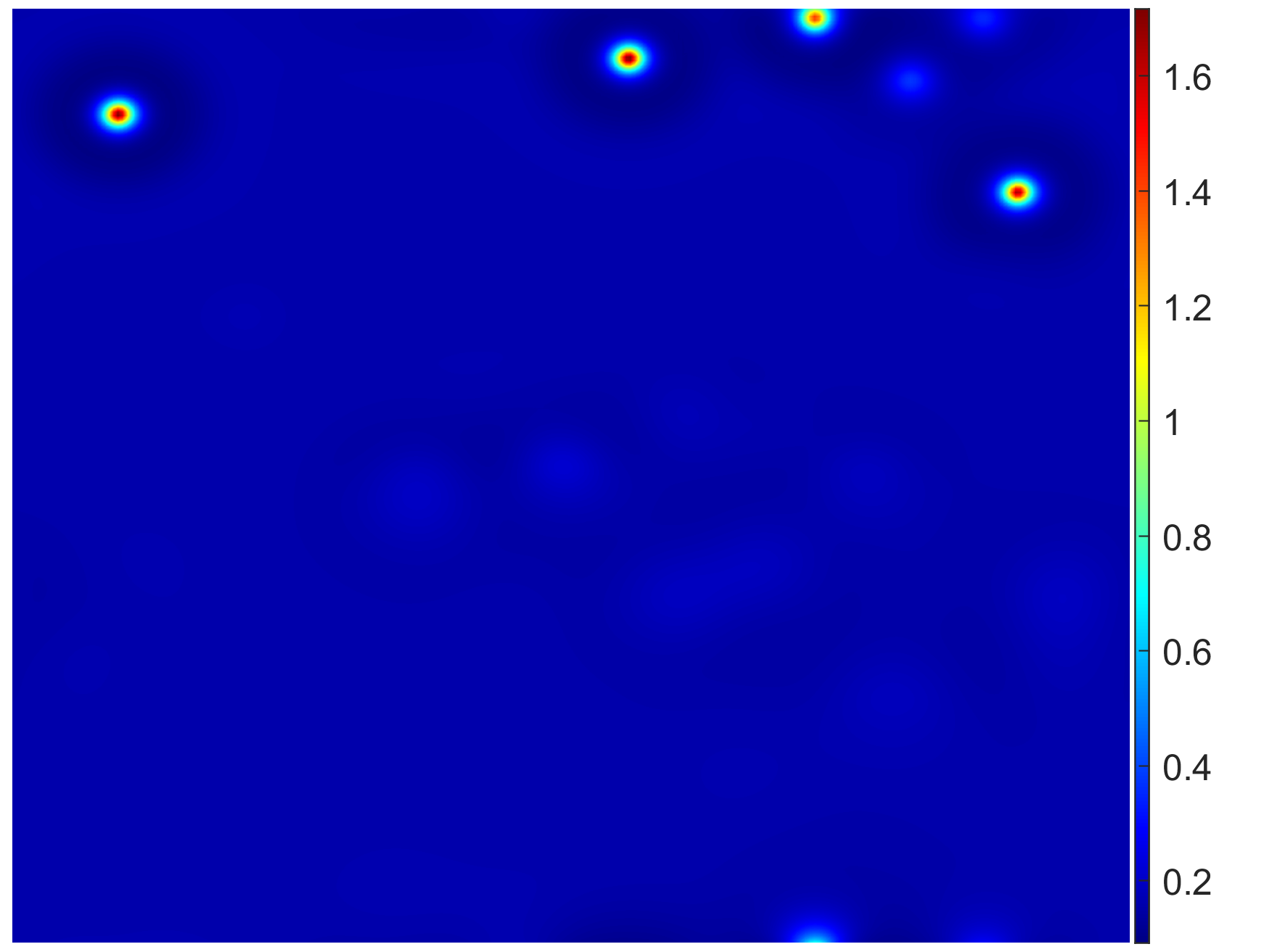} 
				\includegraphics[width=0.25\textwidth]{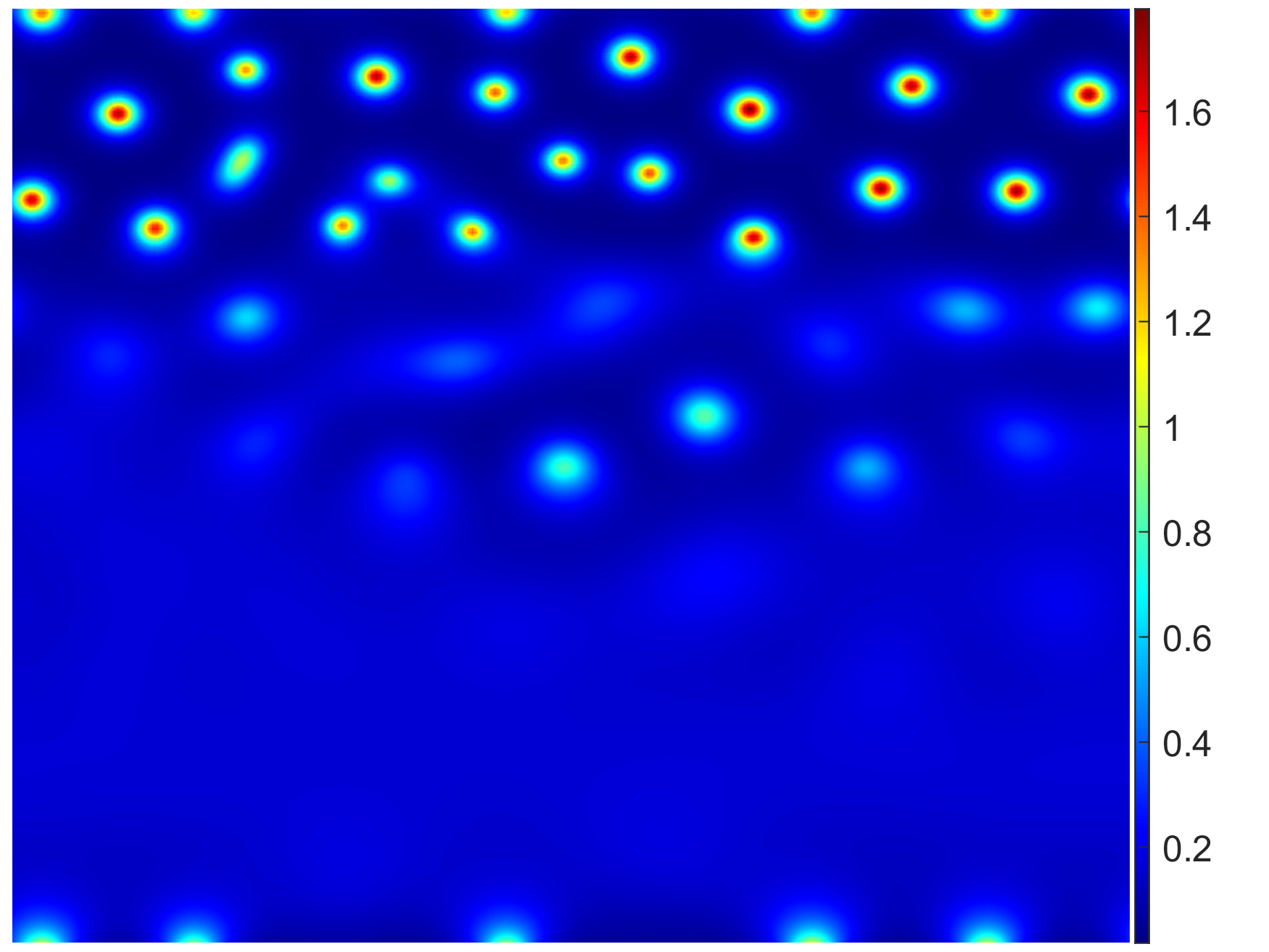} 
				\includegraphics[width=0.25\textwidth]{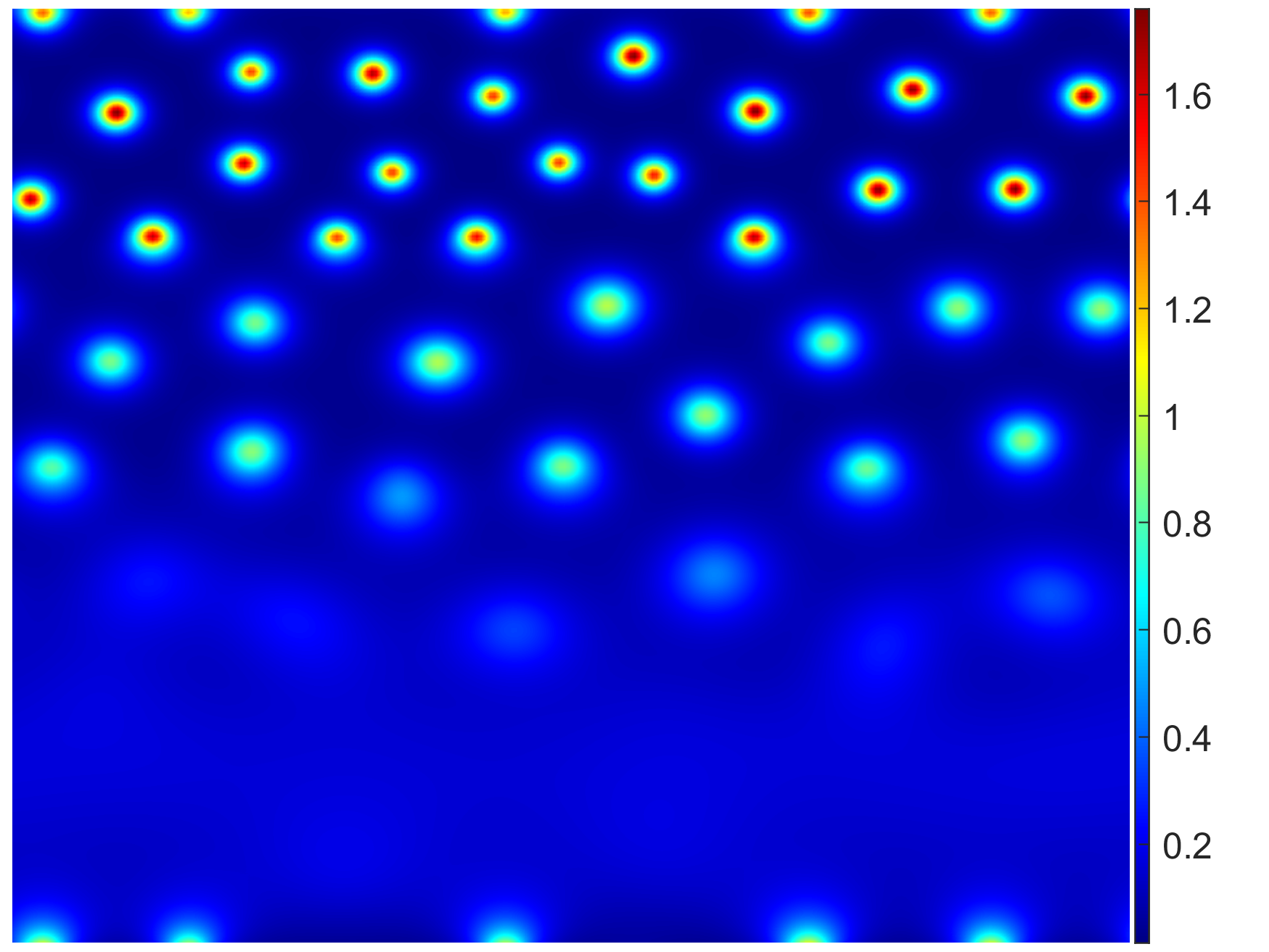}
				\includegraphics[width=0.25\textwidth]{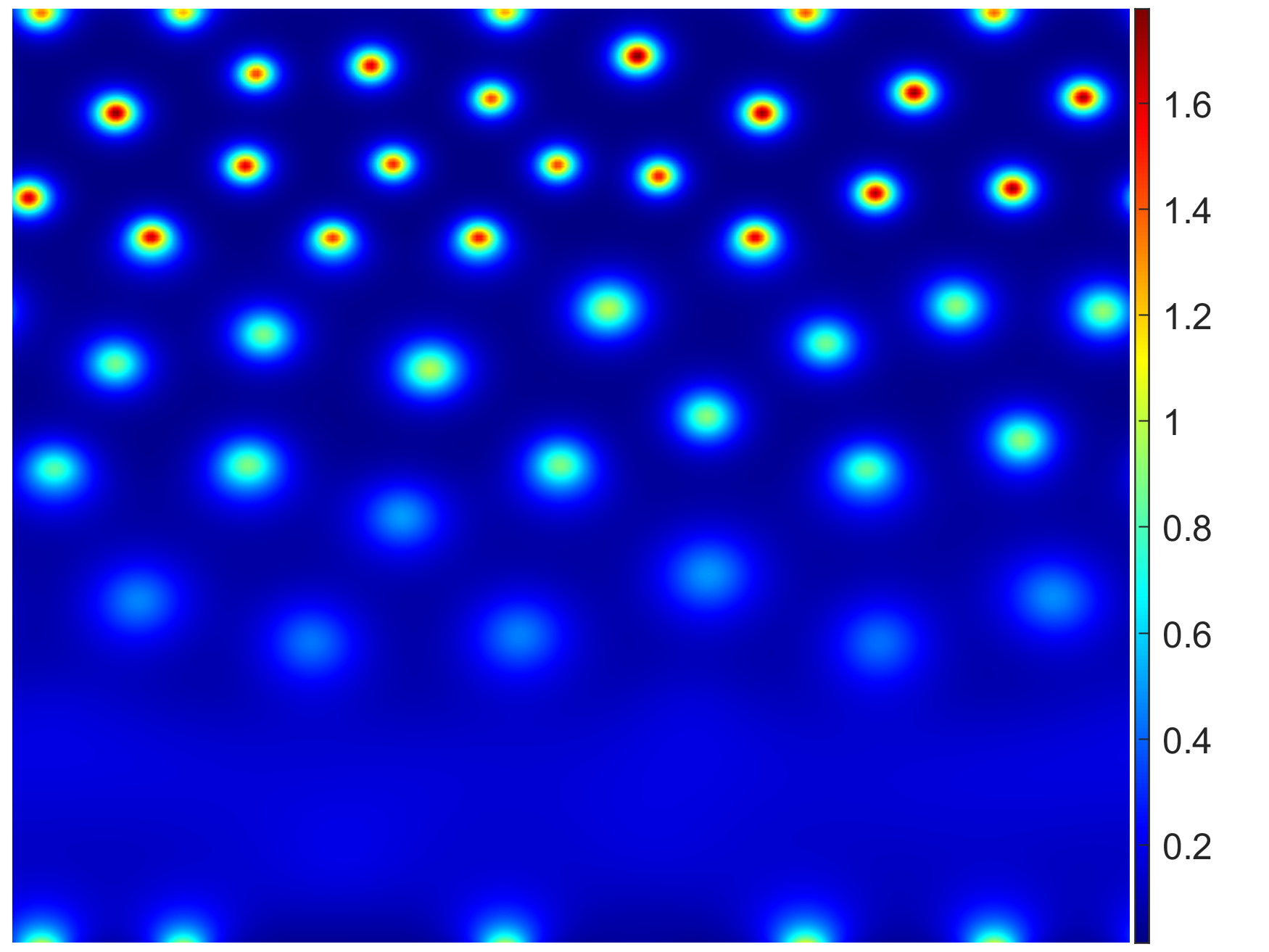} 
			}
			\\\vspace{-5pt}
			\subfigure
			{
				\includegraphics[width=0.25\textwidth]{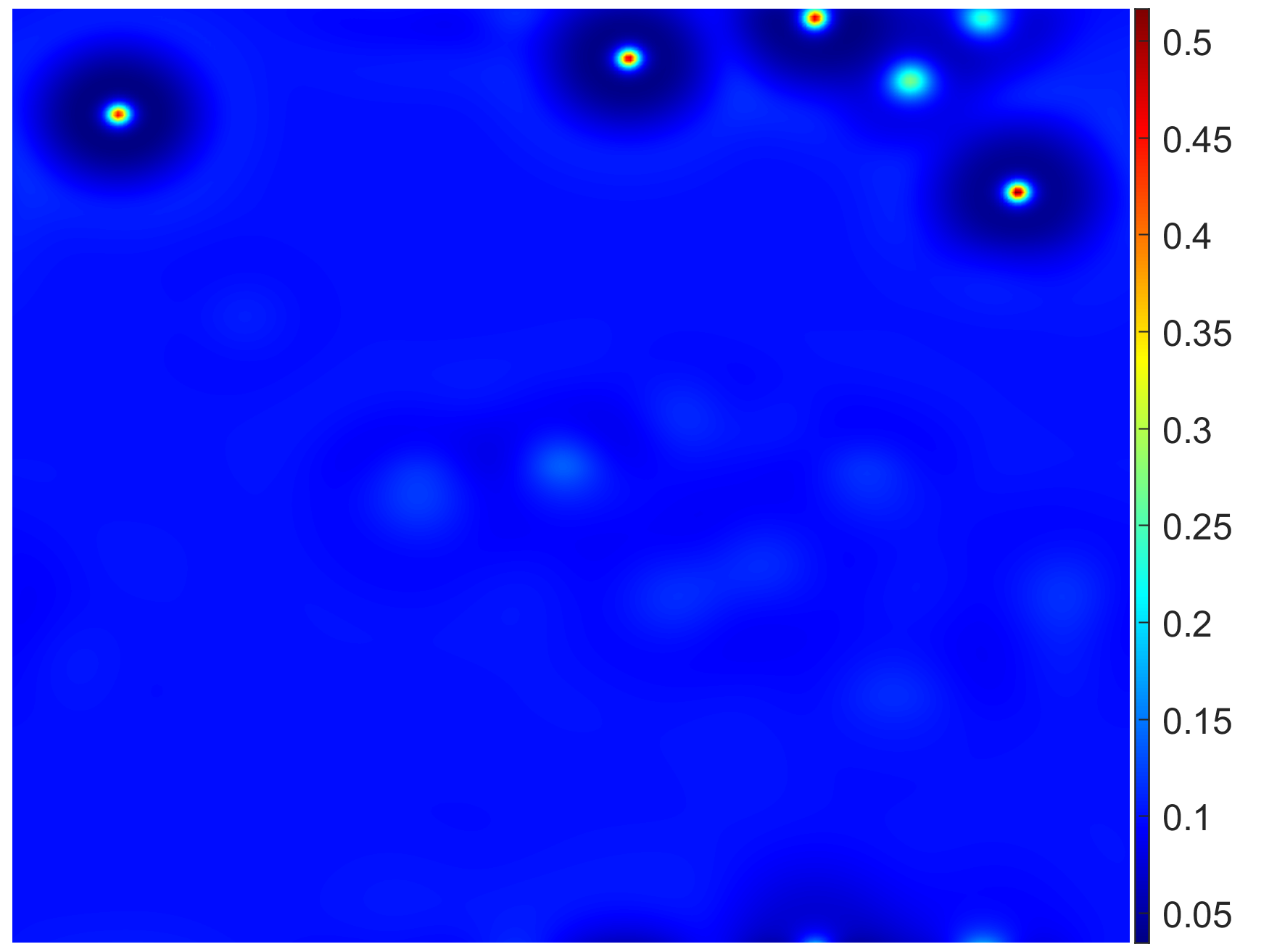}
				\includegraphics[width=0.25\textwidth]{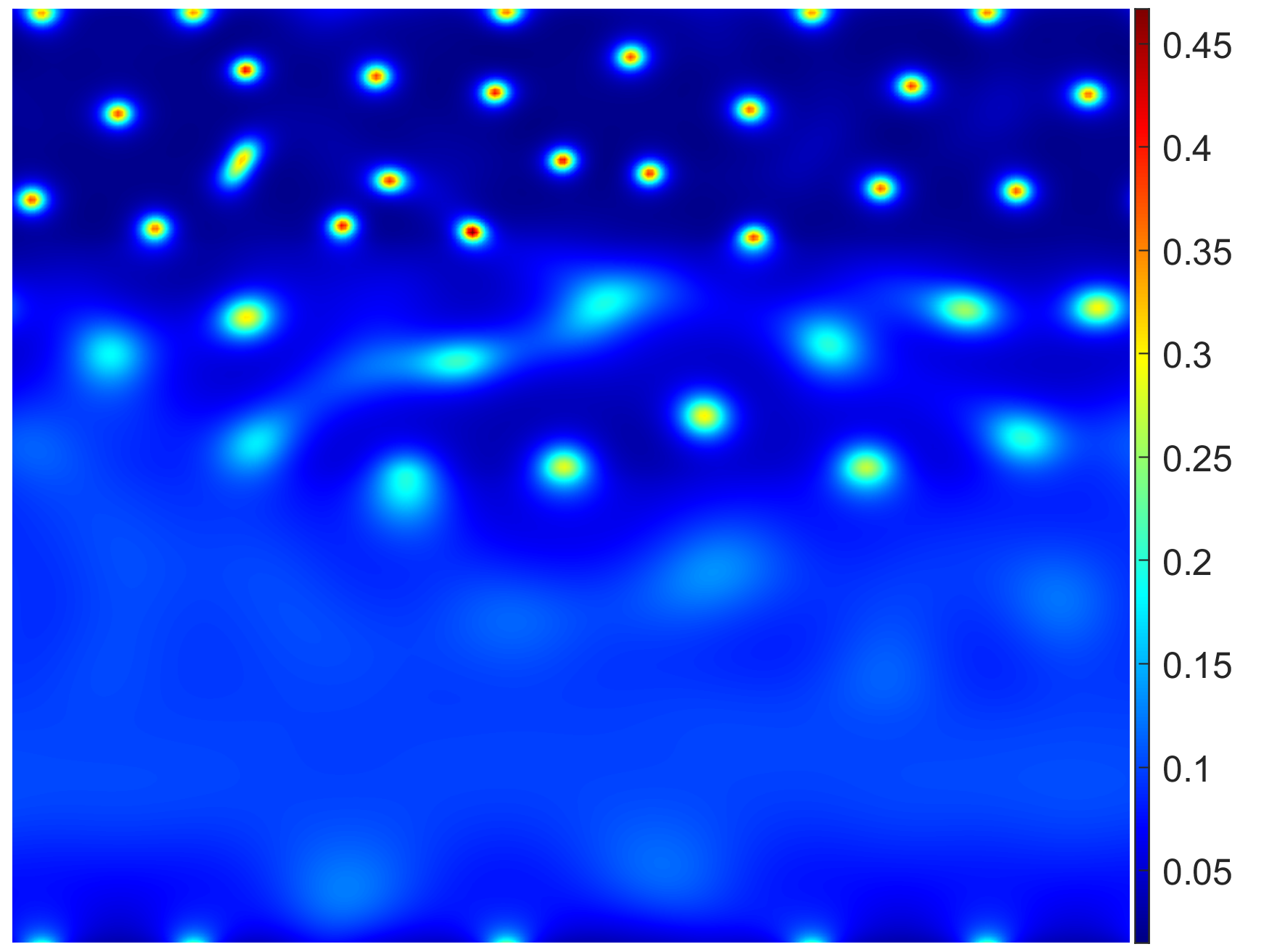} 
				\includegraphics[width=0.25\textwidth]{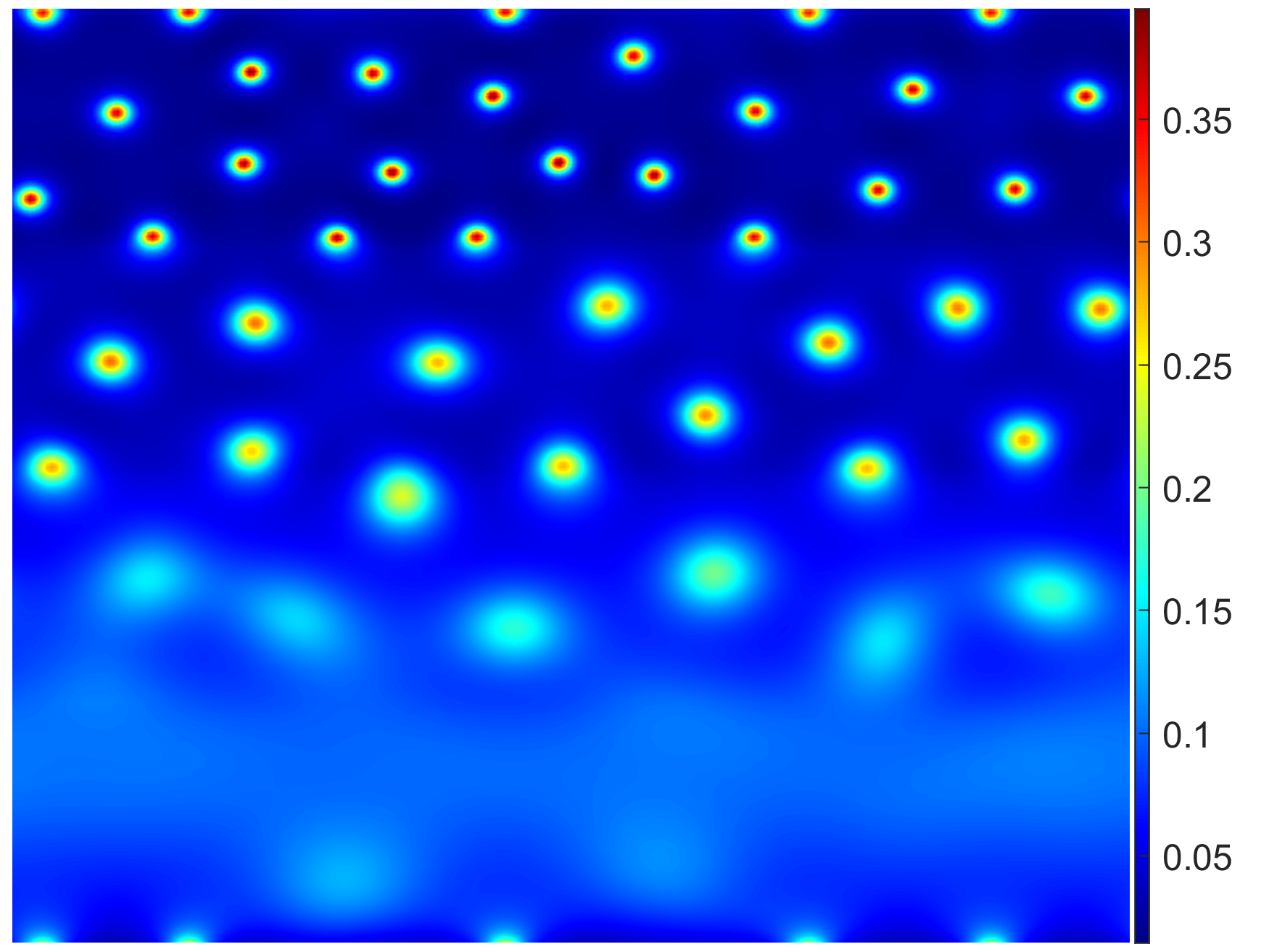}
				\includegraphics[width=0.25\textwidth]{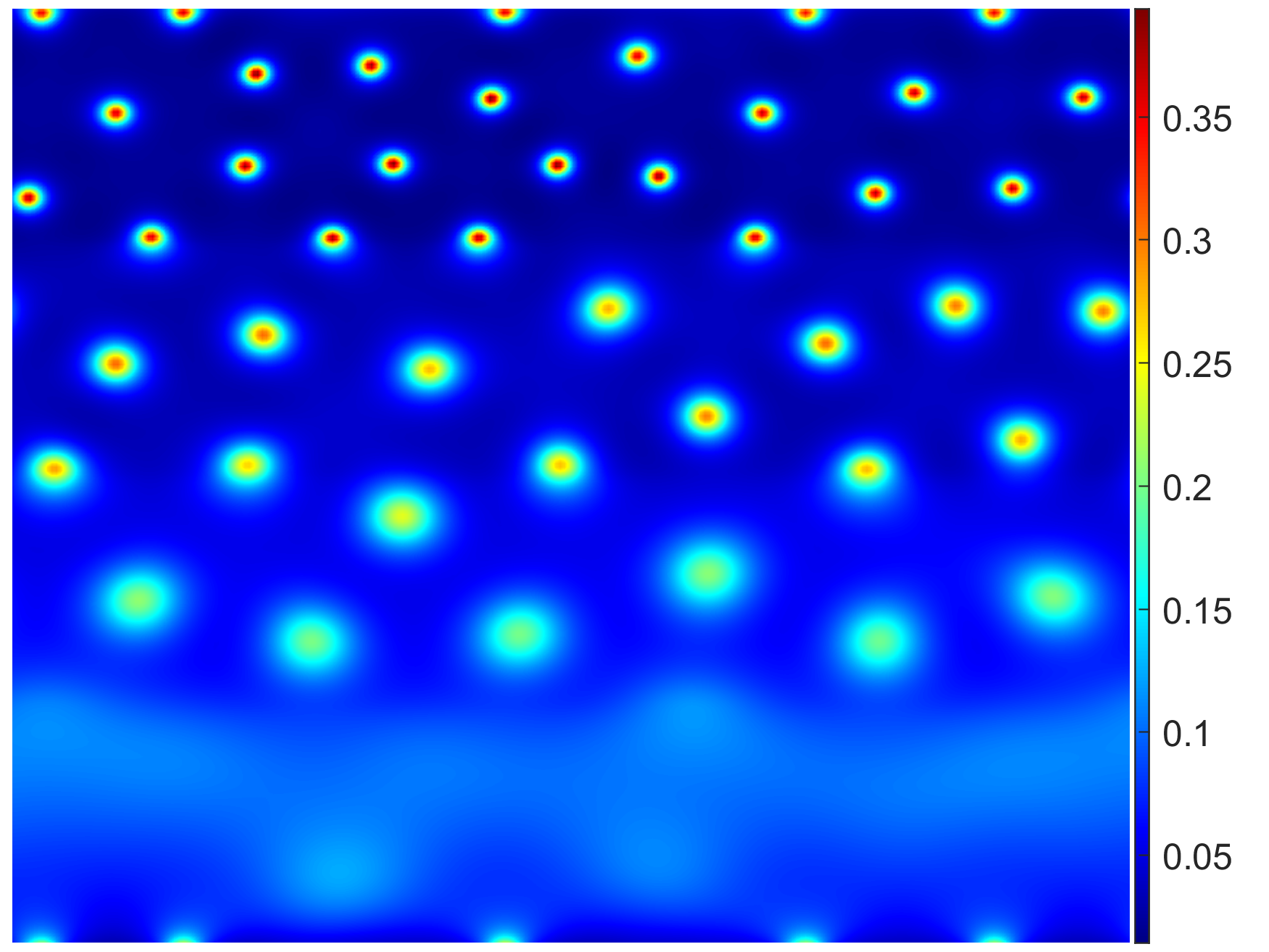}
			}
			\vspace{-15pt}
			\caption{Case with $\eta$ defined piecewise \eqref{eta_change_equation}: Evolution of the level of attractiveness $A$ (top) and
				density $\rho$ (bottom) at $T = 50,200,400,800$.
			}
			\label{fig_eta_change}
		\end{figure}

		Figure~\ref{fig_eta_change} shows the evolution of the
		attractiveness field $A$ (top row) and the density $\rho$
		(bottom row).
		Starting from nearly homogeneous initial data,
		localized hotspot structures gradually emerge and
		self--organize across the domain.	
		A clear spatial variation in the characteristic hotspot scale
		can be observed.
	
		In the regions where $\eta$ is smaller, the diffusion of the
		attractiveness field becomes weaker, which strengthens the
		local aggregation mechanism and leads to small scattered hotspots.
		In contrast, larger values of $\eta$ produce stronger diffusion,
		resulting in fewer but larger hotspots.
		Consequently, the spatial pattern exhibits a gradual transition
		of hotspot scales along the $x$-direction, with hotspot size
		increasing as $\eta$ becomes larger.

			These observations are consistent with the dependence of hotspot
		size on $\eta$ observed in the early experiments with homogeneous coefficients.


\begin{example}
To investigate the influence of spatial heterogeneity in the baseline
attractiveness, we consider the case where $A_0$ is defined piecewisely in $x$-direction:
	\begin{equation}\label{A0_change_equation}
	A_0(x, y)=
	\left\{
	\begin{aligned}
		&\frac{1}{7.5}, &&0 \le x < \pi,\\
		&\frac{1}{30}, &&\pi \le x \leq 2\pi,
	\end{aligned}
	\right.
\end{equation}
with the same initial condition in Eq.~\eqref{init_condition}.
In this setting, the domain is divided into two regions with different
background attractiveness levels.
The upper half ($0\le x<\pi$) has a relatively larger baseline
attractiveness $A_0=1/7.5$, while the lower half
($\pi\le x\le 2\pi$) corresponds to a smaller value $A_0=1/30$.
\end{example}

\begin{figure}[h!]
	\centering
	\subfigure
	{
		\includegraphics[width=0.25\textwidth]{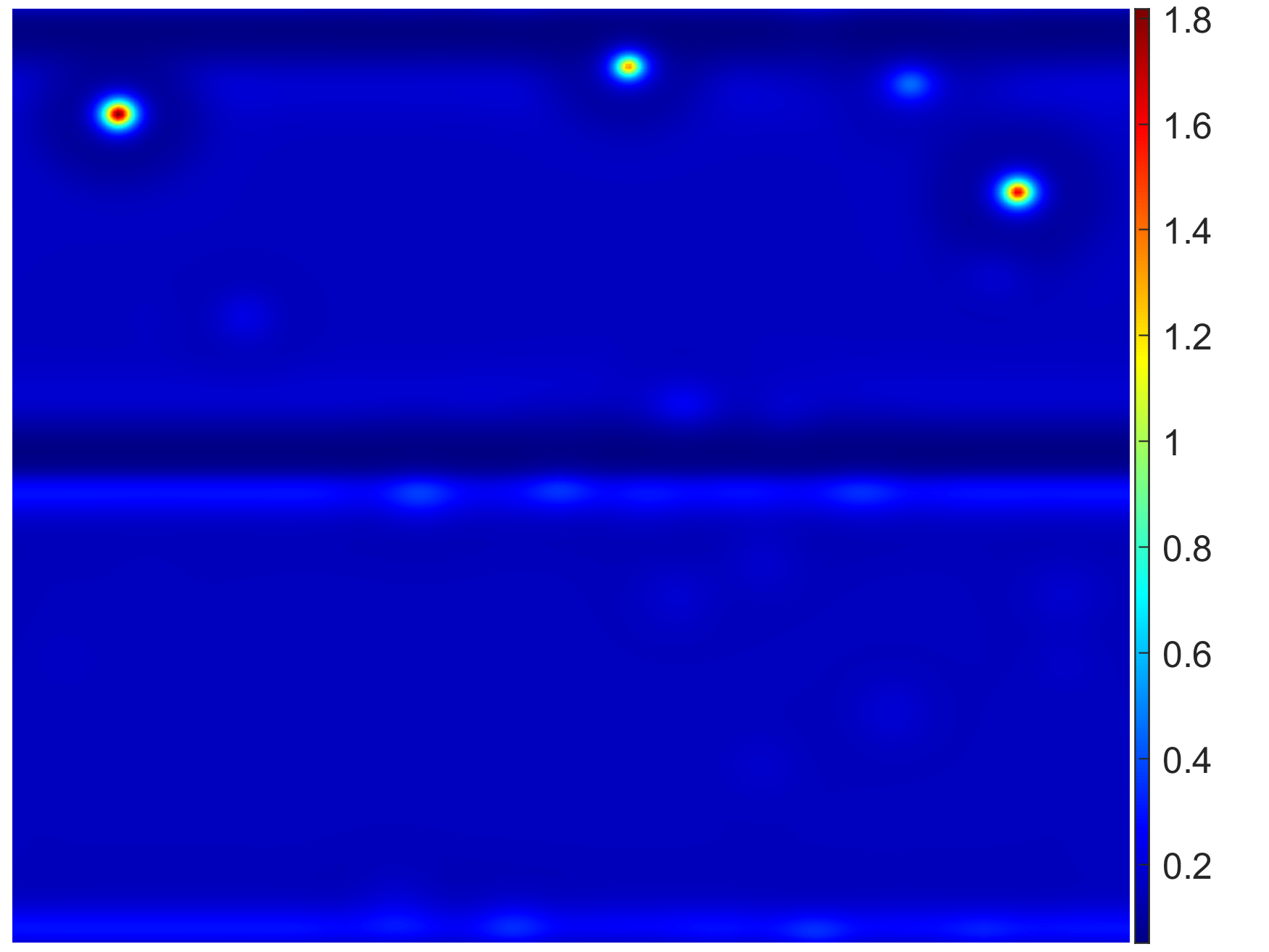} 
		\includegraphics[width=0.25\textwidth]{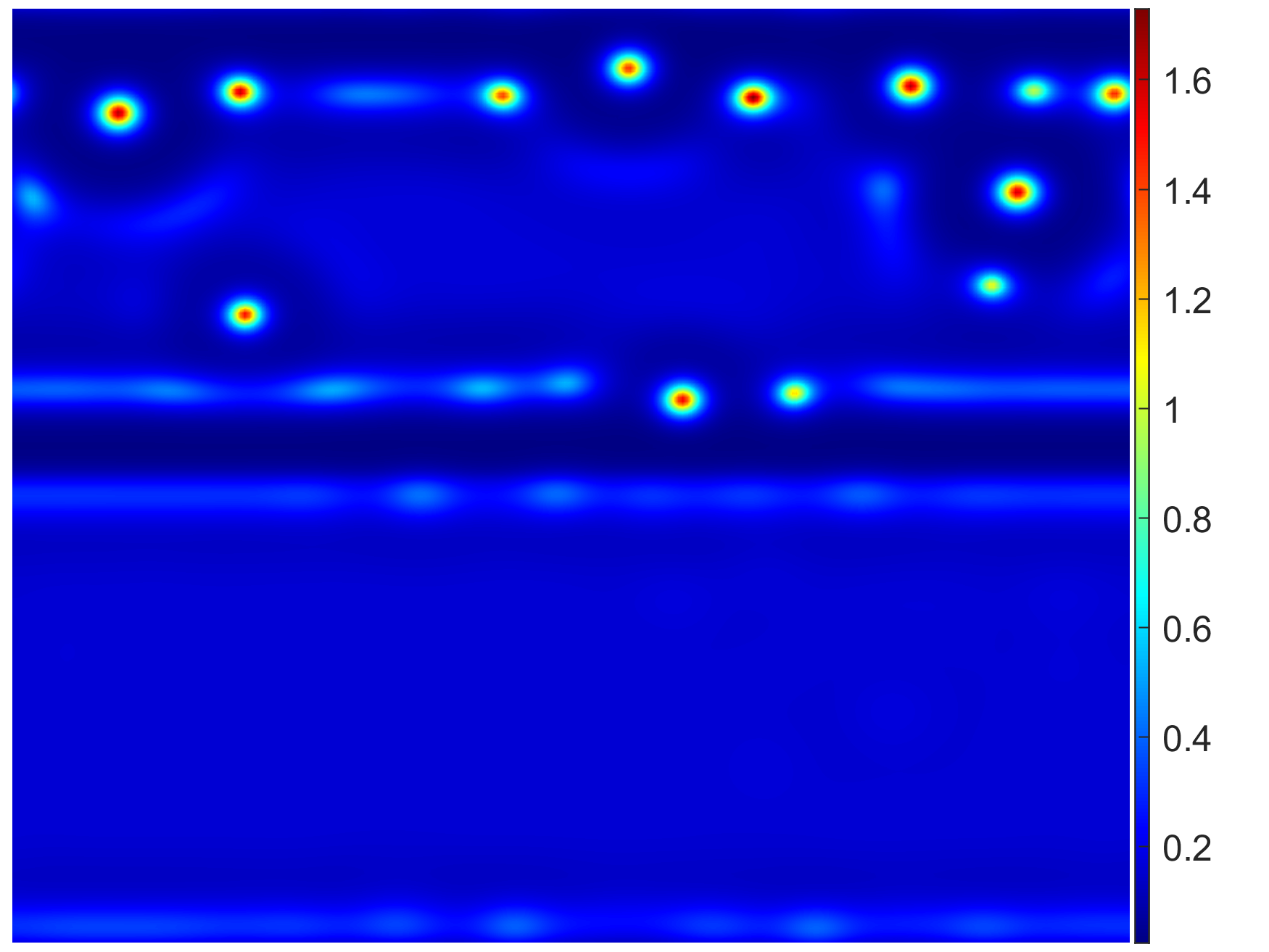} 
		\includegraphics[width=0.25\textwidth]{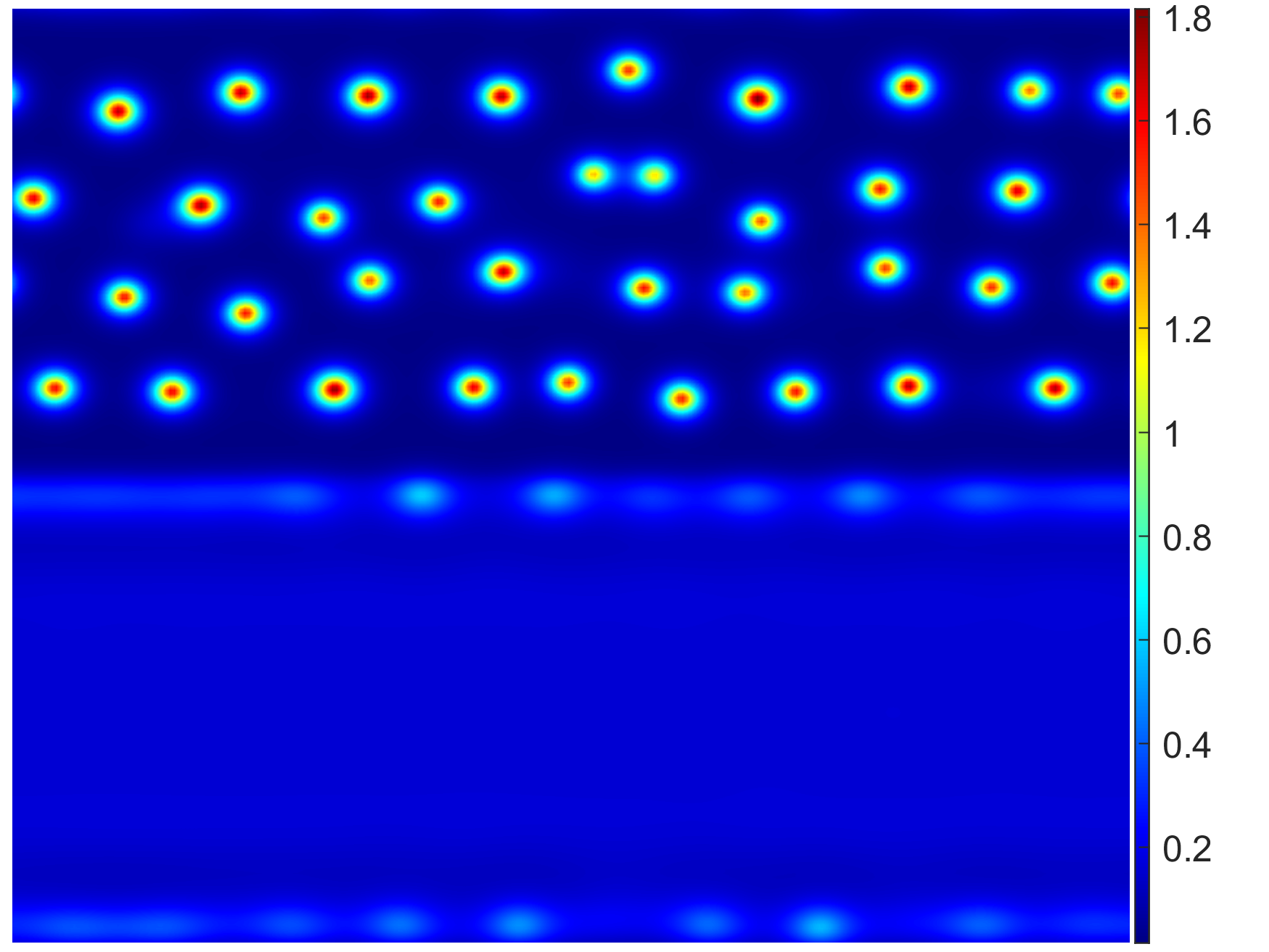}
		\includegraphics[width=0.25\textwidth]{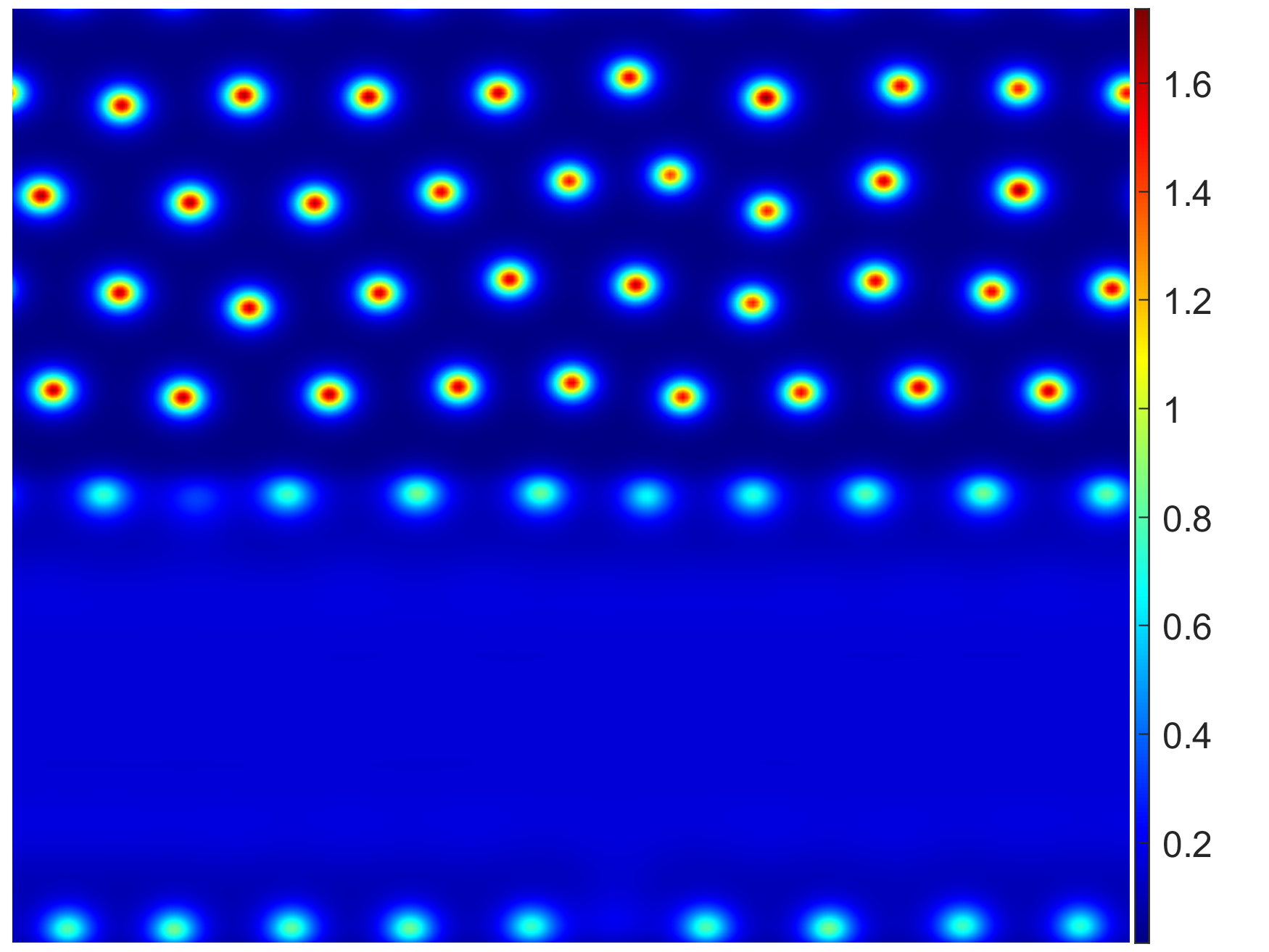} 
	}
	\\
	\subfigure
	{
		\includegraphics[width=0.25\textwidth]{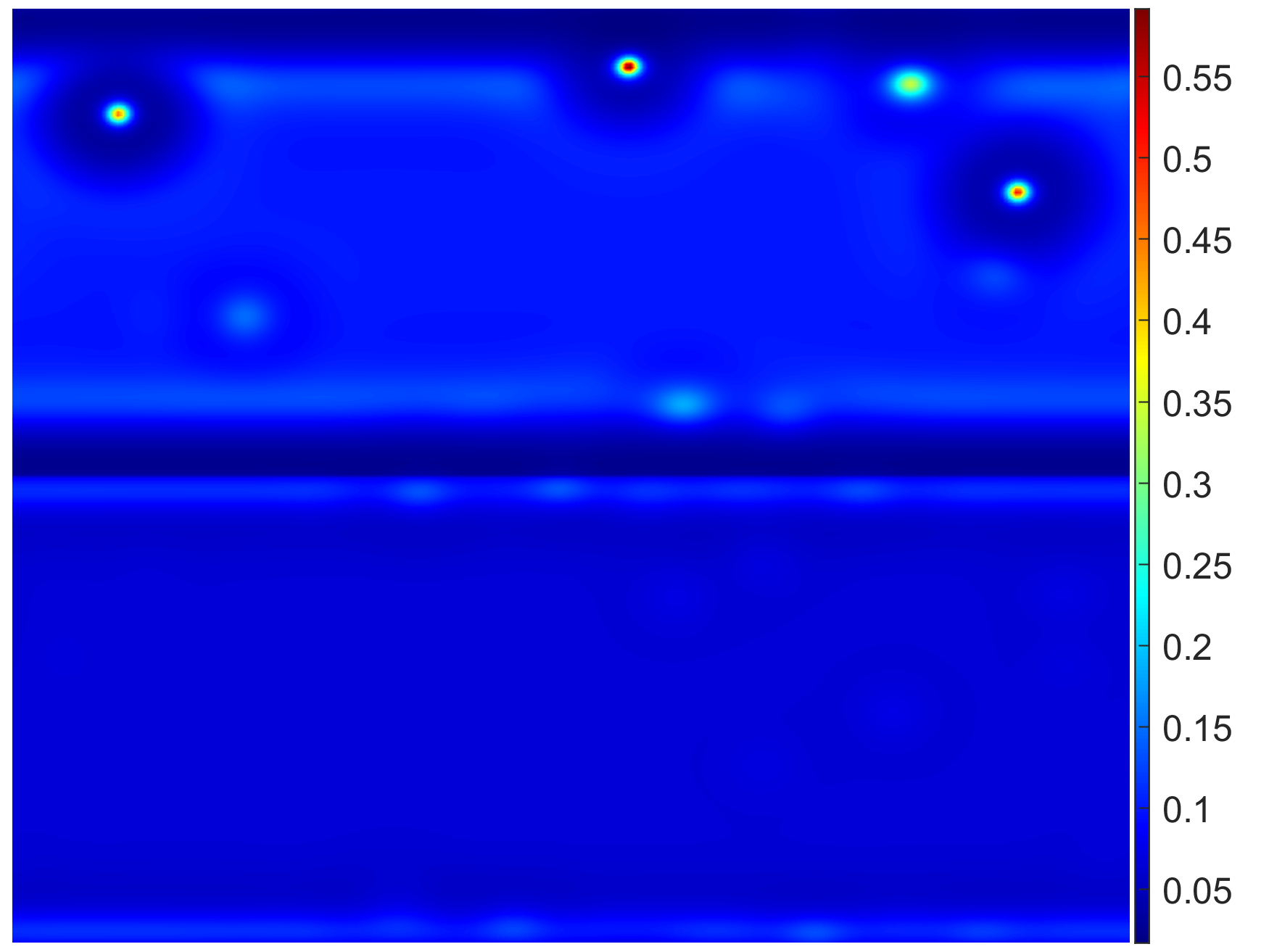}
		\includegraphics[width=0.25\textwidth]{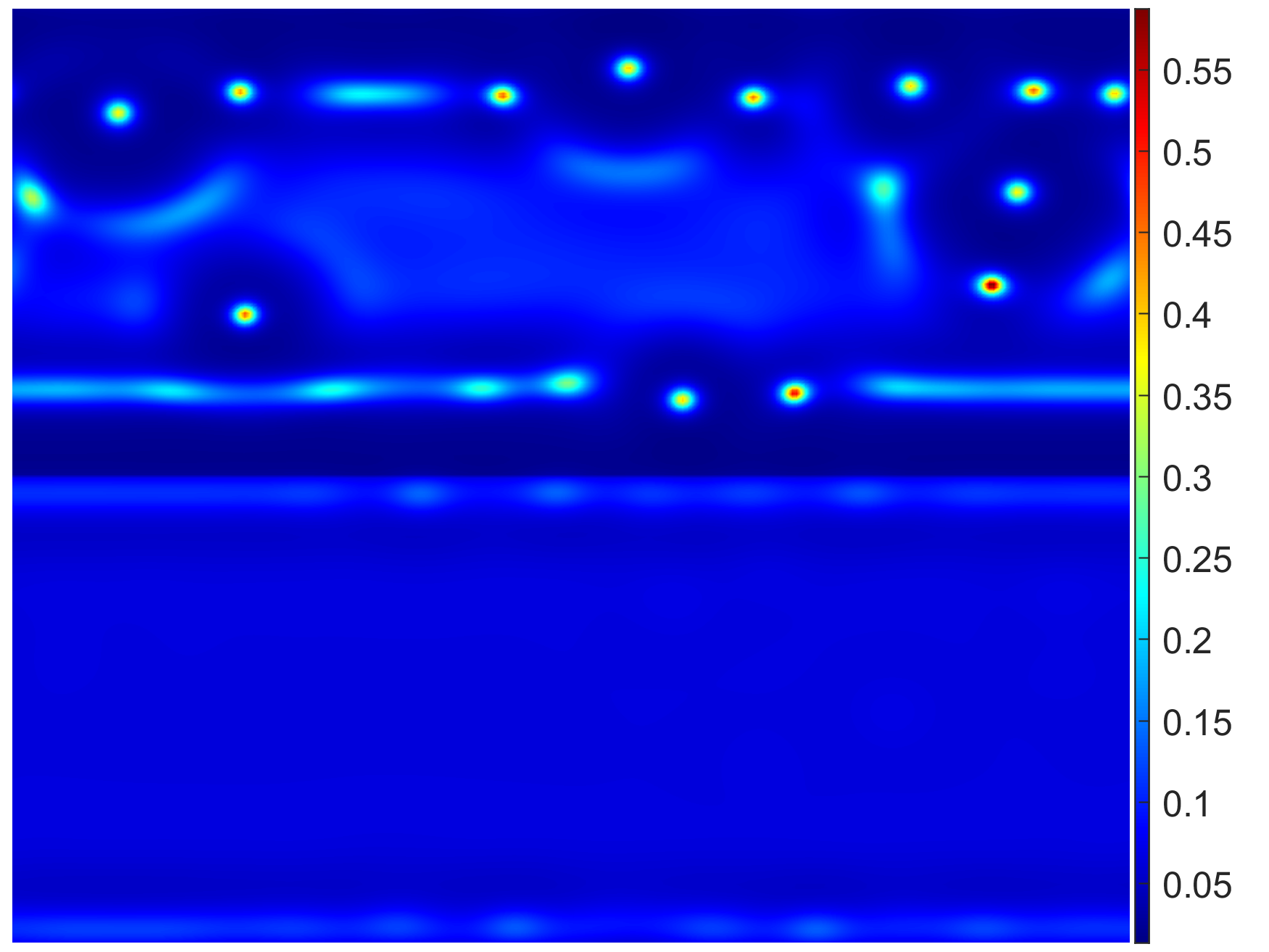} 
		\includegraphics[width=0.25\textwidth]{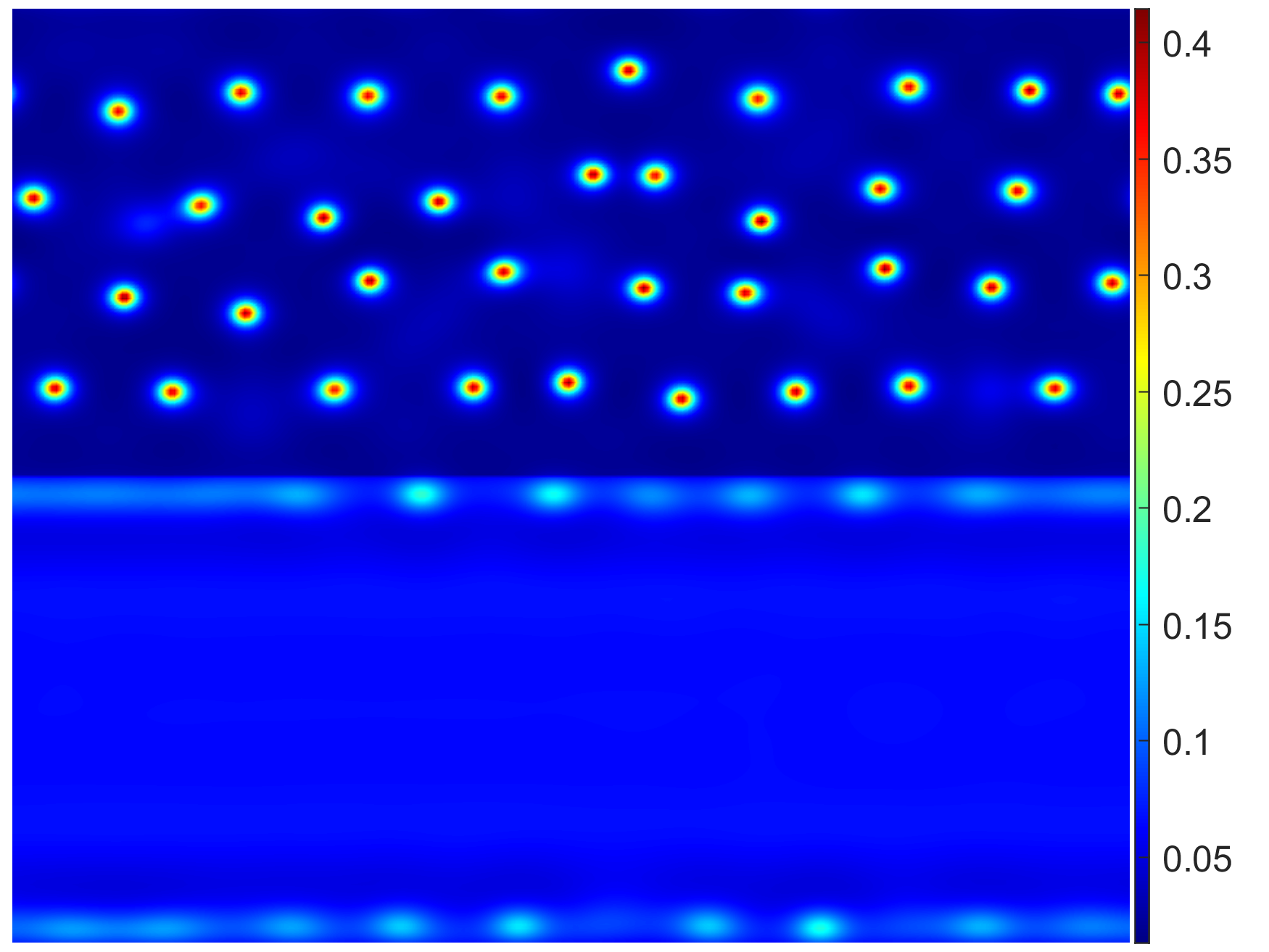}
		\includegraphics[width=0.25\textwidth]{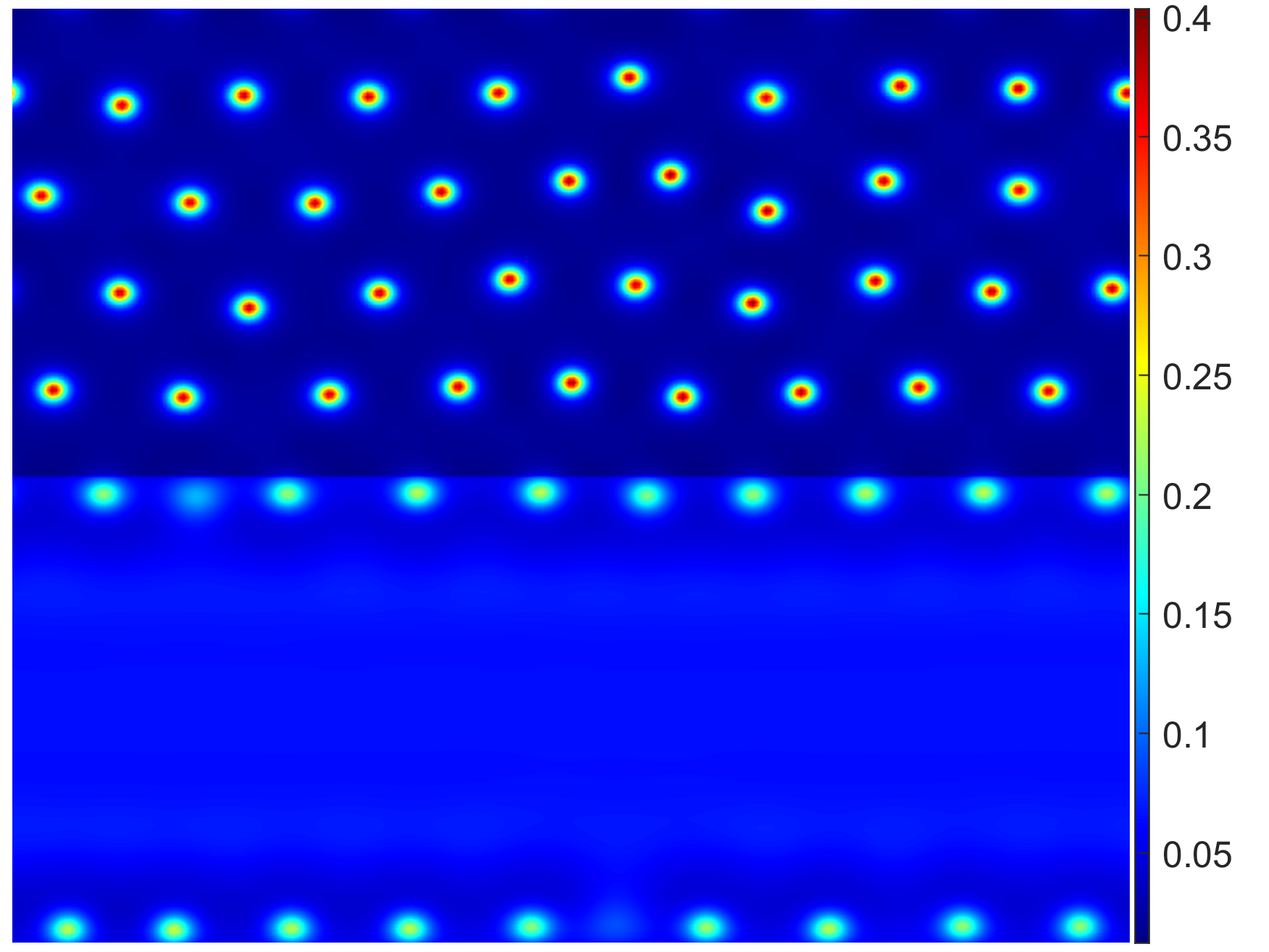}
	}
	\vspace{-15pt}
	\caption{Case with $A_0$ defined piecewise \eqref{A0_change_equation}: Evolution of the level of attractiveness $A$ (top) and
		density $\rho$ (bottom) at $T = 50,100,200,800$.
	}
	\label{fig_A0_change}
\end{figure}

Figure~\ref{fig_A0_change} illustrates the temporal evolution of the
attractiveness field $A$ and the criminal density $\rho$.
Starting from random initial perturbations, localized peaks gradually
emerge and organize into hotspot patterns.
Due to the spatial variation of $A_0$, the distribution of hotspots is
not homogeneous across the domain.
In particular, the region with smaller baseline attractiveness tends to
develop stronger and more pronounced hotspot structures, indicating that
spatial heterogeneity in the environmental attractiveness plays a
significant role in shaping the long--time spatial pattern of crime.

\subsection{Pattern formation and phase transition in epidemic modeling}
Similarly, we first conduct a convergence study under a smooth initial condition and representative parameters.

\begin{example}
The initial data is chosen as
\begin{equation}\label{epidemic_initial_data}
\begin{aligned}
	\psi_{i}(x,y,0) &= 0.1\sin(2\pi (x-0.5))\cos(2\pi (y-0.5))+0.2,\quad i = 1,\dots,8,\\
	p(x,y,0) &= 0.5\sin(2\pi (x-0.5))\cos(2\pi (y-0.5))+1, 
\end{aligned}
\end{equation}
where periodic boundary condition is used and the parameters in the model are: $\Omega = [0, 1] \times [0,1]$, $p^0 = 1/30$, $\eta = 1$, $D = 0.01$. For convenience, all other parameters are set to 0.5.
\end{example}

The simulation parameters are set as the same as in Example~\ref{ex:crime-convergence}.
Table~\ref{epidemic_space_L2_phi} reports the spatial errors and convergence orders, while Table~\ref{epidemic_time_L2_phi} presents the temporal errors and convergence in time, which again verifies our theoretical results.

\begin{table}[h!]
\centering
\caption{\small Spatial error analysis with initial value \eqref{epidemic_initial_data} and $h_{ref}=1/512$.}
\label{epidemic_space_L2_phi}
\begin{tabular}{c|c|c|c|c|c}
\hline\hline
$h$  & $\Vert e_{\bm{\Psi}}\Vert_{L^2}$ &  Order & $\Vert e_p\Vert_{H^1}$ & Order & Convergence \\ \hline
$1/8$   &  $1.7723\times10^{-4}$  & - &  $7.2969\times10^{-3}$ & - & \multirow{4}{*}{\centering\includegraphics[width=0.2\textwidth]{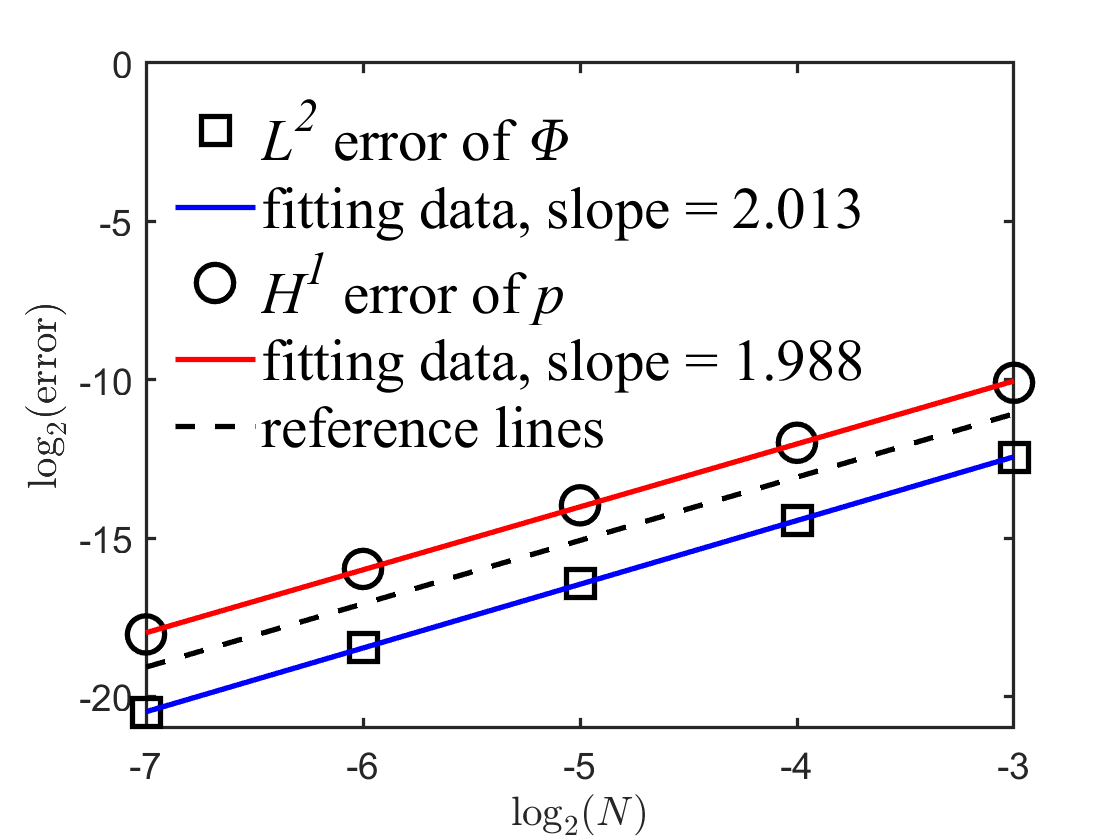}} \\
$1/16$ & $4.5071\times10^{-5}$ &1.9753  &$1.9624\times10^{-3}$ & 1.8947 & \\ 
$1/32$ & $1.1288\times10^{-5}$ & 1.9973 &$4.9928\times10^{-4}$ & 1.9747 & \\ 
$1/64$ & $2.7923\times10^{-6}$ & 2.0153 & $1.2401\times10^{-4}$ & 2.0093 & \\ 
$1/128$ & $6.6503\times10^{-7}$ & 2.0700 & $2.9568\times10^{-5}$ & 2.0684 & \\ 
\hline\hline
\end{tabular}
\end{table}

\begin{table}[h!]
\centering
\caption{\small Temporal error analysis with initial value \eqref{epidemic_initial_data} and $\tau_{ref}=10^{-5}$.}
\label{epidemic_time_L2_phi} 
\begin{tabular}{c|c|c|c|c|c}
\hline\hline
$\tau$  & $\Vert e_{\bm{\Psi}}\Vert_{L^2}$ &  Order & $\Vert e_p\Vert_{H^1}$& Order & Convergence \\ \hline
$2\times10^{-3}$   &    $1.6562\times10^{-8}$  & - &  $1.2829\times10^{-7}$ & - & \multirow{4}{*}{\centering\includegraphics[width=0.2\textwidth]{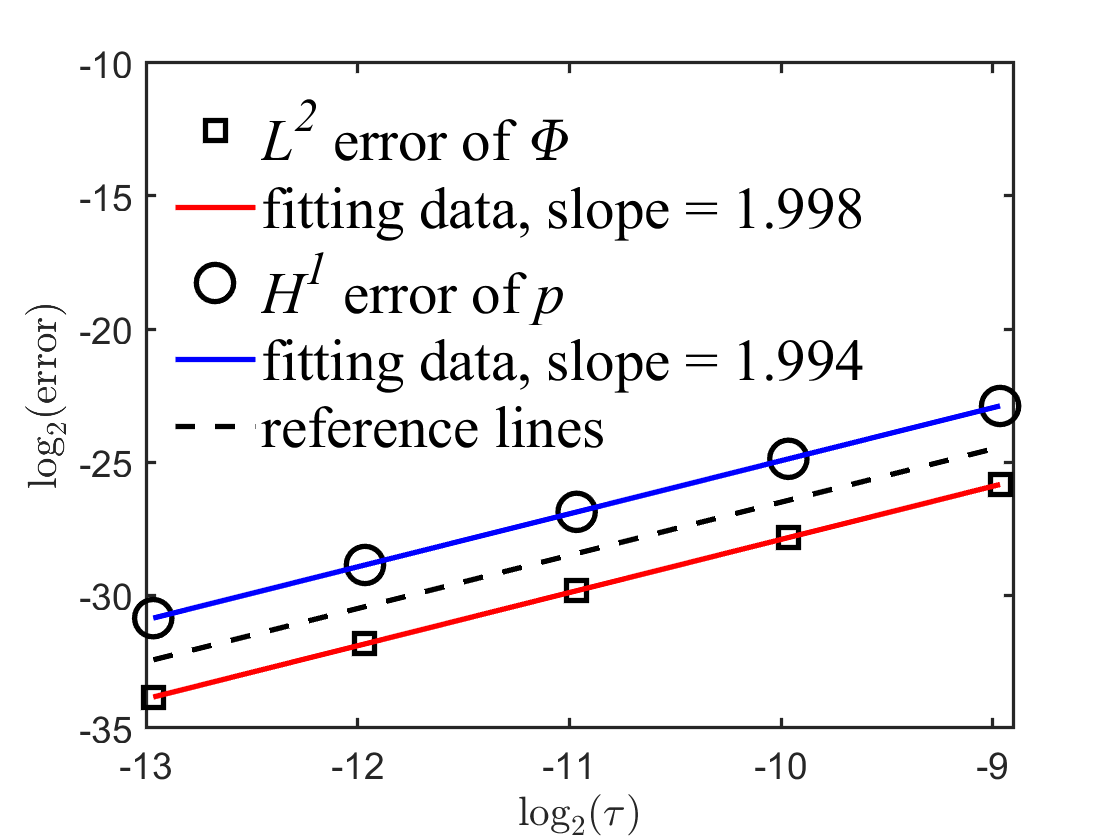}} \\
$1\times10^{-3}$   & $4.1608\times10^{-9}$ &1.9929 & $3.2371\times10^{-8}$ &  1.9866& \\ 
$5\times10^{-4}$   & $1.0426\times10^{-9}$ & 1.9967 & $8.1297\times10^{-9}$ & 1.9934 & \\ 
$2.5\times10^{-4}$  & $2.6083\times10^{-10}$ & 1.9990 & $2.0367\times10^{-9}$ & 1.9969 & \\ 
$1.25\times10^{-4}$  & $6.5077\times10^{-11}$ & 2.0029 & $5.0913\times10^{-10}$ & 2.0002 & \\ 
\hline\hline
\end{tabular}
\end{table}

Now we begin to investigate phase transitions between spatial aggregation and dissipation in epidemic dynamics.

\begin{example}\label{epid_simulation_plot}
Let
$S(\mathbf{x},0),E(\mathbf{x},0),P(\mathbf{x},0),A(\mathbf{x},0),I^-(\mathbf{x},0),I^+(\mathbf{x},0)$
be the initial values  at each point \(\mathbf{x} \in [0,1] \times [0,1]\), respectively. We set
\begin{equation*}
	\begin{split}
		S(\mathbf{x},0) & = 1- \delta_1(\mathbf{x})- \delta_2(\mathbf{x}) -\delta_3(\mathbf{x}) - \delta_4(\mathbf{x}) - \delta_5(\mathbf{x}), \\
		E(\mathbf{x}, 0) & =   \delta_1(\mathbf{x}), \quad P(\mathbf{x}, 0) =  \delta_2(\mathbf{x}), \quad A(\mathbf{x}, 0) =  \delta_3(\mathbf{x}), \quad I^{-}(\mathbf{x},0) = \delta_4(\mathbf{x}),\quad I^{+} = \delta_5(\mathbf{x}),
	\end{split}
\end{equation*}
where each perturbation $\delta_i(\mathbf{x})$ is composed of 30 independent Gaussian functions with randomly chosen centres $(x_1^i,x_2^i)$, heights $h_i$, and widths $\sigma_i$, as defined in Eq.~\eqref{delta_set}. Specifically, $h_i = 0.001 r_i^{(1)}$, $\sigma_i = 0.01 r_i^{(2)}$, and $(x_i, j_i) = (r_i^{(3)}, r_i^{(4)})$, where $r_i^{(1)}, r_i^{(2)}, r_i^{(3)}, r_i^{(4)}$ are independent samples drawn from the uniform distribution on $[0,1)$. Periodic images are added to ensure the periodic boundary condition of all agent and field variables ($L=20$). 			

\end{example}		

\begin{figure}[!h]
	\centering
	\includegraphics[width=1\textwidth]{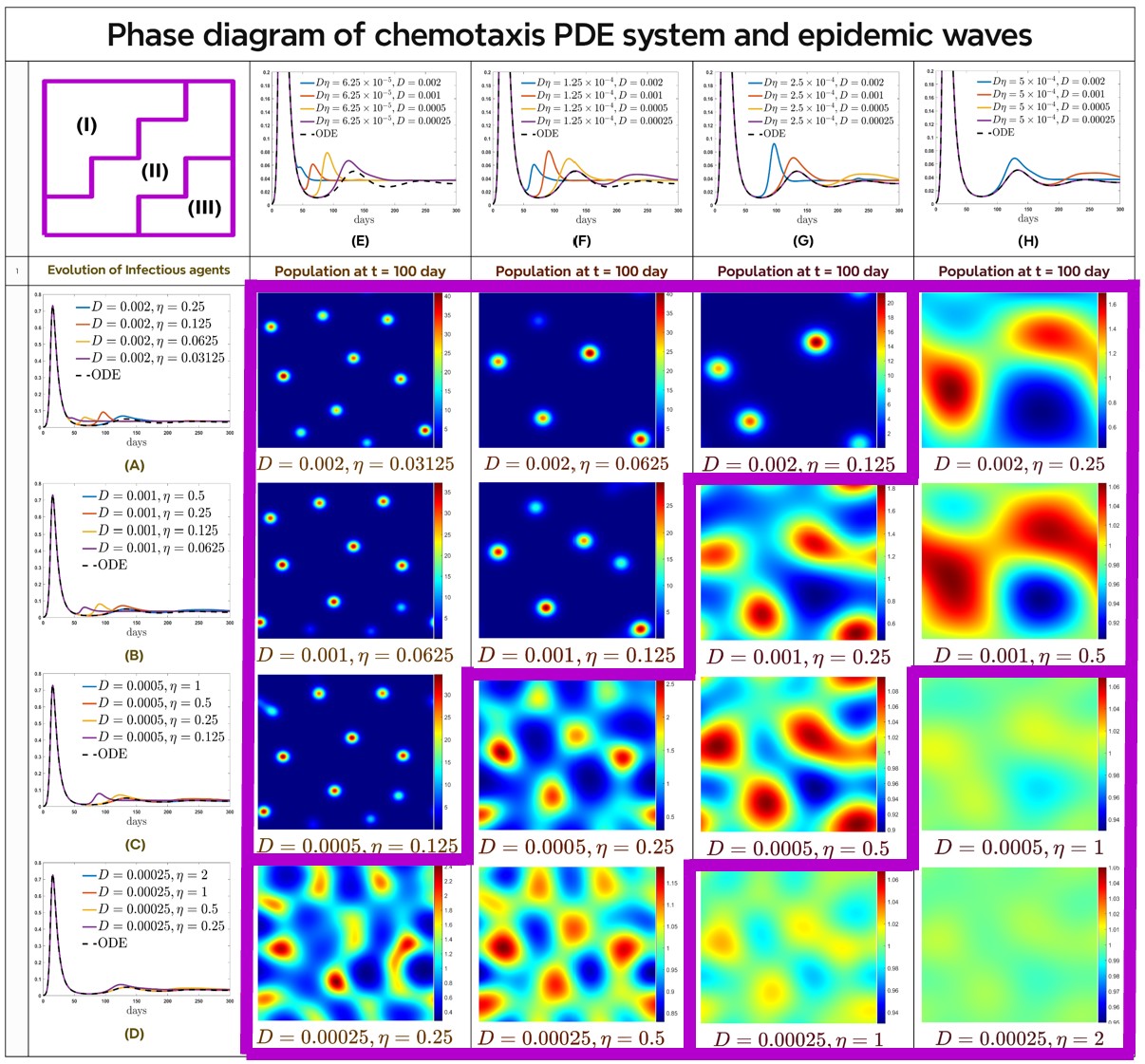}\vspace{-15pt}
	\caption{2-D spatial distribution of mobile agents ($S+E+P+A+I^-+I^++R$) and the time evolution of virus carriers (i.e., the renormalized total population of $E$, $P$, $A$, $I^+$ and $I^-$) under different combinations of $D$ and $\eta$, with the corresponding ODE results used as references. As $\eta$ or $D$ increases,  the phase diagram clearly reveals the transition among three phases: the aggregating phase (I), intermediate phase (II), and dissipative phase (III). Panels (A)–(D) show that aggregation has little impact on the first epidemic peak. In panels (E)–(H), high-density aggregation significantly aggravates the second epidemic peak.}
	\label{epidemic_phase_plot}
\end{figure}

Denote the total density of mobile agents participating in the epidemic
dynamics by
\[
u(\mathbf{x},t):=S(\mathbf{x},t)+E(\mathbf{x},t)+P(\mathbf{x},t)+A(\mathbf{x},t)+I^{-}(\mathbf{x},t)+I^{+}(\mathbf{x},t)+R(\mathbf{x},t).
\]
Figure~\ref{epidemic_phase_plot} presents the
phase diagram of $u(\mathbf{x},t)$ at $t=100$ under different combinations of $\eta$ and $D$, starting from the same
initial condition. Moreover,  we present the time evolution of active virus carriers, i.e.,  $\int_{\Omega} (E+P+A +I^+ + I^-) \mathrm{d} \mathbf{x}$, in Panels (A)-(D) under the same $D$ but different $\eta$ to show the first infectious peak. In Panels (E)-(H), the second epidemic waves of active virus carriers are compared under the fixed product $\eta D$ but varying $D$. In each panel, the results of the corresponding ODE model (i.e., without spatial heterogeneity) are used as the reference.

The phase diagram clearly reveals the transition among three phases: the aggregating phase (I), intermediate phase (II), and dissipative phase (III). 
Under fixed $D$, increasing $\eta$ will suppress the aggregation of agents and impose dissipation on the spatial distribution. For fixed $\eta D$, when agents move faster, it drives the spatial distribution from the aggregating phase to the dissipative phase. Moreover, when the system evolves for a sufficiently long time, the boundary becomes more evident, thereby providing a numerical evidence of phase transitions between aggregation-dominated phase (I) and dissipation-dominated one (III).

If fact, the hotspots in epidemic dynamics that reflect that aggregation of agents may possibly enhance their contact rate, and consequently have a significant influence on the epidemic waves. Panels (A)–(D) show that aggregation has little impact on the first epidemic peak. For a short time, the behavior of PDE system is dominated by the pattern of ODE as the spatial perturbation is not so large. However,  from panels (E)–(H), it is demonstrated that at the level of disease transmission, high-density aggregation significantly aggravates the second epidemic peak because the contact of clustered population density accelerates virus spread.


Motivated by the spatial patterns shown in the third row of
Figure~\ref{epidemic_phase_plot} (corresponding to $D=0.0005$),
we investigate how the spatial mixing parameter $\eta$
affects both the spatial heterogeneity and the epidemic dynamics.  To further quantify the phase transition observed in
Figure~\ref{epidemic_phase_plot}, we introduce an order parameter $Q_{std}(t)$ that
measures the spatial heterogeneity of the population density,
\[
Q_{std}(t)
=
\left(
\frac{1}{|\Omega|}
\int_{\Omega}
\left(u(\mathbf{x},t)-\bar u(t)\right)^2\, \D \mathbf{x}
\right)^{1/2},
\]
where $\bar u(t)$ denotes the spatial average of $u(\mathbf{x},t)$. 

\begin{figure}[h!]
	\centering
	\includegraphics[width=0.32\textwidth]{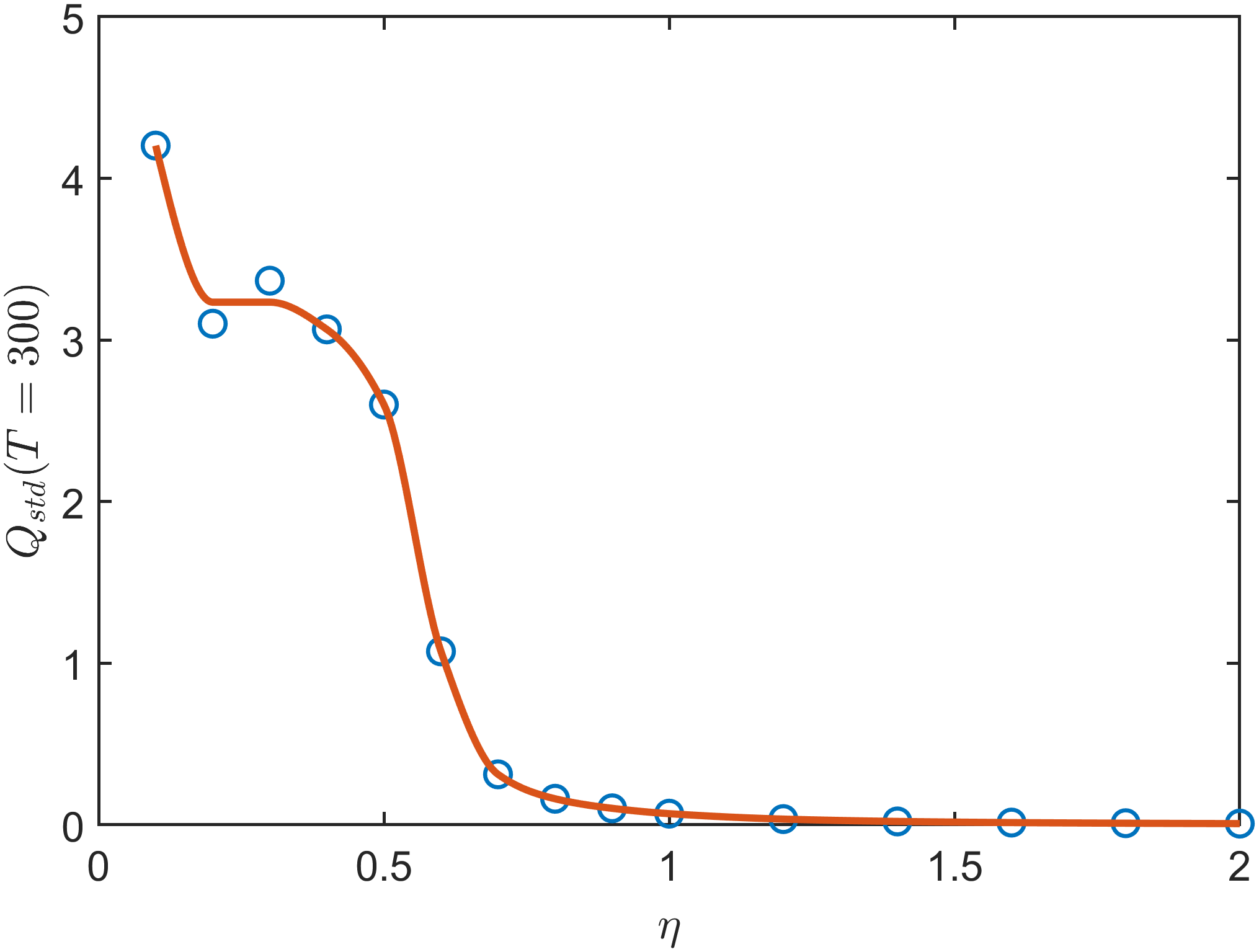}
	\includegraphics[width=0.32\textwidth]{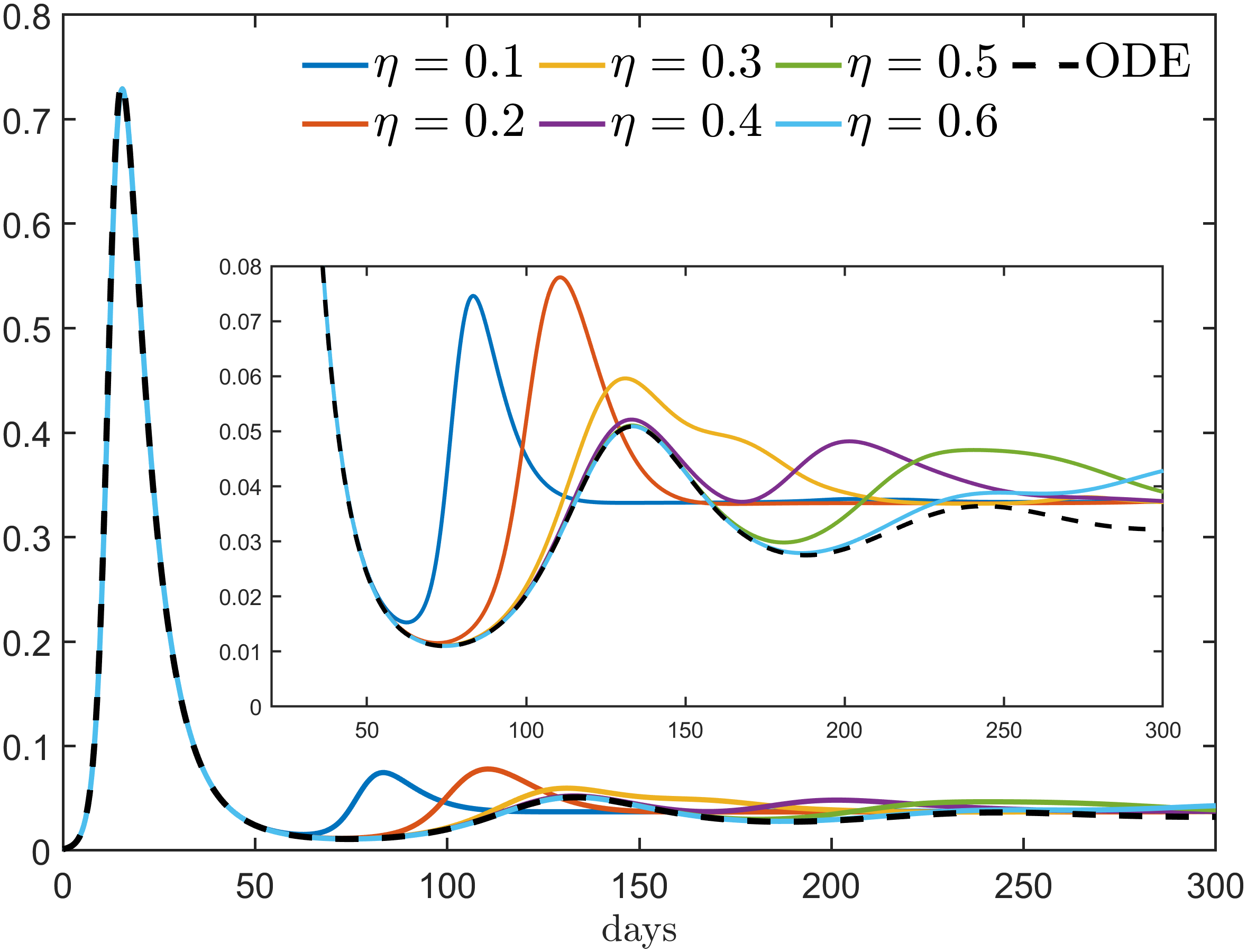}
	\includegraphics[width=0.32\textwidth]{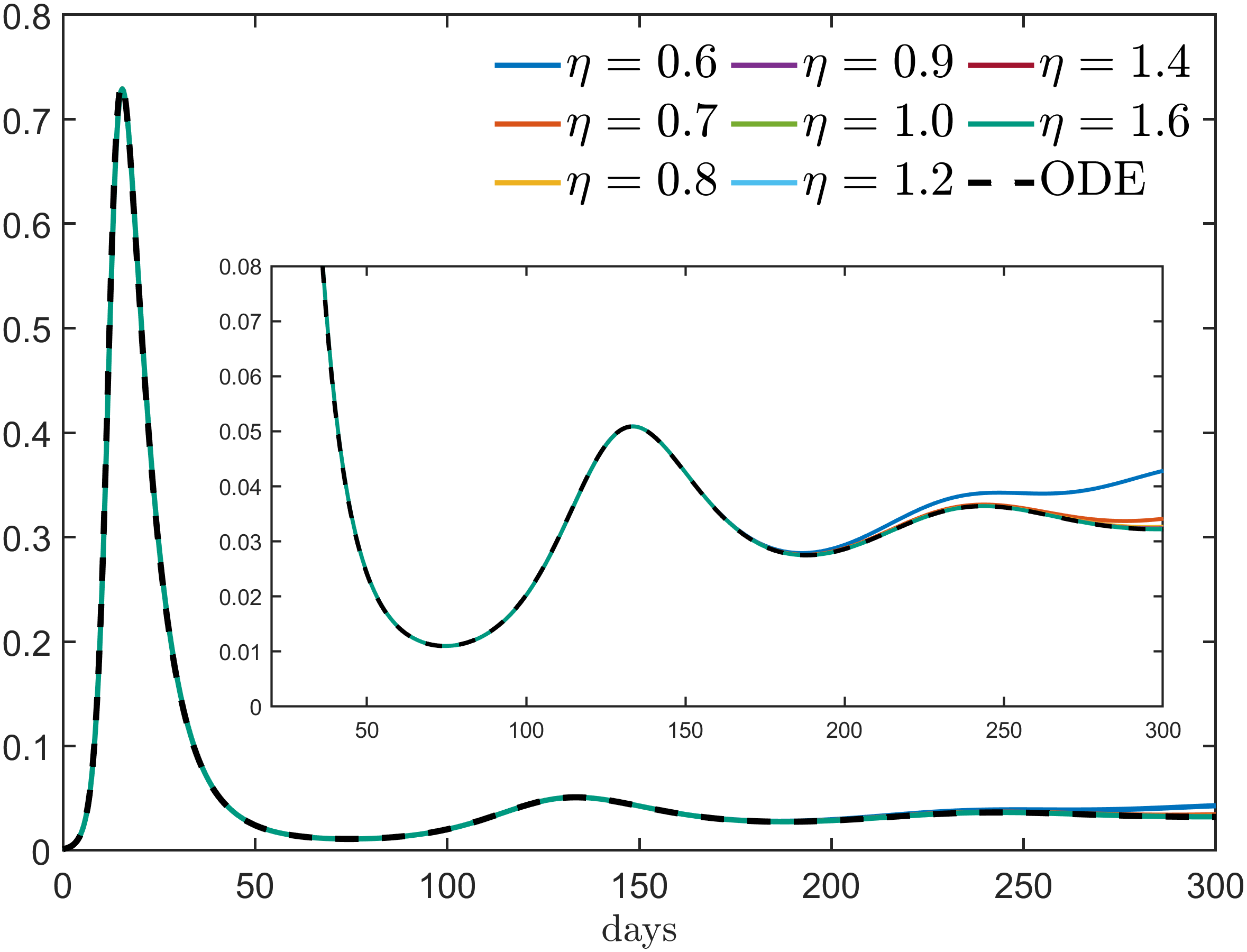}
	\vspace{-10pt}
	\caption{Dependence of the spatial order parameter and epidemic dynamics on
		the mixing parameter $\eta$ for $D=0.0005$ (third row of
		Figure~\ref{epidemic_phase_plot}).
		The spatial order parameter $Q_{std}$ (left) characterizes the aggregation effect. When $\eta$ is smaller, the aggregation of agents is more evident, so that it modifies both the amplitude and the arrival time of the second or third epidemic wave, as measured by the time evolution of virus carriers
		(i.e., the total population of $E$, $P$, $A$, $I^{+}$ and
		$I^{-}$) (middle and right).
	}
	\label{fig:order_parameter}
\end{figure}

The dependence of the spatial order parameter $Q_{std}$ at $T=300$ on
$\eta$ is presented in Figure~\ref{fig:order_parameter}.
The left panel shows that $Q_{std}$ decreases rapidly as
$\eta$ increases, with two-stage transition occurring
around $\eta \approx 0.2$ and $\eta\approx0.5$.
This behavior indicates the existence of three distinct
regimes separated by these thresholds.
To illustrate the dynamical consequences of this transition,
we divide the parameter range into three levels:
$\eta\leq0.6$, $ 0.6 \le \eta \le 1$ and $\eta \ge 1$.

\begin{enumerate}
\item[(1)] For $\eta\leq0.6$, the order parameter remains relatively large,
indicating strong spatial heterogeneity and the persistence
of infection hotspots.
In this regime, variations of $\eta$ mainly influence the
second epidemic wave, as illustrated in the middle panel.
Increasing $\eta$ weakens spatial aggregation and modifies
both the amplitude and the arrival time of the second peak. 

\item[(2)] For $0.6 \le \eta \le 1$, the spatial order parameter decreases and the aggregation only alters the third epidemic wave (see $\eta = 0.6$), with negligible influences on the second one.

\item[(3)] For $\eta \ge 1$, the spatial order parameter becomes very
small. In this regime, the epidemic dynamics becomes close to the
spatially homogeneous ODE limit.

\end{enumerate}

Together, these results indicate that the mixing parameter
$\eta$ controls a transition from a hotspot-dominated regime
to a nearly homogeneous  state, with multiple regimes
that affects different stages of infectious waves.

As a summary, incorporating aggregation behavior into epidemic
dynamics provides a quantitative framework for describing the
interaction between population mobility and infection processes.
The analysis highlights how spatial movement and aggregation can
shape the emergence and suppression of epidemic hotspots.
These findings offer useful insights into the role of spatial
heterogeneity in epidemic spreading and may help inform
spatially targeted intervention strategies.

\subsection{Effectiveness of the Lagrange multiplier correction}
\begin{figure}[h!]
	\centering
	\subfigure[Left: numerical blow-up. Middle: uncorrected phase diagram $(T = 80)$. Right: Corrected phase diagram $(T = 80)$.]
	{
		\includegraphics[width=0.3\textwidth]{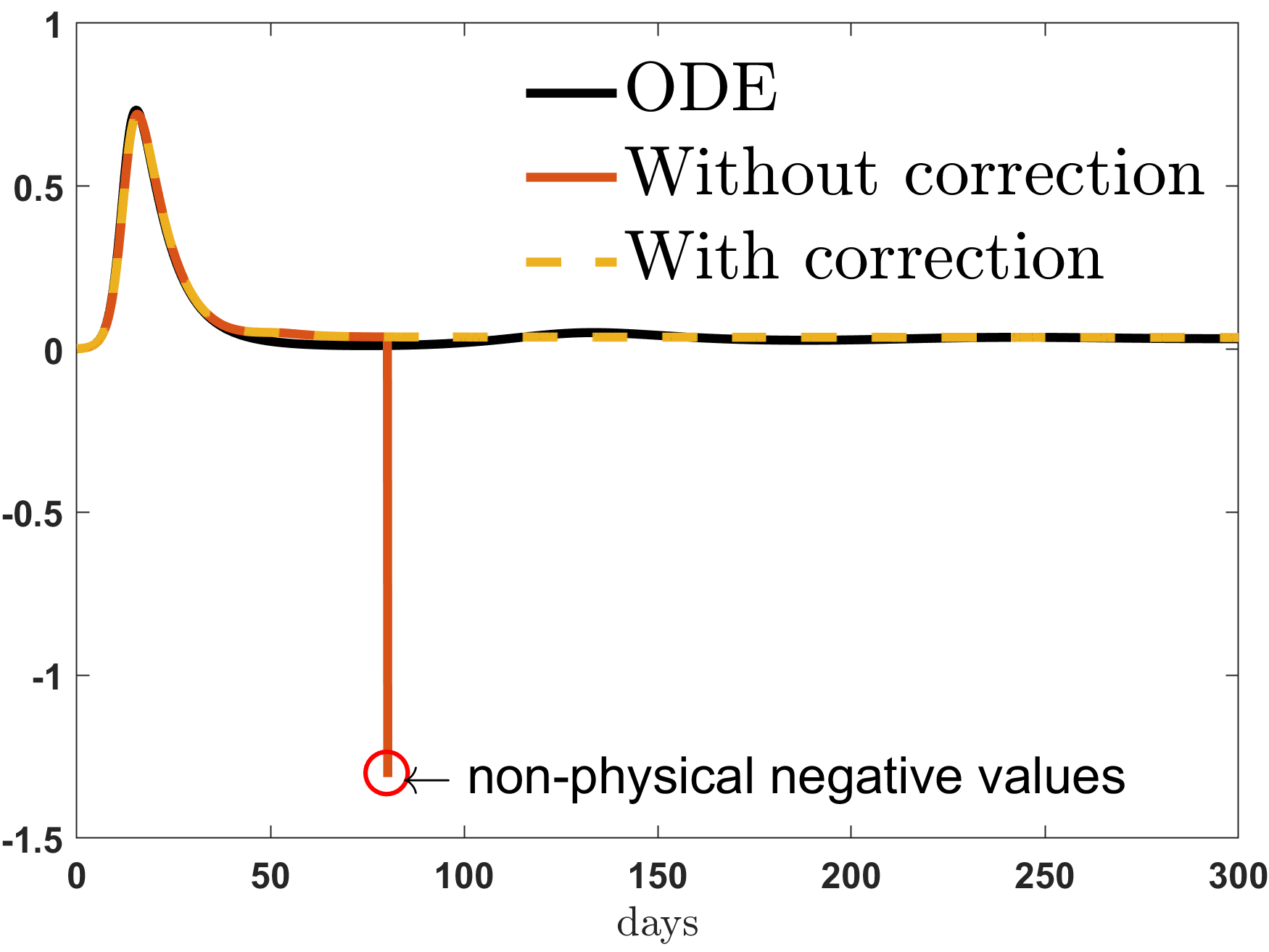} 
		\includegraphics[width=0.3\textwidth]{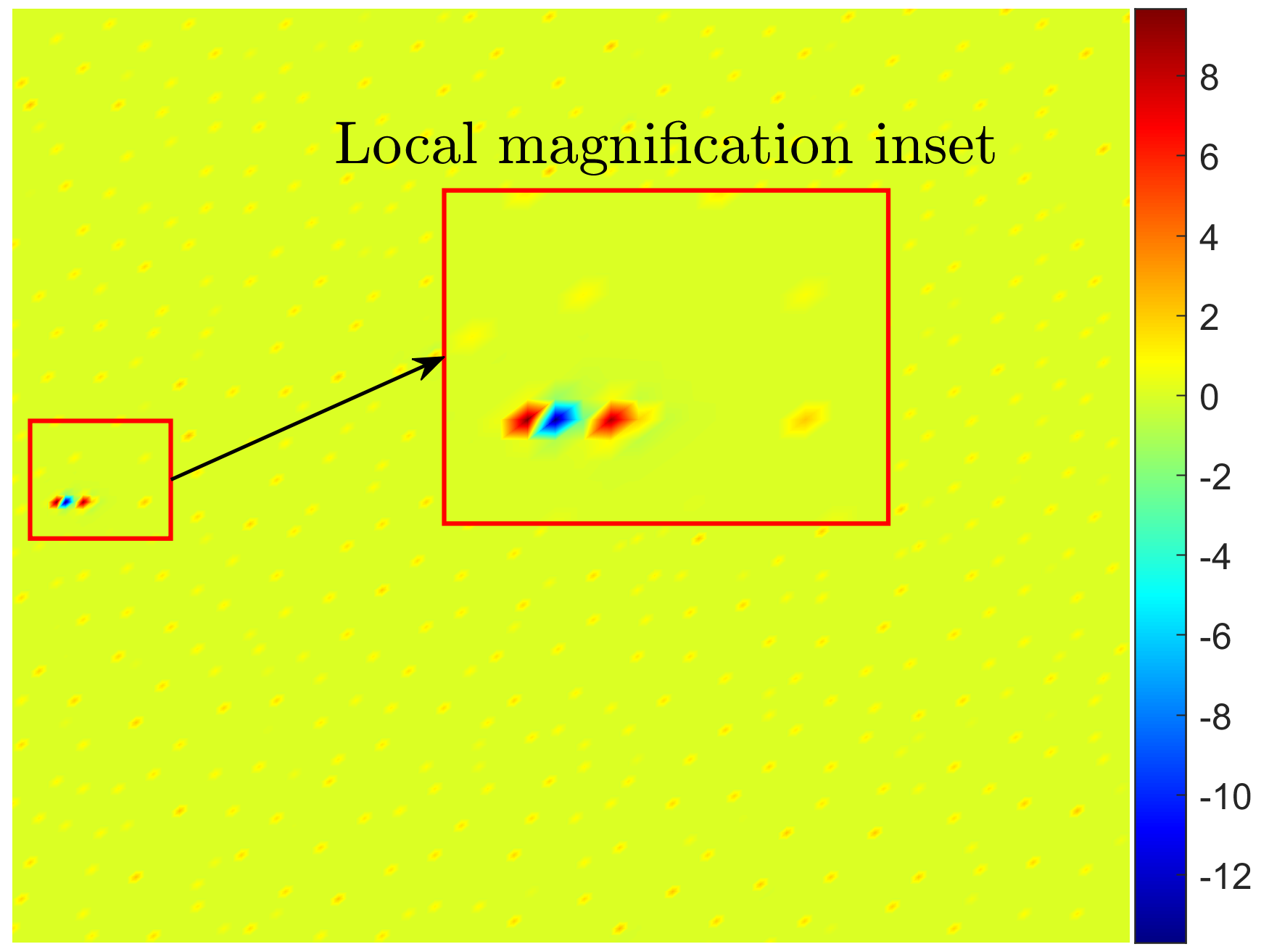} 
		\includegraphics[width=0.3\textwidth]{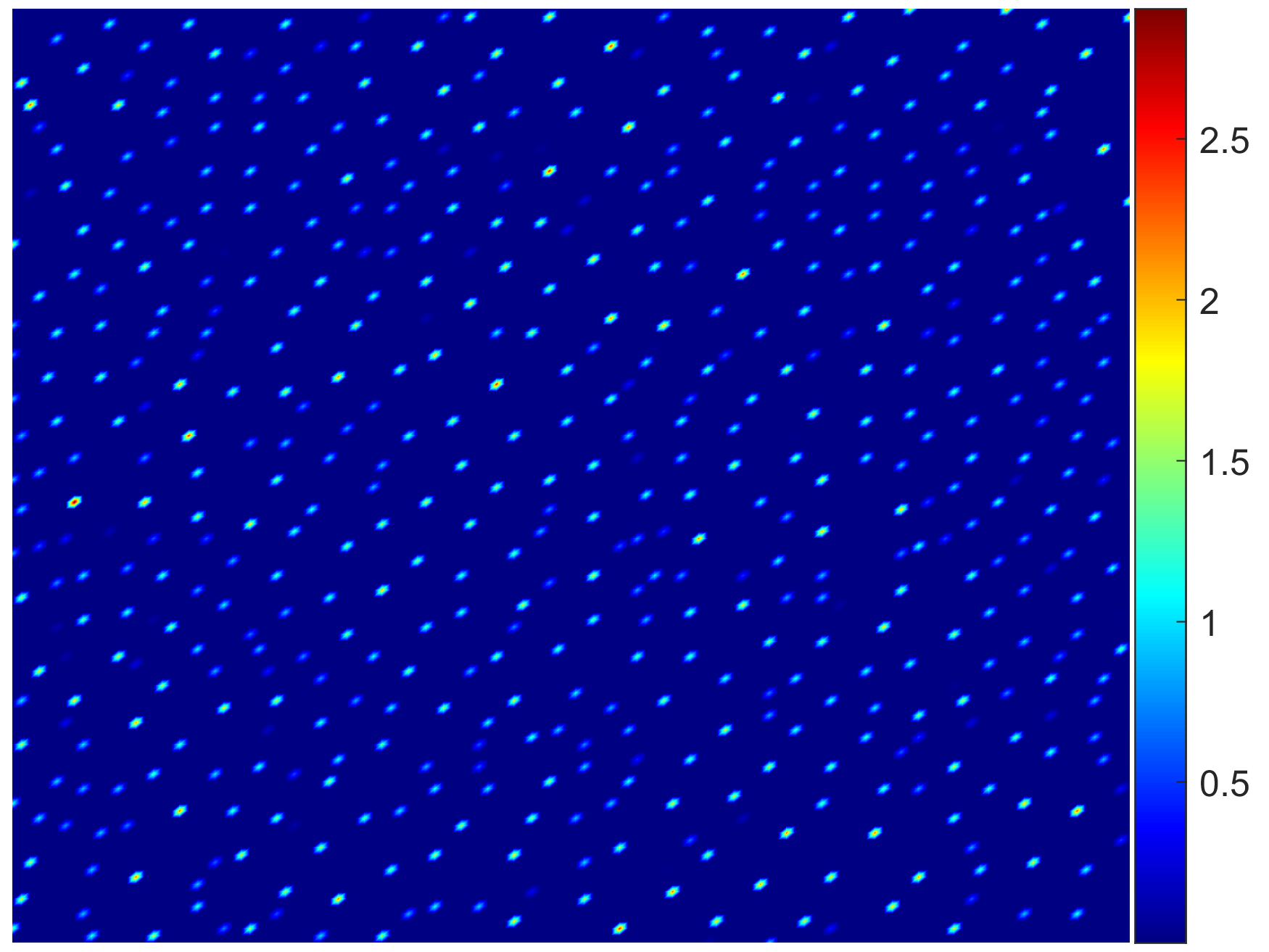}
	}
	\subfigure[Left: lower bounds of $\bm{\Psi}$. Middle:  the semi smooth Newton method with SPIMEX. Right: Conservation of discrete masses for $\bm{\Psi}$.]
	{
		\includegraphics[width=0.3\textwidth]{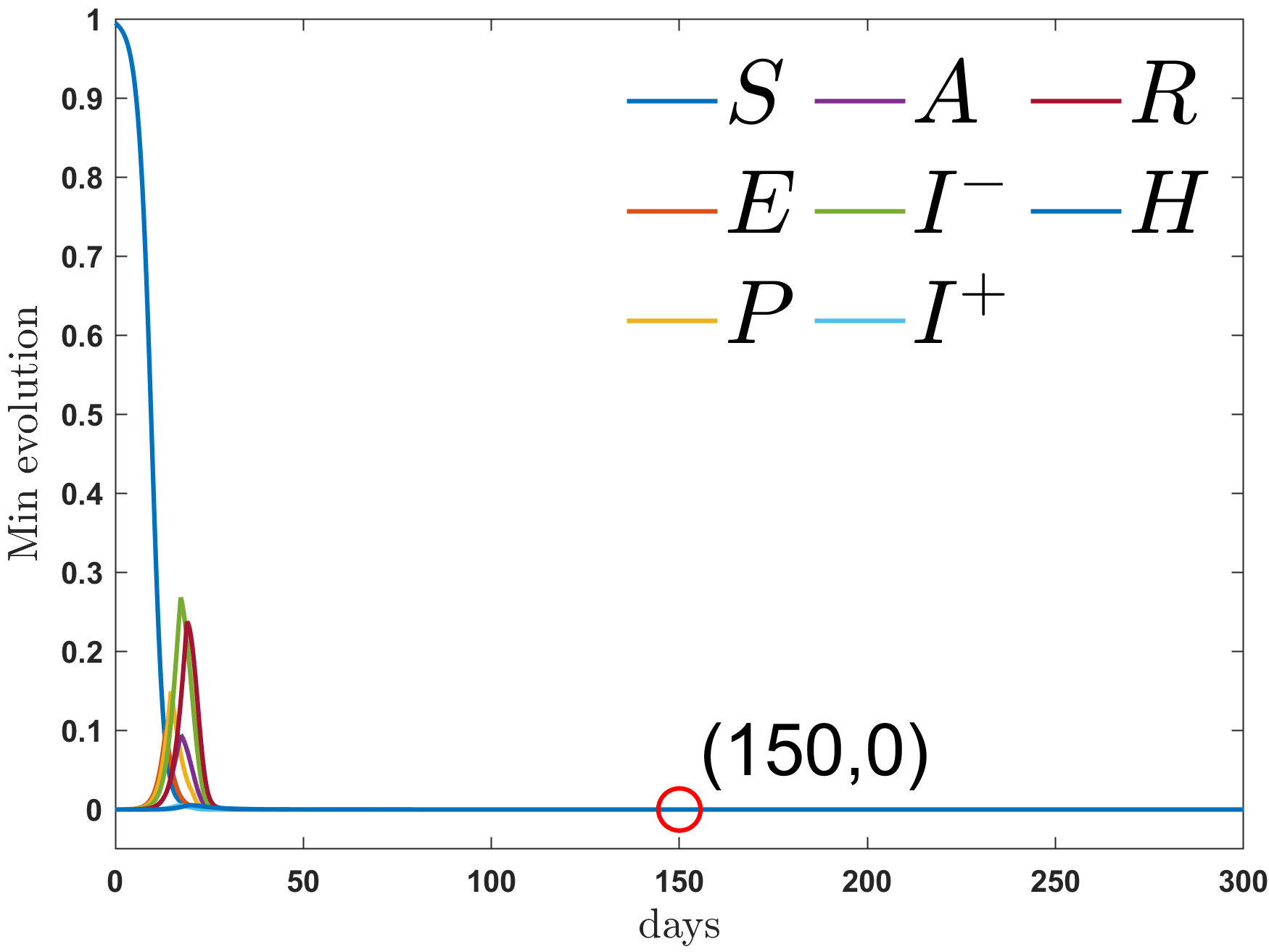}
		\includegraphics[width=0.3\textwidth]{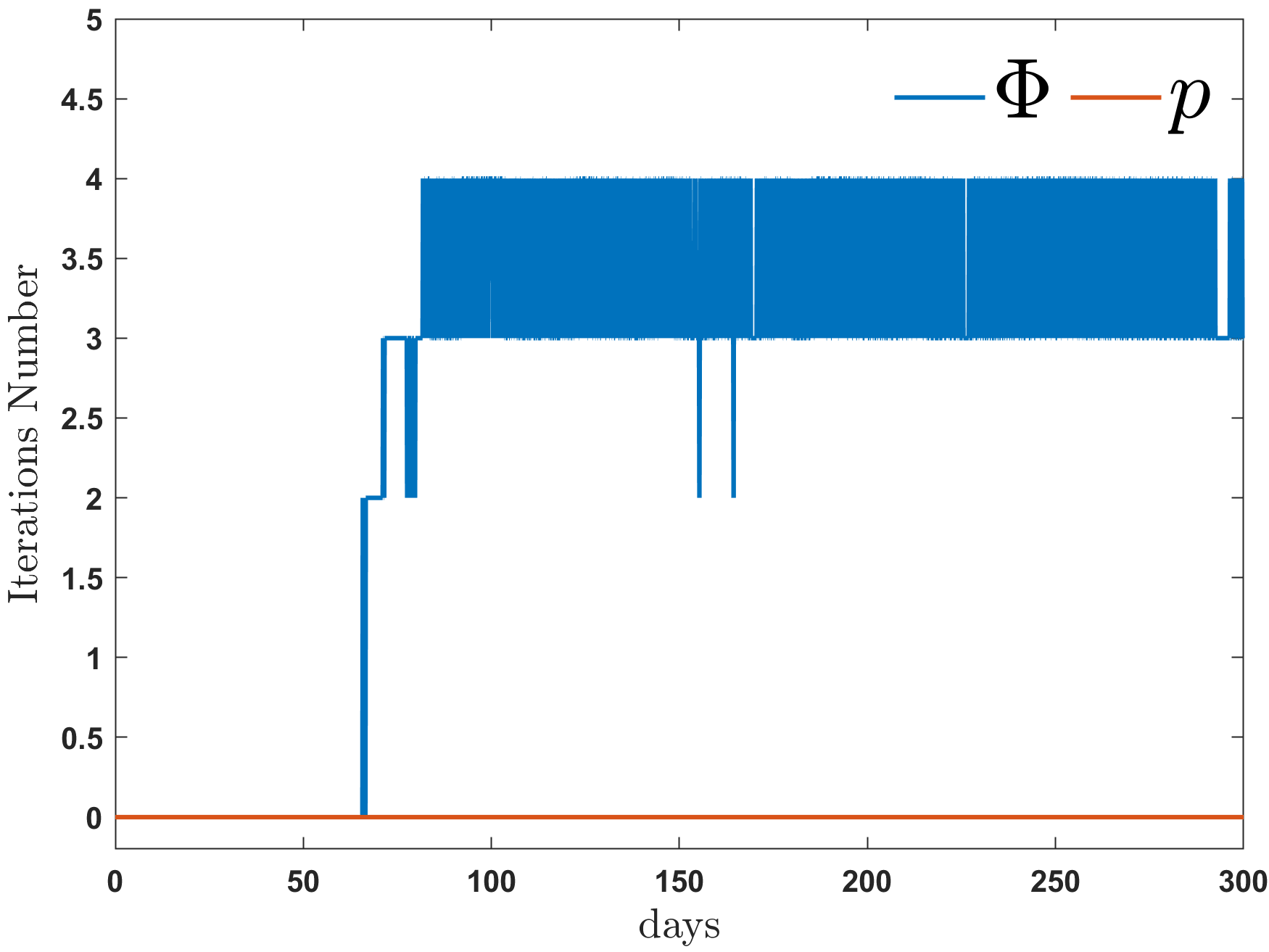} 
		\includegraphics[width=0.3\textwidth]{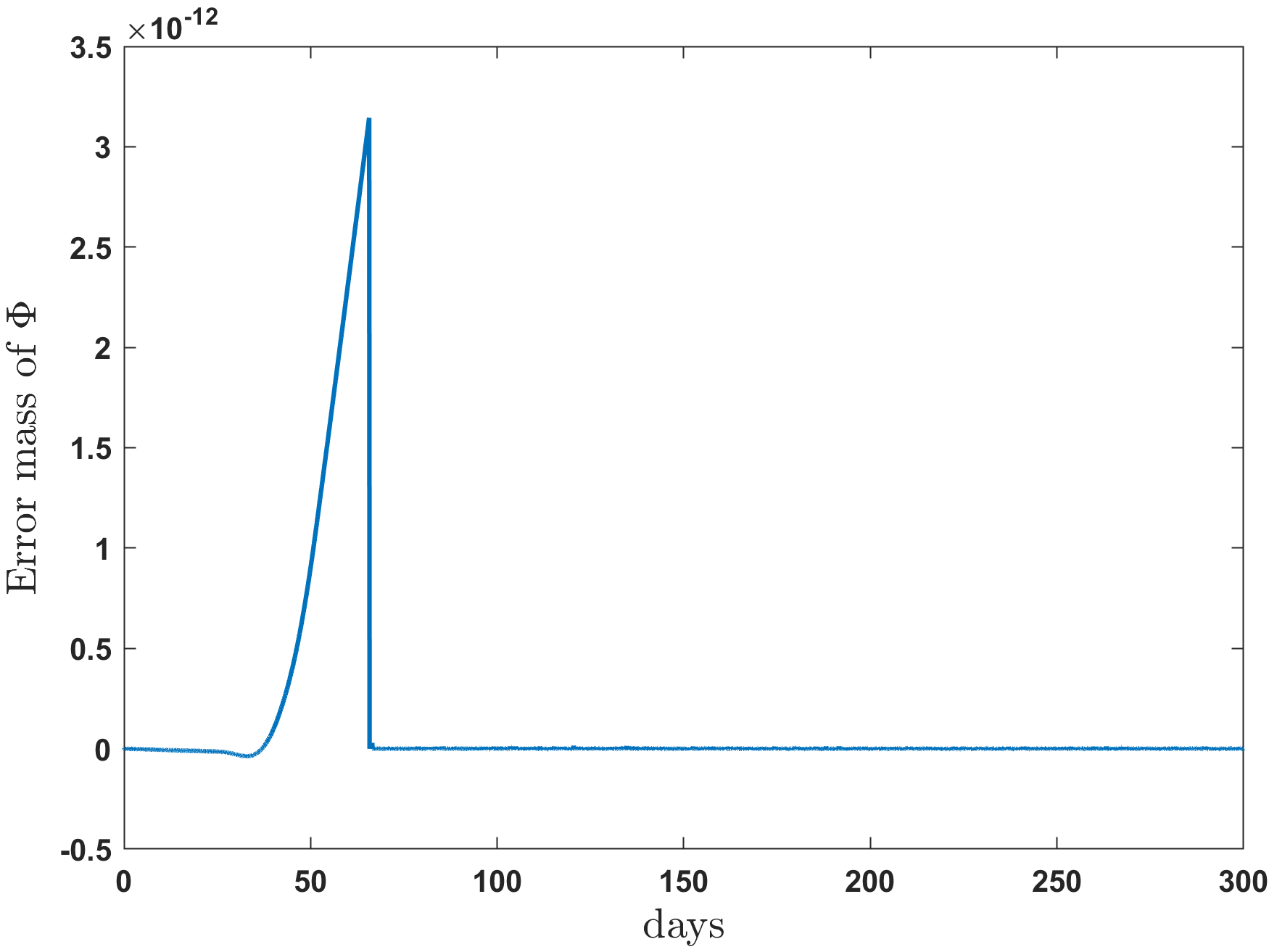}
	}
	\caption{Numerical simulations of an infectious-disease model before and after applying the Lagrange-multiplier correction are compared. In panel (a), the left subplot displays results from three schemes: the uncorrected method produces nonphysical negative values and becomes unstable, whereas the corrected method closely matches the reference ODE solution. The middle and right subplots show the state of the system at $T=80$ without or with corrections, respectively, demonstrating the instability caused by the uncorrected non-physical values. Panel (b) presents only corrected results: the left plot confirms the strict non-negativity of $\psi_i~(i=1,\dots,8)$; the middle plot depicts the number of iterations required by the semi-smooth Newton solver, illustrating its efficiency and robustness; and the right plot verifies discrete mass conservation of $\bm{\Psi}$.}
	\label{correction_plot}
\end{figure}

Finally, to evaluate the effectiveness of the Lagrange multiplier correction strategy, we revisit the epidemic model using the same initial conditions as in Example \ref{epid_simulation_plot}. The parameters are set as follows: diffusion coefficient $D = 0.00035$, $\eta = 0.03125$, time step $\Delta t = 0.1$, and final time $T = 300$, while all other model parameters remain unchanged. The corresponding results are shown in Figure~\ref{correction_plot}.

In the absence of correction, this parameter setting leads to the emergence of non-physical negative values in the solution, which rapidly destabilize the system and force the simulation to break down. Such instabilities are particularly severe in nonlinear coupled systems like the epidemic PDE system, where even small violations of positivity can be amplified in propagation.

Fortunately, when the Lagrange multiplier correction is applied,  it  avoids the numerical blow up and permits a stable long-time simulation. As observed in the previous section, the solution remains strictly non-negative, and the total mass is conserved up to machine accuracy. Notably, the system develops characteristic localized aggregation patterns: small hotspots gradually emerge and coalesce, mirroring the clustering behavior observed in real-world epidemics.

These results demonstrate that the Lagrange multiplier correction not only preserves physical constraints (positivity and mass conservation) but also faithfully captures the pattern-forming dynamics of the infectious disease model.

\section{Conclusion and discussion}\label{sec:conclusion}
In this paper, we propose an efficient structure-preserving implicit-explicit (SPIMEX) scheme with second-order accuracy  to simulate chemotaxis systems with singular sensitivity \eqref{cross_diffusion_PDE}, including both the crime model \eqref{crime_model} and the epidemic dynamics \eqref{Epidemic_PDE}. To overcome the numerical difficulties induced by singular sensitivity and to ensure the preservation of key physical properties, we utilize a posterior error estimation via the Lagrange multiplier correction technique, and derive rigorous error bounds. At the discrete level, an optimization-based projection of the intermediate solutions in the $L^2$-$H^1$ norm is used to enforce positivity and mass conservation. This projection involves only linearized operations and adds negligible computational overhead to the solution of the nonlinear algebraic system. Typical numerical experiments confirm our theoretical results and emphasize the necessity of the correction strategy. The proposed scheme allows us to verify the nucleation, spread, and dissipation of hotspots in crime modeling, as well as confirm that clustering of population density may aggravate the successive infectious wave in epidemic dynamics by accelerating virus transmission. 

The methodology and theoretical underpinnings in this work are readily extended to other dynamical systems featuring singular sensitivity, like opinion dynamics, which will be further investigated in our future work.

\section*{Acknowledgement}
This research was supported by the National Natural Science Foundation of China (NSF: \# 12231003, \# 11871105, \# 12571413). The authors are grateful to Prof. Chuntian Wang at Alabama Unversity and Prof. Yuan Zhang at Renmin University of China for fruitful discussions on epidemic modelings.

\appendix

\section{Epidemic ODE modeling and the corresponding chemotaxis PDEs}\label{app:ODE_PDE}
To study the impact of human's aggregation on the course of the epidemic, Reference  \cite{XiongWangZhang2024} has suggested to incorporate an environmental variable (representing the popularity of the location as perceived by nearby agents) into the above-mentioned agent-based symmetric random walk model, resulting in a biased random walk model.

The population density is divided into eight compartments: the susceptible agent $ {S}$, the exposed agent  ${E}$, the infectious and pre-symptomatic agent ${P}$, the asymptomatic agent $ {A}$, the mildly infectious symptomatic agent $I^-$, the infectious and symptomatic agent $ {I}^+$, the hospitalized agent $ {H}$ and the recovered agent $ {R}$.  The event types and parameter values in the epidemic model and their corresponding physical meanings are shown in Table \ref{tab:parameter1233_1}.
\begin{table}[htbp]
	\small
	\centering
	\caption{Event types, parameter values, sources and references, and corresponding physical meanings in the mathematical model for the epidemic ODE model.}
	\label{tab:parameter1233_1}
	\begin{tabular}{l c c}
		\toprule
		Event types & Parameters values & Physical meanings \\ 
		\midrule
		\multirow{3}{*}{Infectious contacts} 
		& $R_0 = 8.2$ & basic reproduction number \\
		& $\lambda = 1.018$ & rate of infection onset \\
		& $\beta = 1$ & reductive factor on infectivity of asymptomatic carriers \\
		\midrule
		End of latent period & $\alpha = {1}/{(1.2~\textup{days})}$ & inverse of latent period length \\
		\midrule
		\multirow{2}{*}{Symptom onset} 
		& $\eta' = {1}/{(1.8~\textup{days})}$ & rate of symptom onset \\
		& $\rho = 0.745$ & probability of symptomatic infectious cases \\
		\midrule
		\multirow{2}{*}{Hospitalization} 
		& $p_H = 0.0272$ & probability of hospitalization \\
		& $\delta_I^+ = {1}/{(3.8~\textup{days})}$ & rate of hospitalization onset \\
		\midrule
		\multirow{3}{*}{Recovery} 
		& $\delta_I^- = {1}/{(7.5~\textup{days})}$ & rate of virus removal \\
		& $\delta_A = {1}/{(7.5~\textup{days})}$ & rate of virus removal \\
		& $\delta_H = {1}/{(6~\textup{days})}$ & inverse of recovery period \\
		\midrule
		Immunity waning & $\delta_R = {1}/{(268~\textup{days})}$ & Inverse of immunity waning period \\
		\bottomrule
	\end{tabular}
\end{table}

The epidemic ODE system reads that
\begin{equation*}
	\left\{
	\begin{split}
		&\partial_t  {S} =  
		-   { \lambda } ( \beta (    {P}+  {A} ) +  I^-+  I^+ )  {S  } + \delta _{R}   {R}, \\
		&\partial_t  {E}=       \lambda   ( \beta (  {P}  +  {A } ) +    I^-+  I^+ ) { S     }  -\alpha  {E}  ,\\
		&\partial_t  {P} =    \alpha   {E} -\eta^\prime  {P }, \\
		&\partial_t  {A} =    \eta ^\prime (1- \rho)  {P}  - \delta_A   {A }, \\
		&\partial_t  I^-=       \eta ^\prime \rho (1- p_H) { P }-\delta_I ^{-}  I^-, \\
		&\partial_t  I^+ =      \eta ^\prime \rho  p_H  { P }-\delta_I^{+}   I^+, \\
		&\partial_t  H =    \delta_I^{+}   I^+ - \delta _{H}  {H}, \\
		&\partial_t  R =      \delta _A   {A} + \delta_I^{-}   I^-  + \delta_H  { H} - \delta _{R}  {R }.
	\end{split}
	\right.
\end{equation*} 


The event types, parameter values and their corresponding physical meanings of the chemotaxis terms are shown in Table \ref{biase}.
\begin{table}[htbp]
	\centering
	\caption{
		Event types,  parameter values and corresponding physical meanings for the biased random walk  model.}
	\begin{tabular}{lcc}
		\toprule
		Event types  & Parameters    values   & Physical meanings \\ 
		\midrule
		\multirow{1}{*}{Increase of attractiveness}   
		&      {$ \bar{ \delta }  ^{+}_{\mathcal P} = 0.3$}     &      { size of increment  }  \\   
		\midrule
		\multirow{2}{*}{Spread of attractiveness}   
		&       {$\eta $}      &      {magnitude of spread }  \\
		&       {$\eta D$}   &      diffusion rate of attractiveness level  \\
		\midrule
		\multirow{2}{*}{Decrease of attractiveness}  
		&         $\bar{\delta}^{-}_{{\mathcal P}} =0.36$ & size of decay \\
		&        { $\mathcal P_{\min} =1/30$} & minimum value of attractiveness level\\
		\bottomrule
	\end{tabular}
	\label{biase}
\end{table}

\section{Calculation of Lagrange multipliers in epidemic model}\label{app:Lagrange_method}

In this section, we briefly describe the process of calcuating the Lagrange multipliers $\bm{\lambda}_h^{k+1}(\mathbf{x})\in X^8$, $\zeta_h^{k+1}\in X$ and $\xi^{k+1}\in\mathbb{R}$ in epidemic model (see \eqref{KKT_8}). 

\subsection{$L^2$ projection of ${\Psi}_h^{k+1}$} From the complementary relaxation conditions \eqref{KKT_8} in the main text, numerical solutions $\bm{\Psi}_h^{k+1}$ are expressed as
\begin{equation}\label{cut_off}
	(\psi_{k, h}^{k+1}, \lambda_{k, h}^{k+1})
	=\left\{
	\begin{split}
		& (\tilde \psi_{k, h}^{k+1} - {\xi}^{k+1},  0), \quad &\tilde \psi_{k, h}^{k+1} - {\xi}^{k+1} > 0,\\
		& (0, -(\tilde \psi_{k, h}^{k+1} - {\xi}^{k+1})), \quad &\tilde \psi_{k, h}^{k+1} - {\xi}^{k+1} \le 0,
	\end{split}
	\right. \quad k=1,\dots,8,
\end{equation}
where the Lagrange multiplier $\xi^{k+1}$ is determined by the mass constraints
\begin{equation*}
	h^2\sum_{k=1}^{8} \sum_{\tilde{\psi}_{k, h}^{k+1} - \xi^{k+1} > 0} (\tilde{\psi}_{k, h}^{k+1} - \xi^{k+1}) =  \langle \sum_{i=1}^{8}\psi_{k,h}^{0},1\rangle.
\end{equation*}
Hence $\xi^{k+1}\in\mathbb{R}$ is the solution to the nonlinear algebraic equation
\begin{equation}\label{nonlinear_equation}
	F(\xi^{k+1}) = h^2\sum_{k=1}^{8} \sum_{i,j=1}^{N} (\tilde{\psi}_{k, h}^{k+1} - \xi^{k+1})^+ - \langle \sum_{i=k}^{8}\psi_{k,h}^{0},1\rangle=0
\end{equation}
where $f^+ = \max{f,0}$ for $f\in \mathbb{R}$. Semi-smooth Newton methods \cite{FacchineiPang2003,KojimaShindo1986} can be employed to solve \eqref{nonlinear_equation}, i.e. for some initial guess $\xi_0\in\mathbb{R}$, find the root of $F(\xi) = 0$ by updating
\begin{equation}\label{solve_xi}
	\xi_{s+1} = \xi_s - V^{-1}(\xi_s) F(\xi_s),\quad s = 0,1,\dots,
\end{equation}
where $V(\xi_s)$ is a generalized derivate in semi-smooth Newtown methods as
\begin{equation*}
	V(\xi_s) = -h^2\sum_{k=1}^{8}\sum_{i,j=1}^{N}\text{sgn} ((\tilde{\psi}_{k,h}(\mathbf{x}_{i,j})-\xi_s)^+)
\end{equation*}
and $\text{sgn}(\cdot)$ is the sign function with $\text{sgn}(s) = 1 (s > 0)$, $\text{sgn}(0) = 0$ and $\text{sgn}(s) = -1(s < 0)$. Noticing $\xi^{k+1}$ is supposed to be of small magnitude \cite{cheng2022bound,TongFenghua2024Positivity} and $\xi^{k+1}\geq 0$ (Lemma
4.1), we can choose $\xi_0 = 0$ to start the semi-smooth Newton iterations. Once $\xi^{k+1}$
is known, we can update $(\psi_{k,h}^{k+1},\lambda_{k,h}^{k+1})$ according to \eqref{cut_off}. In all our numerical experiments, \eqref{solve_xi} converges in only
one iteration so that the cost of solving \eqref{nonlinear_equation} is negligible. The secant method can be also used to solve \eqref{nonlinear_equation} \cite{cheng2022bound,TongFenghua2024Positivity}.

\subsection{$H^1$ projection of $p_h^{k+1}$}
A semi-smooth Newton method is used to solve the Lagrange multiplier $\zeta_{h}^{k+1}(\mathbf{x})\in X$ under $H^{1}$ projection with only positivity constraint. The KKT condition for $\zeta_{h}^{k+1}(\mathbf{x})\in X$ can be transformed to an equivalent system as follows
\begin{equation}\label{F_p}
	\begin{split}
		& F(U;\mathbf{x}) = -\Delta_{h}U^{+}(\mathbf{x}) + U(\mathbf{x}) - (I-\Delta_{h})\tilde{p}_{h}^{k+1}(\mathbf{x}) = 0, \quad \mathbf{x}\in\Omega_{h},
	\end{split}
\end{equation}
where $U(\mathbf{x})=U^{+}(\mathbf{x})-U^{-}(\mathbf{x})$, $U^{+}=\max(U^{k+1}(\mathbf{x}),0)$, $U^{-}=\max(-U^{k+1}(\mathbf{x}),0)$, and $U^{+}(\mathbf{x})=p_{h}^{k+1}(\mathbf{x})$, $U^{-}(\mathbf{x})=\zeta_{h}^{k+1}(\mathbf{x})$.

The semi-smooth Newton method can be applied to solve \eqref{F_p}. We start with the generalized Jacobian to be used in the semi-smooth Newton iterations. At $U(\mathbf{x})\in X$, the generalized Jacobian $J$ (acting on a vector $V(\mathbf{x})\in X$) can be written as:
\begin{equation}\label{J_V}
	J(V(\mathbf{x})) = I - \Delta_{h}\left(\operatorname{sgn}(U^{+}(\mathbf{x}))V(\mathbf{x})\right).
\end{equation}
By careful computation based on \eqref{F_p}-\eqref{J_V}, a semi-smooth Newton method for solving \eqref{F_p} with initial $U_{0}(\mathbf{x})$ can be designed as follows
\begin{equation*}
	\begin{split}
		& U_{s+1} = U_{s} + V_{s},
	\end{split}
\end{equation*}
where $V_{s}(\mathbf{x}) = (I - \Delta_{h}\operatorname{sgn}(U_{s}^{+}))^{-1}(-F(U_{s}))$. For solving $V_{s}$, we can use $(I-\Delta_{h})^{-1}$ as a pre-conditioner. When Newton iteration converges for some $s\geq 0$, we have $p^{k+1}(\mathbf{x})=U_{s+1}^{+}(\mathbf{x})$. To start the semi-smooth Newton iteration, we may set $U_{0}(\mathbf{x})=\tilde{p}_{h}^{k+1}(\mathbf{x})$.






\end{document}